\documentclass[10pt]{article}
\font\smallit=cmti10

\usepackage{amssymb,latexsym,amsmath,epsfig,amsthm,graphicx,tikz-network,array,mathbbol,subcaption,verbatim,enumitem} 

\makeatletter

\renewcommand\section{\@startsection{section}{1}{\z@}
{-30pt \@plus -1ex \@minus -.2ex}
{2.3ex \@plus.2ex}
{\normalfont\normalsize\bfseries}}

\renewcommand\subsection{\@startsection{subsection}{2}{\z@}
{-3.25ex\@plus -1ex \@minus -.2ex}
{1.5ex \@plus .2ex}
{\normalfont\normalsize\bfseries}}

\renewcommand\subsubsection{\@startsection{subsubsection}{2}{\z@}
{-3.25ex\@plus -1ex \@minus -.2ex}
{1.5ex \@plus .2ex}
{\normalfont\normalsize\bfseries}}

\renewcommand{\@seccntformat}[1]{\csname the#1\endcsname. }

\makeatother

\begin{document}

\begin{center}
\uppercase{\bf Optimally reconnecting graphs against an edge-destroying adversary}
\vskip 20pt
{\bf Daniel C. McDonald}\\
{\smallit Wolfram Research Inc., Champaign, Illinois}\\
{\tt daniel.cooper.mcdonald@gmail.com}\\
\end{center}
\vskip 30pt

\centerline{\bf Abstract}
\noindent
We introduce a model involving two adversaries Buster and Fixer taking turns modifying a connected graph, where each round consists of Buster deleting a subset of edges and Fixer responding by adding edges from a finite reserve set of weighted edges to leave the graph connected, with Buster limited by the total number of edges he is allowed to delete throughout the game. Fixer wins if she can reconnect the graph after Buster has reached his limit of edges to delete, while Buster wins if he can delete edges in such a way that Fixer cannot reconnect the graph using the remaining edges in reserve. With the weights representing the cost for Fixer to use specific reserve edges to reconnect the graph, we prove that a greedy strategy for Fixer always results in an optimal result for Fixer: victory, if possible, for as cheaply as can be guaranteed against any Buster strategy, and if defeat cannot be avoided, the cheapest possible loss that can be guaranteed against any Buster strategy.

\pagestyle{myheadings} 
\thispagestyle{empty} 
\baselineskip=12.875pt 
\vskip 30pt

\newtheorem{thm}{Theorem}[section]
\newtheorem{lem}[thm]{Lemma}
\newtheorem{prop}[thm]{Proposition}
\newtheorem{cor}[thm]{Corollary}
\newtheorem{hyp}[thm]{Hypothesis}
\newtheorem{case}[thm]{Case}
\theoremstyle{definition}
\newtheorem{sce}{Scenario}
\newtheorem*{ack}{Acknowledgment}

\newcommand{\mc}[1]{\ensuremath{\mathbb{#1}}}
\newcommand{\mcp}[2]{\ensuremath{\mathbb{#1}[{#2}]}}
\newcommand{\his}[2]{\ensuremath{\mathbb{H}_{#2}^{\mathbb{#1}}}}
\newcommand{\ess}[2]{\ensuremath{#1_{#2}}}
\newcommand{\es}[3]{\ensuremath{#1^{#2}_{#3}}}
\newcommand{\pair}[4]{\ensuremath{\langle #1,#2 \mid #3,#4\rangle}}
\newcommand{\pay}[1]{\ensuremath{\pi(#1)}}

\section{Introduction} 
\label{intro}
Suppose a network must stay connected in the face of some adversary that periodically destroys subsets of its edges, with the network being reconnected after each attack by adding replacement edges, each having its own cost. Beyond the requirement that each individual selection of replacement edges must reconnect the network, ideally taken together these selections should keep the network connected for as long as possible as cheaply as possible. In Subsection \ref{model} we introduce a graph-theoretic framework to formally model this situation, and in Subsection \ref{example} we give an example of the model in action. In Subsection \ref{statement} we state our main theorem, that a greedy strategy for selecting replacement edges is optimal. In Subsection \ref{past} we relate greedy strategies in our model to finding minimum spanning trees of graphs, as well as compare and contrast our model to the well-studied Maker-Breaker game, before giving an overview of several related models. We outline some directions for future research in Subsection \ref{future}. Using tools developed in Section \ref{prelim}, we prove our main theorem in Section \ref{proof}. We extract out the details of certain technical proofs to appendices at the end of the paper.

\subsection{Description of Model}
\label{model}

For this paper, all vertex sets of multigraphs will be unambiguous, so we view multigraphs simply as multisets of edges; in particular, for any multigraph $G$ we shall let $|G|$ denote the number of edges of $G$, rather than the number of vertices. To initialize an instance of our model, we are given a finite multigraph $G$, a finite multiset $R$ of weighted "reserve" edges between vertices of $G$ satisfying $G\cap R=\emptyset$ (we consider all edges distinct, even if they have the same endpoints), and a nonnegative integer $\ell$ called the \emph{limit}. Each edge $r\in R$ has some nonnegative, finite weight $w(r)$; for $R'\subseteq R$, define $w(R')=\sum _{r\in R'}w(r)$.

There are two parties in this model: a positive actor Fixer and an antagonist Buster, who play against each other in a series of rounds. The $k$th round starts with a multigraph \ess{G}{k} and multiset \ess{R}{k} of ``reserve'' edges, with initializations of $\ess{G}{1}=G$ and $\ess{R}{1}=R$. For $1\leq j<k$, the $j$th round consisted of Buster ``busting'' the graph \ess{G}{j} by removing a set \ess{B}{j} of its edges, under the limitation that the total number of edges he removes throughout the game cannot exceed the game's limit $\ell$, followed by Fixer ``fixing'' the graph by reconnecting it with a set \ess{F}{j} of edges taken from the reserve set \ess{R}{j}. At the beginning of the $k$th round, if $\sum _{j=1}^{k-1}|\ess{B}{j}|=\ell$ (i.e. Buster cannot remove a nonempty set of edges in the $k$th round without violating his upper limit for total edges removed), then say \emph{Fixer wins in the $k$th round}. If Fixer does not win in the $k$th round, then Buster removes a nonempty set $\ess{B}{k}\subseteq\ess{G}{k}$ of edges from the connected graph \ess{G}{k} satisfying the limiting inequality $\sum _{j=1}^{k}|\ess{B}{j}|\leq\ell$. If $(\ess{G}{k}-\ess{B}{k})\cup\ess{R}{k}$ is disconnected, then Fixer cannot reconnect the graph, and we say \emph{Buster wins in the $k$th round}. If $(\ess{G}{k}-\ess{B}{k})\cup\ess{R}{k}$ is connected, then Fixer adds to the graph $\ess{G}{k}-\ess{B}{k}$ a set $\ess{F}{k}\subseteq\ess{R}{k}$ of edges in such a way that the graph $(\ess{G}{k}-\ess{B}{k})\cup\ess{F}{k}$ is connected, in which case we set $\ess{G}{k+1}=(\ess{G}{k}-\ess{B}{k})\cup\ess{F}{k}$ and $\ess{R}{k+1}=\ess{R}{k}-\ess{F}{k}$ before continuing on to the $(k+1)$st round. Since our original graph $G$, reserve edge set $R$, and limit $\ell$ were finite, the game ends eventually, with its payout calculated as follows: if Fixer wins in the $k$th round, the payout for Fixer is $1+w(\ess{R}{k})$ and the payout for Buster is $-1-w(\ess{R}{k})$, and if Buster wins in the $k$th round, the payout for Fixer is $-1-w(\ess{R}{1}-\ess{R}{k})$ and the payout for Buster is $1+w(\ess{R}{1}-\ess{R}{k})$. That is, the game is zero-sum, with Fixer's winning share equal to one more than the total cost of the reserve edges Fixer did not have to activate, and Buster's winning share equal to one more than the total cost of the edges Fixer did have to activate; the one is added/subtracted so that winning payouts are always strictly positive and losing payouts are always strictly negative. Since $w(\ess{R}{k})=w(\ess{R}{1})-\sum _{j=1}^{k-1} w(\ess{F}{j})$ and $w(\ess{R}{1}-\ess{R}{k})=\sum _{j=1}^{k-1} w(\ess{F}{j})$, we see that the primary goal for Fixer is to win the game, with the secondary goal of spending as little as possible on reserve edges no matter the final result, and that the primary goal for Buster is to win the game, with the secondary goal of forcing Fixer to spend as much as possible on reserve edges no matter the final result.

Let a \emph{history} be a (finite) gameplay sequence of the form $((\ess{G}{1},\ess{R}{1}),\ess{B}{1},\ess{F}{1},(\ess{G}{2},\ess{R}{2}),\ess{B}{2},\ess{F}{2},\ldots)$, implicitly endowed with the limit $\ell$ of the game, which we also refer to as the limit of the history. Let \his{B}{} denote the set of histories ending in a move by Buster, with \his{B}{k} being the subset of those ending in a move by Buster in the $k$th round, and let \his{F}{} denote the set of histories ending in a move by Fixer, with \his{F}{k} being the subset of those ending in a move by Fixer in the $k$th round. Let the graph, reserve set, Buster move, and Fixer move of the $k$th round of a history \mc{h} be denoted by, respectively, \es{G}{\mc{h}}{k}, \es{R}{\mc{h}}{k}, \es{B}{\mc{h}}{k}, and \es{F}{\mc{h}}{k}. Call \mc{h} \emph{complete} if Fixer or Buster wins \mc{h}; note that if \mc{h} is complete and ends with Fixer winning in the $k$th round, then $\mc{h}\in\his{F}{k-1}$ and the final entry in \mc{h} is $(\ess{G}{k},\ess{R}{k})$, while if \mc{h} is complete and ends with Buster winning in the $k$th round, then $\mc{h}\in\his{B}{k}$ and the final entry in \mc{h} is $\ess{B}{k}$. For a complete history \mc{h}, let $|\mc{h}|=k$ denote that Fixer or Buster wins in the $k$th round. For a complete history \mc{h}, define the \emph{Fixer payout function} $\pay{\mc{h}}$ to equal the Fixer payout defined in the previous paragraph, i.e. $\pay{\mc{h}}=1+w(\es{R}{\mc{h}}{|\mc{h}|})$ if Fixer wins and $\pay{\mc{h}}=-1-w(\es{R}{\mc{h}}{1}-\es{R}{\mc{h}}{|\mc{h}|})$ if Buster wins (since the game is zero-sum, the Buster payout of \mc{h} is always equal to $-\pay{\mc{h}}$). For complete histories \mc{h} and \mc{h'} having the same limit as well as satisfying $\es{G}{h}{1}=\es{G}{h'}{1}$ and $\es{R}{h}{1}=\es{R}{h'}{1}$, say \mc{h} is \emph{Fixer-superior} to \mc{h'} if $\pay{\mc{h}}\geq\pay{\mc{h'}}$.

A \emph{Buster strategy} is a function $\beta$ from \his{F}{} to a valid Buster action. That is, if $\mc{h}=((\ess{G}{1},\ess{R}{1}),\ess{B}{1},\ess{F}{1},\ldots,(\ess{G}{k},\ess{R}{k}))\in\his{F}{k-1}$ and the limit $\ell$ of the game satisfies $\sum _{j=1}^{k-1}|\ess{B}{j}|<\ell$ (i.e. Fixer does not win in the $k$th round), then $\beta(\mc{h})=\ess{B}{k}\subseteq\ess{G}{k}$ with $1\leq |\ess{B}{k}|\leq \ell-\sum _{j=1}^{k-1}|\ess{B}{j}|$, so that $((\ess{G}{1},\ess{R}{1}),\ess{B}{1},\ess{F}{1},\ldots,(\ess{G}{k},\ess{R}{k}),\ess{B}{k})\in\his{B}{k}$. A \emph{Fixer strategy} is a function $\phi$ from \his{B}{} to a valid Fixer action. That is, if $\mc{h}=((\ess{G}{1},\ess{R}{1}),\ess{B}{1},\ess{F}{1},\ldots,(\ess{G}{k},\ess{R}{k}),\ess{B}{k})\in\his{B}{k}$ and $(\ess{G}{k}-\ess{B}{k})\cup\ess{R}{k}$ is connected (i.e. Buster does not win in the $k$th round), then $\phi(\mc{h})=\ess{F}{k}\subseteq\ess{R}{k}$ with $(\ess{G}{k}-\ess{B}{k})\cup\ess{F}{k}$ connected, so that for $\ess{G}{k+1}=(\ess{G}{k}-\ess{B}{k})\cup\ess{F}{k}$ and $\ess{R}{k+1}=\ess{R}{k}-\ess{F}{k}$ we have $((\ess{G}{1},\ess{R}{1}),\ess{B}{1},\ess{F}{1},\ldots,(\ess{G}{k},\ess{R}{k}),\ess{B}{k},\ess{F}{k},(\ess{G}{k+1},\ess{R}{k+1}))\in\his{F}{k}$.

Let an \emph{extension} of a history \mc{h} be a history \mc{h'} having \mc{h} as a prefix and inheriting the limit $\ell$ of \mc{h}; say \mc{h'} is a \emph{proper extension} of \mc{h} if \mc{h'} is an extension of \mc{h} but does not equal \mc{h}. Given a history \mc{h} and a sequence $\sigma$ of moves, we slightly abuse notation by letting $(\mc{h},\sigma)$ denote the extension of \mc{h} obtained by appending $\sigma$ to the end of \mc{h}. If $\mc{h}\in\his{F}{}$, define $\mcp{h}{\beta}=(\mc{h},\beta(\mc{h}))$. Similarly, if $\mc{h}\in\his{B}{}$, define $\mcp{h}{\phi}=(\mc{h},\phi(\mc{h}))$. For any history \mc{h}, let \mcp{h}{\beta,\phi} denote the history obtained from \mc{h} by letting Buster play according to $\beta$ and Fixer play according to $\phi$ until either Fixer or Buster wins.

Say Fixer strategy $\phi$ \emph{Fixer-dominates} Fixer strategy $\phi'$ at history $\mc{h}\in\his{B}{}$ if for any Buster strategy $\beta$ there exists a Buster strategy $\beta'$ such that \mcp{h}{\beta,\phi} is Fixer-superior to \mcp{h}{\beta',\phi'} (noting that $\beta$ and $\beta'$ obey the same given limit $\ell$). Say $\phi$ is \emph{Fixer-dominant} at \mc{h} if it Fixer-dominates all other Fixer strategies at \mc{h}; that is, $\phi$ is Fixer-dominant at \mc{h} if for any other Fixer strategy $\phi'$, the minimum over all Buster strategies $\beta$ of $\pay{\mcp{h}{\beta,\phi}}$ is greater than or equal to the minimum over all Buster strategies $\beta'$ of $\pay{\mcp{h}{\beta',\phi'}}$. Note that every complete history is Fixer-superior to itself, so any Fixer strategy Fixer-dominates itself everywhere.

\subsection{Example Gameplay}
\label{example}

Suppose $G=\{e_1,e_2,e_3\}$, where these edges form a triangle, and $R=\{e_4,e_5\}$, where $e_4$ has the same endpoints as $e_1$ and satisfies $w(e_4)=1$, while $e_5$ has the same endpoints as $e_2$ and satisfies and $w(e_5)=2$; see Figure \ref{exampleGraph}.

\begin{figure}[htb]
\centering
\begin{tikzpicture}
\Vertex[x=0,y=2,fontsize=\large]{a}
\Vertex[x=-1.7,y=-1,fontsize=\large]{b}
\Vertex[x=1.7,y=-1,fontsize=\large]{c}
\Edge[label=$e_1$,fontsize=\large](a)(b)
\Edge[label=$e_2$,fontsize=\large](b)(c)
\Edge[label=$e_3$,fontsize=\large](c)(a)
\Edge[label=$e_4$,bend=-45,style={dashed},fontsize=\large](a)(b)
\Edge[label=$e_5$,bend=-45,style={dashed},fontsize=\large](b)(c)
\end{tikzpicture}
\caption{The solid edges form the graph $G$, while the dashed edges form the reserve set $R$, with $w(e_4)=1$ and $w(e_5)=2$.}\label{exampleGraph}
\end{figure}

Consider the history $\mc{h}=((\ess{G}{1}=G,\ess{R}{1}=R),\ess{B}{1}=\{e_1,e_2\})\in\his{B}{1}$ having limit $\ell$. We detail how three Fixer strategies $\phi$, $\phi'$, and $\phi''$ extend \mc{h} for any value of $\ell$ in Figures \ref{example1}, \ref{example2}, and \ref{example3}, respectively. These Fixer strategies are displayed as decision trees going from top to bottom, with the vertices being the shaded boxes containing information about the state of the game between rounds and each edge representing the gameplay of some round $k$ consisting of Buster move \ess{B}{k} and Fixer response \ess{F}{k}. The internal vertices follow rounds during which Buster does not win (except for the root, which starts the game) and include the present values of the sets \ess{G}{k+1} and \ess{R}{k+1} after the $k$th round, as well as the total number $\sum _{j=1}^{k}|\ess{B}{j}|$ of edges removed by Buster through the first $k$ rounds, plus the value of the Fixer payout function $\pi$ for the complete history represented by the path from the root to that vertex if Fixer wins in the $(k+1)$st round because $\ell=\sum _{j=1}^{k}|\ess{B}{j}|$. The leaf vertices follow rounds during which Buster does win and include the set of edges in the disconnected graph $(\ess{G}{k}-\ess{B}{k})\cup\ess{R}{k}$ left after the $k$th round, as well as the total number $\sum _{j=1}^{k}|\ess{B}{j}|$ of edges removed by Buster, plus the value of the Fixer payout function $\pi$ for the complete history represented by the path from the root to that leaf. In other words, the path from the root to a given internal vertex whose last edge is the gameplay from the $k$th round represents a Fixer victory in the $(k+1)$st round, with \ess{G}{k+1} connected, the limit $\ell$ satisfying $\ell=\sum _{j=1}^{k}|\ess{B}{j}|$, and Fixer's payout $\pi$ equaling $1+w(\ess{R}{k+1})$, while the path from the root to a given leaf vertex whose last edge is the gameplay from the $k$th round represents a Buster victory in the $k$th round, with $(\ess{G}{k}-\ess{B}{k})\cup\ess{R}{k}$ disconnected, the limit $\ell$ satisfying $\ell\geq\sum _{j=1}^{k}|\ess{B}{j}|$, and Fixer's payout $\pi$ equaling $-1-w(\ess{R}{1}-\ess{R}{k})=w(\ess{R}{k})-4$.

Note that $\phi(\mc{h})=\{e_4\}$, $\phi'(\mc{h})=\{e_5\}$, $\phi''(\mc{h})=\{e_4,e_5\}$ are Fixer's only choices for \ess{F}{1} in responding to Buster's move $\ess{B}{1}=\{e_1,e_2\}$ in the first round (i.e. the only subsets of \ess{R}{1} whose addition to $\ess{G}{1}-\ess{B}{1}$ creates a connected graph \ess{G}{2}), and further note that the rest of $\phi'$ and $\phi''$ are forced (i.e. past the first round Fixer has only one choice for \ess{F}{k}, as detailed in the captions of Figures \ref{example2} and \ref{example3}). Also see that for any limit $\ell$ of \mc{h}, for each path from the root to a vertex in Figure \ref{example1} that could represent \mcp{h}{\beta,\phi} for some Buster strategy $\beta$ (i.e. $\ell=\sum _{j=1}^{|\mcp{h}{\beta,\phi}|-1}|\es{B}{\mcp{h}{\beta,\phi}}{j}|$ for paths to internal vertices and $\ell\geq\sum _{j=1}^{|\mcp{h}{\beta,\phi}|}|\es{B}{\mcp{h}{\beta,\phi}}{j}|$ for paths to leaves), the path from the root to the vertex in the same relative position in Figure \ref{example2} represents \mcp{h}{\beta',\phi'} for some Buster strategy $\beta'$ such that \mcp{h}{\beta,\phi} is Fixer-superior to \mcp{h}{\beta',\phi'} (i.e. $\pay{\mcp{h}{\beta,\phi}}\geq\pay{\mcp{h}{\beta',\phi'}}$). Similarly, it can be seen that for each path from the root to a vertex in Figure \ref{example1}, representing \mcp{h}{\beta,\phi} for some Buster strategy $\beta$, there is a path from the root to a vertex somewhere in Figure \ref{example3} representing \mcp{h}{\beta'',\phi''} for some Buster strategy $\beta''$ such that \mcp{h}{\beta,\phi} is Fixer-superior to \mcp{h}{\beta'',\phi''}.

Thus we have shown that regardless of the value of the limit $\ell$, for any Buster strategy $\beta$ and the only other Fixer strategies $\phi'$ and $\phi''$ besides $\phi$, there exist Buster strategies $\beta'$ and $\beta''$ such that \mcp{h}{\beta,\phi} is Fixer-superior to \mcp{h}{\beta',\phi'} and \mcp{h}{\beta'',\phi''}. Hence $\phi$ is Fixer-dominant at \mc{h}, regardless of the limit of \mc{h}. Finally, see that $\phi'$ and $\phi''$ are not Fixer-dominant when the limit is $\ell=2$, because after any legal play by Fixer at \mc{h} (which must exist since $(\ess{G}{1}-\ess{B}{1})\cup\ess{R}{1}=\{e_3,e_4,e_5\}$ is connected, preventing Buster from winning in the first round), Fixer wins in the second round (since $|\ess{B}{1}|=2=\ell$), and for any Buster strategy $\beta^*$ neither \mcp{h}{\beta^*,\phi'} nor \mcp{h}{\beta^*,\phi''} is Fixer-superior to \mcp{h}{\beta^*,\phi} (as these gameplays correspond respectively to the paths from the root to its only child in each of Figures \ref{example2}, \ref{example3}, and \ref{example1}, with $\pay{ \mcp{h}{\beta^*,\phi}}=3$, $\pay{ \mcp{h}{\beta^*,\phi'}}=2$, and $\pay{ \mcp{h}{\beta^*,\phi''}}=1$).

\begin{figure}[htb]
\centering
\begin{tikzpicture}
\Vertex[x=0,y=0,shape=rectangle,style={minimum width=2.25cm,minimum height=1.25cm},label={\begin{tabular}{c}$\ess{G}{1}=\{e_1,e_2,e_3\}$ \tabularnewline $\ess{R}{1}=\{e_4,e_5\}$ \tabularnewline $\sum |\ess{B}{i}|=0$ \tabularnewline $\pi=4$ if $\ell=0$ \end{tabular}}]{v1}
\Vertex[x=0,y=-1.25,shape=rectangle,style={minimum width=2.25cm,minimum height=.75cm,opacity=0},label={\begin{tabular}{c}$\ess{B}{1}=\{e_1,e_2\}$ \tabularnewline $\ess{F}{1}=\{e_4\}$\end{tabular}}]{e11}
\Vertex[x=0,y=-2.5,shape=rectangle,style={minimum width=2.25cm,minimum height=1.25cm},label={\begin{tabular}{c}$\ess{G}{2}=\{e_3,e_4\}$ \tabularnewline $\ess{R}{2}=\{e_5\}$ \tabularnewline $\sum |\ess{B}{i}|=2$ \tabularnewline $\pi=3$ if $\ell=2$ \end{tabular}}]{v11}
\Edge(v1)(e11)
\Edge(e11)(v11)
\Vertex[x=-2.375,y=-3.75,shape=rectangle,style={minimum width=2.25cm,minimum height=.75cm,opacity=0},label={\begin{tabular}{c}$\ess{B}{2}=\{e_3\}$ \tabularnewline $\ess{F}{2}=\{e_5\}$\end{tabular}}]{e111}
\Vertex[x=-4.75,y=-5,shape=rectangle,style={minimum width=2.25cm,minimum height=1.25cm},label={\begin{tabular}{c}$\ess{G}{3}=\{e_4,e_5\}$ \tabularnewline $\ess{R}{3}=\{\}$ \tabularnewline $\sum |\ess{B}{i}|=3$ \tabularnewline $\pi=1$ if $\ell=3$\end{tabular}}]{v111}
\Edge(v11)(e111)
\Edge(e111)(v111)
\Vertex[x=0,y=-3.75,shape=rectangle,style={minimum width=2.25cm,minimum height=.75cm,opacity=0},label={$\ess{B}{2}=\{e_3,e_4\}$}]{e112}
\Vertex[x=0,y=-5,shape=rectangle,style={minimum width=2.25cm,minimum height=1.25cm},label={\begin{tabular}{c}$(\ess{G}{2}-\ess{B}{2})\cup\ess{R}{2}$ \tabularnewline $=\{e_5\}$ \tabularnewline $\sum |\ess{B}{i}|=4$ \tabularnewline $\pi=-2$ if $\ell\geq 4$\end{tabular}}]{v112}
\Edge(v11)(e112)
\Edge(e112)(v112)
\Vertex[x=2.375,y=-3.75,shape=rectangle,style={minimum width=2.25cm,minimum height=.75cm,opacity=0},label={\begin{tabular}{c}$\ess{B}{2}=\{e_4\}$ \tabularnewline $\ess{F}{2}=\{e_5\}$\end{tabular}}]{e113}
\Vertex[x=4.75,y=-5,shape=rectangle,style={minimum width=2.25cm,minimum height=1.25cm},label={\begin{tabular}{c}$\ess{G}{3}=\{e_3,e_5\}$ \tabularnewline $\ess{R}{3}=\{\}$ \tabularnewline $\sum |\ess{B}{i}|=3$ \tabularnewline $\pi=1$ if $\ell=3$\end{tabular}}]{v113}
\Edge(v11)(e113)
\Edge(e113)(v113)
\Vertex[x=-6.75,y=-6.25,shape=rectangle,style={minimum width=2.25cm,minimum height=.75cm,opacity=0},label={$\ess{B}{3}=\{e_4\}$}]{e1111}
\Vertex[x=-8.25,y=-7.5,shape=rectangle,style={minimum width=2.25cm,minimum height=1.25cm},label={\begin{tabular}{c}$(\ess{G}{3}-\ess{B}{3})\cup\ess{R}{3}$ \tabularnewline $=\{e_5\}$ \tabularnewline $\sum |\ess{B}{i}|=4$ \tabularnewline $\pi=-4$ if $\ell\geq 4$\end{tabular}}]{v1111}
\Edge(v111)(e1111)
\Edge(e1111)(v1111)
\Vertex[x=-4.75,y=-6.25,shape=rectangle,style={minimum width=2.25cm,minimum height=.75cm,opacity=0},label={$\ess{B}{3}=\{e_5\}$}]{e1112}
\Vertex[x=-4.75,y=-7.5,shape=rectangle,style={minimum width=2.25cm,minimum height=1.25cm},label={\begin{tabular}{c}$(\ess{G}{3}-\ess{B}{3})\cup\ess{R}{3}$ \tabularnewline $=\{e_4\}$ \tabularnewline $\sum |\ess{B}{i}|=4$ \tabularnewline $\pi=-4$ if $\ell\geq 4$\end{tabular}}]{v1112}
\Edge(v111)(e1112)
\Edge(e1112)(v1112)
\Vertex[x=-3,y=-6.25,shape=rectangle,style={minimum width=2.25cm,minimum height=.75cm,opacity=0},label={$\ess{B}{3}=\{e_4,e_5\}$}]{e1113}
\Vertex[x=-1.25,y=-7.5,shape=rectangle,style={minimum width=2.25cm,minimum height=1.25cm},label={\begin{tabular}{c}$(\ess{G}{3}-\ess{B}{3})\cup\ess{R}{3}$ \tabularnewline $=\{\}$ \tabularnewline $\sum |\ess{B}{i}|=5$ \tabularnewline $\pi=-4$ if $\ell\geq 5$\end{tabular}}]{v1113}
\Edge(v111)(e1113)
\Edge(e1113)(v1113)
\Vertex[x=3,y=-6.25,shape=rectangle,style={minimum width=2.25cm,minimum height=.75cm,opacity=0},label={$\ess{B}{3}=\{e_3\}$}]{e1131}
\Vertex[x=1.25,y=-7.5,shape=rectangle,style={minimum width=2.25cm,minimum height=1.25cm},label={\begin{tabular}{c}$(\ess{G}{3}-\ess{B}{3})\cup\ess{R}{3}$ \tabularnewline $=\{e_5\}$ \tabularnewline $\sum |\ess{B}{i}|=4$ \tabularnewline $\pi=-4$ if $\ell\geq 4$\end{tabular}}]{v1131}
\Edge(v113)(e1131)
\Edge(e1131)(v1131)
\Vertex[x=4.75,y=-6.25,shape=rectangle,style={minimum width=2.25cm,minimum height=.75cm,opacity=0},label={$\ess{B}{3}=\{e_5\}$}]{e1132}
\Vertex[x=4.75,y=-7.5,shape=rectangle,style={minimum width=2.25cm,minimum height=1.25cm},label={\begin{tabular}{c}$(\ess{G}{3}-\ess{B}{3})\cup\ess{R}{3}$ \tabularnewline $=\{e_3\}$ \tabularnewline $\sum |\ess{B}{i}|=4$ \tabularnewline $\pi=-4$ if $\ell\geq 4$\end{tabular}}]{v1132}
\Edge(v113)(e1132)
\Edge(e1132)(v1132)
\Vertex[x=6.75,y=-6.25,shape=rectangle,style={minimum width=2.25cm,minimum height=.75cm,opacity=0},label={$\ess{B}{3}=\{e_3,e_5\}$}]{e1133}
\Vertex[x=8.25,y=-7.5,shape=rectangle,style={minimum width=2.25cm,minimum height=1.25cm},label={\begin{tabular}{c}$(\ess{G}{3}-\ess{B}{3})\cup\ess{R}{3}$ \tabularnewline $=\{\}$ \tabularnewline $\sum |\ess{B}{i}|=5$ \tabularnewline $\pi=-4$ if $\ell\geq 5$\end{tabular}}]{v1133}
\Edge(v113)(e1133)
\Edge(e1133)(v1133)
\end{tikzpicture}
\caption{A decision tree for $\phi$, where $\phi(h)=\{e_4\}$.}\label{example1}
\end{figure}

\begin{figure}[htb]
\centering
\begin{tikzpicture}
\Vertex[x=0,y=0,shape=rectangle,style={minimum width=2.25cm,minimum height=1.25cm},label={\begin{tabular}{c}$\ess{G}{1}=\{e_1,e_2,e_3\}$ \tabularnewline $\ess{R}{1}=\{e_4,e_5\}$ \tabularnewline $\sum |\ess{B}{i}|=0$ \tabularnewline $\pi=4$ if $\ell=0$\end{tabular}}]{v1}
\Vertex[x=0,y=-1.25,shape=rectangle,style={minimum width=2.25cm,minimum height=.75cm,opacity=0},label={\begin{tabular}{c}$\ess{B}{1}=\{e_1,e_2\}$ \tabularnewline $\ess{F}{1}=\{e_5\}$\end{tabular}}]{e11}
\Vertex[x=0,y=-2.5,shape=rectangle,style={minimum width=2.25cm,minimum height=1.25cm},label={\begin{tabular}{c}$\ess{G}{2}=\{e_3,e_5\}$ \tabularnewline $\ess{R}{2}=\{e_4\}$ \tabularnewline $\sum |\ess{B}{i}|=2$ \tabularnewline $\pi=2$ if $\ell=2$\end{tabular}}]{v11}
\Edge(v1)(e11)
\Edge(e11)(v11)
\Vertex[x=-2.375,y=-3.75,shape=rectangle,style={minimum width=2.25cm,minimum height=.75cm,opacity=0},label={\begin{tabular}{c}$\ess{B}{2}=\{e_3\}$ \tabularnewline $\ess{F}{2}=\{e_4\}$\end{tabular}}]{e111}
\Vertex[x=-4.75,y=-5,shape=rectangle,style={minimum width=2.25cm,minimum height=1.25cm},label={\begin{tabular}{c}$\ess{G}{3}=\{e_4,e_5\}$ \tabularnewline $\ess{R}{3}=\{\}$ \tabularnewline $\sum |\ess{B}{i}|=3$ \tabularnewline $\pi=1$ if $\ell=3$\end{tabular}}]{v111}
\Edge(v11)(e111)
\Edge(e111)(v111)
\Vertex[x=0,y=-3.75,shape=rectangle,style={minimum width=2.25cm,minimum height=.75cm,opacity=0},label={$\ess{B}{2}=\{e_3,e_5\}$}]{e112}
\Vertex[x=0,y=-5,shape=rectangle,style={minimum width=2.25cm,minimum height=1.25cm},label={\begin{tabular}{c}$(\ess{G}{2}-\ess{B}{2})\cup\ess{R}{2}$ \tabularnewline $=\{e_4\}$ \tabularnewline $\sum |\ess{B}{i}|=4$ \tabularnewline $\pi=-3$ if $\ell\geq 4$\end{tabular}}]{v112}
\Edge(v11)(e112)
\Edge(e112)(v112)
\Vertex[x=2.375,y=-3.75,shape=rectangle,style={minimum width=2.25cm,minimum height=.75cm,opacity=0},label={\begin{tabular}{c}$\ess{B}{2}=\{e_5\}$ \tabularnewline $\ess{F}{2}=\{e_4\}$\end{tabular}}]{e113}
\Vertex[x=4.75,y=-5,shape=rectangle,style={minimum width=2.25cm,minimum height=1.25cm},label={\begin{tabular}{c}$\ess{G}{3}=\{e_3,e_4\}$ \tabularnewline $\ess{R}{3}=\{\}$ \tabularnewline $\sum |\ess{B}{i}|=3$ \tabularnewline $\pi=1$ if $\ell=3$\end{tabular}}]{v113}
\Edge(v11)(e113)
\Edge(e113)(v113)
\Vertex[x=-6.75,y=-6.25,shape=rectangle,style={minimum width=2.25cm,minimum height=.75cm,opacity=0},label={$\ess{B}{3}=\{e_4\}$}]{e1111}
\Vertex[x=-8.25,y=-7.5,shape=rectangle,style={minimum width=2.25cm,minimum height=1.25cm},label={\begin{tabular}{c}$(\ess{G}{3}-\ess{B}{3})\cup\ess{R}{3}$ \tabularnewline $=\{e_5\}$ \tabularnewline $\sum |\ess{B}{i}|=4$ \tabularnewline $\pi=-4$ if $\ell\geq 4$\end{tabular}}]{v1111}
\Edge(v111)(e1111)
\Edge(e1111)(v1111)
\Vertex[x=-4.75,y=-6.25,shape=rectangle,style={minimum width=2.25cm,minimum height=.75cm,opacity=0},label={$\ess{B}{3}=\{e_5\}$}]{e1112}
\Vertex[x=-4.75,y=-7.5,shape=rectangle,style={minimum width=2.25cm,minimum height=1.25cm},label={\begin{tabular}{c}$(\ess{G}{3}-\ess{B}{3})\cup\ess{R}{3}$ \tabularnewline $=\{e_4\}$ \tabularnewline $\sum |\ess{B}{i}|=4$ \tabularnewline $\pi=-4$ if $\ell\geq 4$\end{tabular}}]{v1112}
\Edge(v111)(e1112)
\Edge(e1112)(v1112)
\Vertex[x=-3,y=-6.25,shape=rectangle,style={minimum width=2.25cm,minimum height=.75cm,opacity=0},label={$\ess{B}{3}=\{e_4,e_5\}$}]{e1113}
\Vertex[x=-1.25,y=-7.5,shape=rectangle,style={minimum width=2.25cm,minimum height=1.25cm},label={\begin{tabular}{c}$(\ess{G}{3}-\ess{B}{3})\cup\ess{R}{3}$ \tabularnewline $=\{\}$ \tabularnewline $\sum |\ess{B}{i}|=5$ \tabularnewline $\pi=-4$ if $\ell\geq 5$\end{tabular}}]{v1113}
\Edge(v111)(e1113)
\Edge(e1113)(v1113)
\Vertex[x=3,y=-6.25,shape=rectangle,style={minimum width=2.25cm,minimum height=.75cm,opacity=0},label={$\ess{B}{3}=\{e_3\}$}]{e1131}
\Vertex[x=1.25,y=-7.5,shape=rectangle,style={minimum width=2.25cm,minimum height=1.25cm},label={\begin{tabular}{c}$(\ess{G}{3}-\ess{B}{3})\cup\ess{R}{3}$ \tabularnewline $=\{e_4\}$ \tabularnewline $\sum |\ess{B}{i}|=4$ \tabularnewline $\pi=-4$ if $\ell\geq 4$\end{tabular}}]{v1131}
\Edge(v113)(e1131)
\Edge(e1131)(v1131)
\Vertex[x=4.75,y=-6.25,shape=rectangle,style={minimum width=2.25cm,minimum height=.75cm,opacity=0},label={$\ess{B}{3}=\{e_4\}$}]{e1132}
\Vertex[x=4.75,y=-7.5,shape=rectangle,style={minimum width=2.25cm,minimum height=1.25cm},label={\begin{tabular}{c}$(\ess{G}{3}-\ess{B}{3})\cup\ess{R}{3}$ \tabularnewline $=\{e_3\}$ \tabularnewline $\sum |\ess{B}{i}|=4$ \tabularnewline $\pi=-4$ if $\ell\geq 4$\end{tabular}}]{v1132}
\Edge(v113)(e1132)
\Edge(e1132)(v1132)
\Vertex[x=6.75,y=-6.25,shape=rectangle,style={minimum width=2.25cm,minimum height=.75cm,opacity=0},label={$\ess{B}{3}=\{e_3,e_4\}$}]{e1133}
\Vertex[x=8.25,y=-7.5,shape=rectangle,style={minimum width=2.25cm,minimum height=1.25cm},label={\begin{tabular}{c}$(\ess{G}{3}-\ess{B}{3})\cup\ess{R}{3}$ \tabularnewline $=\{\}$ \tabularnewline $\sum |\ess{B}{i}|=5$ \tabularnewline $\pi=-4$ if $\ell\geq 5$\end{tabular}}]{v1133}
\Edge(v113)(e1133)
\Edge(e1133)(v1133)
\end{tikzpicture}
\caption{A decision tree for $\phi'$, where $\phi'(h)=\{e_5\}$. Observe that the rest of the decision tree is forced. Indeed, $\ess{G}{2}=\{e_3,e_5\}$, which is a path graph, and $\ess{R}{2}=\{e_4\}$, where $e_4$ spans the endpoints of that path. Thus Buster deleting any single edge from \ess{G}{2} necessitates Fixer reconnecting the graph using $e_4$, from which point on any further deletions by Buster cannot be countered by Fixer, and Buster initially deleting both edges from \ess{G}{2} also cannot be countered by Fixer.}\label{example2}
\end{figure}

\begin{figure}[htb]
\centering
\begin{tikzpicture}
\Vertex[x=0,y=0,shape=rectangle,style={minimum width=2.25cm,minimum height=1.25cm},label={\begin{tabular}{c}$\ess{G}{1}=\{e_1,e_2,e_3\}$ \tabularnewline $\ess{R}{1}=\{e_4,e_5\}$ \tabularnewline $\sum |\ess{B}{i}|=0$ \tabularnewline $\pi=4$ if $\ell=0$\end{tabular}}]{v1}
\Vertex[x=0,y=-1.25,shape=rectangle,style={minimum width=2.25cm,minimum height=.75cm,opacity=0},label={\begin{tabular}{c}$\ess{B}{1}=\{e_1,e_2\}$ \tabularnewline $\ess{F}{1}=\{e_4,e_5\}$\end{tabular}}]{e11}
\Vertex[x=0,y=-2.5,shape=rectangle,style={minimum width=2.25cm,minimum height=1.25cm},label={\begin{tabular}{c}$\ess{G}{2}=\{e_3,e_4,e_5\}$ \tabularnewline $\ess{R}{2}=\{\}$ \tabularnewline $\sum |\ess{B}{i}|=2$ \tabularnewline $\pi=1$ if $\ell=2$\end{tabular}}]{v11}
\Edge(v1)(e11)
\Edge(e11)(v11)
\Vertex[x=-5.30357,y=-4.75,shape=rectangle,style={minimum width=2.25cm,minimum height=.75cm,opacity=0},label={$\ess{B}{2}=\{e_3,e_4\}$}]{e111}
\Vertex[x=-8.25,y=-6,shape=rectangle,style={minimum width=2.25cm,minimum height=1.25cm},label={\begin{tabular}{c}$(\ess{G}{2}-\ess{B}{2})\cup\ess{R}{2}$ \tabularnewline $=\{e_5\}$ \tabularnewline $\sum |\ess{B}{i}|=4$ \tabularnewline $\pi=-4$ if $\ell\geq 4$\end{tabular}}]{v111}
\Edge(v11)(e111)
\Edge(e111)(v111)
\Vertex[x=-3.53571,y=-4.75,shape=rectangle,style={minimum width=2.25cm,minimum height=.75cm,opacity=0},label={\begin{tabular}{c}$\ess{B}{2}=\{e_3\}$ \tabularnewline $\ess{F}{2}=\{\}$\end{tabular}}]{e112}
\Vertex[x=-5.5,y=-6,shape=rectangle,style={minimum width=2.25cm,minimum height=1.25cm},label={\begin{tabular}{c}$\ess{G}{3}=\{e_4,e_5\}$ \tabularnewline $\ess{R}{3}=\{\}$ \tabularnewline $\sum |\ess{B}{i}|=3$ \tabularnewline $\pi=1$ if $\ell=3$\end{tabular}}]{v112}
\Edge(v11)(e112)
\Edge(e112)(v112)
\Vertex[x=-1.76786,y=-4.75,shape=rectangle,style={minimum width=2.25cm,minimum height=.75cm,opacity=0},label={$\ess{B}{2}=\{e_3,e_5\}$}]{e113}
\Vertex[x=-2.75,y=-6,shape=rectangle,style={minimum width=2.25cm,minimum height=1.25cm},label={\begin{tabular}{c}$(\ess{G}{2}-\ess{B}{2})\cup\ess{R}{2}$ \tabularnewline $=\{e_4\}$ \tabularnewline $\sum |\ess{B}{i}|=4$ \tabularnewline $\pi=-4$ if $\ell\geq 4$\end{tabular}}]{v113}
\Edge(v11)(e113)
\Edge(e113)(v113)
\Vertex[x=0,y=-4.75,shape=rectangle,style={minimum width=2.25cm,minimum height=.75cm,opacity=0},label={\begin{tabular}{c}$\ess{B}{2}=\{e_4\}$ \tabularnewline $\ess{F}{2}=\{\}$\end{tabular}}]{e114}
\Vertex[x=0,y=-6,shape=rectangle,style={minimum width=2.25cm,minimum height=1.25cm},label={\begin{tabular}{c}$\ess{G}{3}=\{e_3,e_5\}$ \tabularnewline $\ess{R}{3}=\{\}$ \tabularnewline $\sum |\ess{B}{i}|=3$ \tabularnewline $\pi=1$ if $\ell=3$\end{tabular}}]{v114}
\Edge(v11)(e114)
\Edge(e114)(v114)
\Vertex[x=1.76786,y=-4.75,shape=rectangle,style={minimum width=2.25cm,minimum height=.75cm,opacity=0},label={$\ess{B}{2}=\{e_4,e_5\}$}]{e115}
\Vertex[x=2.75,y=-6,shape=rectangle,style={minimum width=2.25cm,minimum height=1.25cm},label={\begin{tabular}{c}$(\ess{G}{2}-\ess{B}{2})\cup\ess{R}{2}$ \tabularnewline $=\{e_3\}$ \tabularnewline $\sum |\ess{B}{i}|=4$ \tabularnewline $\pi=-4$ if $\ell\geq 4$\end{tabular}}]{v115}
\Edge(v11)(e115)
\Edge(e115)(v115)
\Vertex[x=3.53571,y=-4.75,shape=rectangle,style={minimum width=2.25cm,minimum height=.75cm,opacity=0},label={\begin{tabular}{c}$\ess{B}{2}=\{e_5\}$ \tabularnewline $\ess{F}{2}=\{\}$\end{tabular}}]{e116}
\Vertex[x=5.5,y=-6,shape=rectangle,style={minimum width=2.25cm,minimum height=1.25cm},label={\begin{tabular}{c}$\ess{G}{3}=\{e_3,e_4\}$ \tabularnewline $\ess{R}{3}=\{\}$ \tabularnewline $\sum |\ess{B}{i}|=3$ \tabularnewline $\pi=1$ if $\ell=3$\end{tabular}}]{v116}
\Edge(v11)(e116)
\Edge(e116)(v116)
\Vertex[x=5.30357,y=-4.75,shape=rectangle,style={minimum width=2.25cm,minimum height=.75cm,opacity=0},label={$\ess{B}{2}=\{e_3,e_4,e_5\}$}]{e117}
\Vertex[x=8.25,y=-6,shape=rectangle,style={minimum width=2.25cm,minimum height=1.25cm},label={\begin{tabular}{c}$(\ess{G}{2}-\ess{B}{2})\cup\ess{R}{2}$ \tabularnewline $=\{\}$ \tabularnewline $\sum |\ess{B}{i}|=5$ \tabularnewline $\pi=-4$ if $\ell\geq 5$\end{tabular}}]{v117}
\Edge(v11)(e117)
\Edge(e117)(v117)
\Vertex[x=-6.1875,y=-7.25,shape=rectangle,style={minimum width=2.25cm,minimum height=.75cm,opacity=0},label={$\ess{B}{3}=\{e_4\}$}]{e1121}
\Vertex[x=-6.875,y=-8.5,shape=rectangle,style={minimum width=2.25cm,minimum height=1.25cm},label={\begin{tabular}{c}$(\ess{G}{3}-\ess{B}{3})\cup\ess{R}{3}$\tabularnewline $=\{e_5\}$ \tabularnewline $\sum |\ess{B}{i}|=4$ \tabularnewline $\pi=-4$ if $\ell\geq 4$\end{tabular}}]{v1121}
\Edge(v112)(e1121)
\Edge(e1121)(v1121)
\Vertex[x=-4.8125,y=-7.25,shape=rectangle,style={minimum width=2.25cm,minimum height=.75cm,opacity=0},label={$\ess{B}{3}=\{e_5\}$}]{e1122}
\Vertex[x=-4.125,y=-8.5,shape=rectangle,style={minimum width=2.25cm,minimum height=1.25cm},label={\begin{tabular}{c}$(\ess{G}{3}-\ess{B}{3})\cup\ess{R}{3}$ \tabularnewline $=\{e_4\}$ \tabularnewline $\sum |\ess{B}{i}|=4$ \tabularnewline $\pi=-4$ if $\ell\geq 4$\end{tabular}}]{v1122}
\Edge(v112)(e1122)
\Edge(e1122)(v1122)
\Vertex[x=-.6875,y=-7.25,shape=rectangle,style={minimum width=2.25cm,minimum height=.75cm,opacity=0},label={$\ess{B}{3}=\{e_3\}$}]{e1141}
\Vertex[x=-1.375,y=-8.5,shape=rectangle,style={minimum width=2.25cm,minimum height=1.25cm},label={\begin{tabular}{c}$(\ess{G}{3}-\ess{B}{3})\cup\ess{R}{3}$\tabularnewline $=\{e_5\}$ \tabularnewline $\sum |\ess{B}{i}|=4$ \tabularnewline $\pi=-4$ if $\ell\geq 4$\end{tabular}}]{v1141}
\Edge(v114)(e1141)
\Edge(e1141)(v1141)
\Vertex[x=.6875,y=-7.25,shape=rectangle,style={minimum width=2.25cm,minimum height=.75cm,opacity=0},label={$\ess{B}{3}=\{e_5\}$}]{e1142}
\Vertex[x=1.375,y=-8.5,shape=rectangle,style={minimum width=2.25cm,minimum height=1.25cm},label={\begin{tabular}{c}$(\ess{G}{3}-\ess{B}{3})\cup\ess{R}{3}$\tabularnewline $=\{e_3\}$ \tabularnewline $\sum |\ess{B}{i}|=4$ \tabularnewline $\pi=-4$ if $\ell\geq 4$\end{tabular}}]{v1142}
\Edge(v114)(e1142)
\Edge(e1142)(v1142)
\Vertex[x=4.8125,y=-7.25,shape=rectangle,style={minimum width=2.25cm,minimum height=.75cm,opacity=0},label={$\ess{B}{3}=\{e_3\}$}]{e1161}
\Vertex[x=4.125,y=-8.5,shape=rectangle,style={minimum width=2.25cm,minimum height=1.25cm},label={\begin{tabular}{c}$(\ess{G}{3}-\ess{B}{3})\cup\ess{R}{3}$ \tabularnewline $=\{e_4\}$ \tabularnewline $\sum |\ess{B}{i}|=4$ \tabularnewline $\pi=-4$ if $\ell\geq 4$\end{tabular}}]{v1161}
\Edge(v116)(e1161)
\Edge(e1161)(v1161)
\Vertex[x=6.1875,y=-7.25,shape=rectangle,style={minimum width=2.25cm,minimum height=.75cm,opacity=0},label={$\ess{B}{3}=\{e_4\}$}]{e1162}
\Vertex[x=6.875,y=-8.5,shape=rectangle,style={minimum width=2.25cm,minimum height=1.25cm},label={\begin{tabular}{c}$(\ess{G}{3}-\ess{B}{3})\cup\ess{R}{3}$ \tabularnewline $=\{e_3\}$ \tabularnewline $\sum |\ess{B}{i}|=4$ \tabularnewline $\pi=-4$ if $\ell\geq 4$\end{tabular}}]{v1162}
\Edge(v116)(e1162)
\Edge(e1162)(v1162)
\Vertex[x=-5.5,y=-9.75,shape=rectangle,style={minimum width=2.25cm,minimum height=.75cm,opacity=0},label={$\ess{B}{3}=\{e_4,e_5\}$}]{e11200}
\Vertex[x=-5.5,y=-11,shape=rectangle,style={minimum width=2.25cm,minimum height=1.25cm},label={\begin{tabular}{c}$(\ess{G}{3}-\ess{B}{3})\cup\ess{R}{3}$ \tabularnewline $=\{\}$ \tabularnewline $\sum |\ess{B}{i}|=5$ \tabularnewline $\pi=-4$ if $\ell\geq 5$\end{tabular}}]{v11200}
\Edge(v112)(e11200)
\Edge(e11200)(v11200)
\Vertex[x=0,y=-9.75,shape=rectangle,style={minimum width=2.25cm,minimum height=.75cm,opacity=0},label={$\ess{B}{3}=\{e_3,e_5\}$}]{e11400}
\Vertex[x=0,y=-11,shape=rectangle,style={minimum width=2.25cm,minimum height=1.25cm},label={\begin{tabular}{c}$(\ess{G}{3}-\ess{B}{3})\cup\ess{R}{3}$ \tabularnewline $=\{\}$ \tabularnewline $\sum |\ess{B}{i}|=5$ \tabularnewline $\pi=-4$ if $\ell\geq 5$\end{tabular}}]{v11400}
\Edge(v114)(e11400)
\Edge(e11400)(v11400)
\Vertex[x=5.5,y=-9.75,shape=rectangle,style={minimum width=2.25cm,minimum height=.75cm,opacity=0},label={$\ess{B}{3}=\{e_3,e_4\}$}]{e11600}
\Vertex[x=5.5,y=-11,shape=rectangle,style={minimum width=2.25cm,minimum height=1.25cm},label={\begin{tabular}{c}$(\ess{G}{3}-\ess{B}{3})\cup\ess{R}{3}$ \tabularnewline $=\{\}$ \tabularnewline $\sum |\ess{B}{i}|=5$ \tabularnewline $\pi=-4$ if $\ell\geq 5$\end{tabular}}]{v11600}
\Edge(v116)(e11600)
\Edge(e11600)(v11600)
\end{tikzpicture}
\caption{A decision tree for $\phi''$, where $\phi''(h)=\{e_4,e_5\}$. Again, observe that the rest of the decision tree is forced. Indeed, \ess{R}{2} is empty, leaving Fixer with no options besides playing the empty set when possible.}\label{example3}
\end{figure}

\subsection{Statement of Main Theorem}
\label{statement}

Fixer strategy $\phi$ is \emph{greedy} if $w(\phi(\mc{h}))\leq w(\phi'(\mc{h}))$ for any history $\mc{h}\in\his{B}{}$ and any other Fixer strategy $\phi'$. That is, Fixer plays greedily in response to a move by Buster by adding no reserve edge if the graph remains connected, and otherwise adding some cheapest set of reserve edges that connects the graph. We may also refer to individual Fixer moves that come from a greedy Fixer strategy as being greedy.

Note that any Fixer-dominant strategy necessarily makes a greedy move responding to any Buster move that reaches the limit $\ell$. Indeed, if $w(\phi(\mc{h})) >  w(\phi'(\mc{h}))$ for $\mc{h}\in\his{B}{k-1}$ with $\sum _{j=1}^{k-1}|\es{B}{\mc{h}}{j}|=\ell$ and $(\es{G}{\mc{h}}{k-1}-\es{B}{\mc{h}}{k-1})\cup\es{R}{\mc{h}}{k-1}$ connected, then for any two Buster strategies $\beta$ and $\beta'$, Fixer wins both \mcp{h}{\beta,\phi} and \mcp{h}{\beta',\phi'} in the $k$th round; hence \mcp{h}{\beta,\phi} is not Fixer-superior to \mcp{h}{\beta',\phi'} because $\pay{\mcp{h}{\beta,\phi}}=1+w(\es{R}{\mcp{h}{\beta,\phi}}{k})=1+w(\es{R}{\mc{h}}{k-1})-w(\phi(\mc{h}))<1+w(\es{R}{\mc{h}}{k-1})-w(\phi'(\mc{h}))=1+w(\es{R}{\mcp{h}{\beta',\phi'}}{k})=\pay{\mcp{h}{\beta',\phi'}}$. Our main theorem states that greediness is sufficient for Fixer-dominance.

\begin{thm}\label{mainthm}
Any greedy Fixer strategy is Fixer-dominant at any history \mc{h}.
\end{thm}

After establishing some facts about Fixer-superiority and Fixer-dominance in Section \ref{prelim}, we use them to prove Theorem \ref{mainthm} in Section \ref{proof}. The difficulty of the proof is in showing that if $\phi(\mc{h})$ and $\phi'(\mc{h})$ both reconnect $\es{G}{\mc{h}}{k}-\es{B}{\mc{h}}{k}$ for some $\mc{h}\in\his{B}{k}$, then Fixer immediately playing $\phi'(\mc{h})$ rather than $\phi(\mc{h})$ in the $k$th round does not provide a long-term advantage for maintaining graph connectedness in the face of Buster's future edge deletions (so Fixer might as well play $\phi(\mc{h})$ if it's cheaper). Our proof will be by induction, split into cases by the number of components of $\es{G}{\mc{h}}{k}-\es{B}{\mc{h}}{k}$. We shall see that the case of a single component is mostly trivial, the case of two components is the most difficult and requires case analysis of each move to verify that certain invariants are maintained, and the case of at least three components follows from a more direct application of the inductive hypothesis.

Note that the Fixer-dominance of a greedy Fixer strategy $\phi$ at a history \mc{h} such that $\es{G}{\mc{h}}{k}-\es{B}{\mc{h}}{k}$ has $c\geq 3$ components would not follow from iteratively applying the two component case $c-1$ times during the first round (e.g. for $c=3$ decomposing a greedy Fixer response reconnecting the graph in the first round into a greedy Fixer response reconnecting two of the components plus a greedy Fixer response reconnecting that now-connected graph with the third component). Indeed, this argument fails to prove that Fixer can't somehow exploit Buster removing edges all at once in a single round rather than over the course of multiple rounds.

\subsection{Related Work}
\label{past}

A \emph{bridge} in a multigraph $M$ is an edge $e$ such $M-\{e\}$ has one more component than $M$; equivalently, $e$ is part of no cycle in $M$. A \emph{spanning tree} of a connected multigraph $M$ is a subgraph $T$ of $M$ such that $T$ is a tree (i.e. connected and acyclic) whose vertex set matches that of $M$. A \emph{minimum spanning tree} of an edge-weighted multigraph $M$ is a spanning tree of $M$ minimizing the total weight of the edges. Minimum spanning trees are of interest to us because for $\mc{h}\in\his{B}{k}$, $\phi(\mc{h})$ is greedy if and only if it is a minimum spanning tree of the multigraph whose vertices are the components of $\es{G}{\mc{h}}{k}-\es{B}{\mc{h}}{k}$ and whose edges are the edges of \es{R}{\mc{h}}{k} (identifying each endpoint of the edges in \es{R}{\mc{h}}{k} with the component of $\es{G}{\mc{h}}{k}-\es{B}{\mc{h}}{k}$ within which it lies). Prim's Algorithm (first discovered by Jarnik \cite{J} and later by Prim \cite{P} and Dijkstra \cite{D}) finds a minimum spanning tree $T$ of a weighted connected multigraph $M$ one edge at a time by the following construction: with $T$ initialized as any vertex, iteratively add to $T$ any cheapest edge of $M$ joining a vertex in $T$ to one not yet in $T$, until all vertices of $M$ are in $T$.

For a family $\mathcal{F}$ of subgraphs of the complete graph $K_n$, the unbiased \emph{Maker-Breaker game on $\mathcal{F}$} consists of players Maker and Breaker taking turns claiming edges of $K_n$ (see \cite{CE} for some notable early results, and \cite{HKSS} for a more recent survey). Maker wins by claiming all edges of some graph in $\mathcal{F}$, while Breaker wins if all edges are claimed before Maker wins (equivalently, Breaker wins by claiming an edge from each minimal member of $\mathcal{F}$). The family $\mathcal{F}$ most relevant for comparison of the Maker-Breaker game on $\mathcal{F}$ with our game is the family of connected spanning subgraphs of $K_n$ (it is shown in \cite{L} that Maker can win such a game in $n-1$ moves if $n\geq 4$; see \cite{CFGHL,FKK,KK} for more Maker-Breaker results on spanning trees). The gameplay of Maker-Breaker differs from that of Buster-Fixer in several obvious ways, including that Maker only needs to end up with a connected graph (while Fixer must maintain connectedness after each turn), Maker cannot replace edges claimed by Breaker (whereas the reserve edges Fixer may select could include edges with the same endpoints as the edges deleted by Buster), and the Maker-Breaker game does not typically include weighted edges (which do occur in the Buster-Fixer model).

Our Buster-Fixer model is also related to the various models of \emph{network interdiction} and \emph{fortification}. Network interdiction models typically involve an interdictor and evader/operator; the interdictor modifies the graph (in ways such as reducing capacity or increasing flow cost along edges, or placing detection devices) for purposes of preventing the evader/operator from maximizing flow of objects from a source to a sink. Network fortification models involve the additional step, before the interdictor acts, of the evader/operator fortifying the graph (by such means as increasing capacity or reducing flow cost along edges, or defending certain vertices and edges from sabotage by the interdictor). Interdictors might be rational actors (e.g. terrorists) or random ones (e.g. natural disasters).

Specific examples of network interdiction include the following (see \cite{SPG} for more details). \emph{Shortest path interdiction} has the interdictor increase the cost of certain edges in the graph in order to maximize the evader's shortest path from a source vertex to a destination vertex. \emph{Maximum flow interdiction} has the interdictor reduce the capacity of certain edges in the graph in order to minimize the operator's maximum flow from a source vertex to a destination vertex. \emph{Facility assignment interdiction} has the interdictor remove $r$ vertices representing supply facilities in order to maximize the minimum weighted distance required for the operator to meet demand at each vertex representing a demand point. \emph{Most reliable path interdiction} has the interdictor place detection devices on the edges in order to maximize the minimum probability of detection sought by the evader. \emph{Minimum cost flow interdiction} has the interdictor reduce capacity on certain edges in the graph in order to maximize the minimum cost of the operator's flow that must satisfy capacity constraints along edges as well as conservation and demand constraints at vertices.

The literature on the network interdiction models mentioned above mostly revolves around seeking solutions via integer and linear programming. Related models, including one in \cite{BCD} that studies a game of attack and interception in a network where an attacker chooses a target and a path and each vertex chooses how much to invest in protection, are framed as more game-theoretic in nature; for instance, the main result of \cite{BCD} is the existence of a unique Nash equilibrium.

Also related are pursuit-evasion games such as \emph{Cops and Robbers} and \emph{Revolutionaries and Spies}. In Cops and Robbers (see \cite{ADG} for details and \cite{DG} for recent results), a cop chooses a vertex in a given graph, followed by the robber doing the same. In full view of both players, the cop and robber take turns either staying at their present vertex or moving to an adjacent one. The cop wins by moving to the same vertex as the robber, while the robber wins by avoiding the cop forever. In Revolutionaries and Spies (see \cite{BCPWZ} for details), $r$ revolutionaries choose vertices of a given graph to occupy, followed by $s$ spies doing the same, before the game proceeds in a series of rounds. In full view of all players, a round begins with each revolutionary either staying at their present vertex or moving to an adjacent one, followed by each spy doing the same. The revolutionaries win by accumulating $k$ revolutionaries at a single vertex without the any spy being able to occupy that vertex by the end of the round, while the spies win if they can continue the game forever.

Finally, we mention \emph{all-terminal network reliability}, which refers to the problem of determining the probability that a graph is connected when it is derived from an edge-weighted graph by randomly deleting edges according to the probability given by their weights. See \cite{GJ} for a polynomial-time approximation algorithm.

\subsection{Avenues for Future Research}
\label{future}

Many avenues exist for future research into variations on our Buster-Fixer model. In particular, we wonder about Fixer-dominant strategies for alternative games, where the condition that Fixer must maintain on the graph through each round is changed from maintaining connectedness to one of the following conditions:
\begin{enumerate}
\item Two given vertices $s$ and $t$ must stay in the same component.
\item The graph must stay $k$-connected for a given $k>1$.
\item Instead of a simple graph, the graph is directed, and Fixer must maintain one of the following conditions:
\begin{enumerate}
\item The directed graph must stay strongly connected.
\item The directed graph must have directed paths from (or to) a given vertex $s$ to (or from) all other vertices.
\item The directed graph must have a directed path from a given vertex $s$ to a given vertex $t$.
\item The directed graph must have directed paths in both directions between given vertices $s$ and $t$.
\end{enumerate}
\end{enumerate}

\section{Facts about Fixer-superiority and Fixer-dominance}
\label{prelim}

We first verify that Fixer-superiority and Fixer-dominance are transitive.

\begin{prop}\label{fixSupTransitivityProp}
Let \mc{h}, \mc{h'}, and \mc{h''} be complete histories such that \mc{h} is Fixer-superior to \mc{h'} and \mc{h'} is Fixer-superior to \mc{h''}. Then \mc{h} is Fixer-superior to \mc{h''}.
\end{prop}
\begin{proof}
$\pay{\mc{h}}\geq\pay{\mc{h'}}\geq\pay{\mc{h''}}$
\end{proof}

\begin{prop}\label{fixDomTransitivityCor}
For any $\mc{h}\in\his{B}{}$, if $\phi$ Fixer-dominates $\phi'$ at \mc{h} and $\phi'$ Fixer-dominates $\phi''$ at \mc{h}, then $\phi$ Fixer-dominates $\phi''$ at \mc{h}.
\end{prop}
\begin{proof}
Let $\beta$ be an arbitrary Buster strategy. Since $\phi$ Fixer-dominates $\phi'$ at \mc{h}, there exists a Buster strategy $\beta'$ such that \mcp{h}{\beta,\phi} is Fixer-superior to \mcp{h}{\beta',\phi'}. Since $\phi'$ Fixer-dominates $\phi''$ at \mc{h}, there exists a Buster strategy $\beta''$ such that \mcp{h}{\beta',\phi'} is Fixer-superior to \mcp{h}{\beta'',\phi''}. By Proposition \ref{fixSupTransitivityProp}, \mcp{h}{\beta,\phi} is Fixer-superior to \mcp{h}{\beta'',\phi''}, and thus $\phi$ Fixer-dominates $\phi''$ at \mc{h} because $\beta$ was arbitrary.
\end{proof}

We next show that Fixer following a Fixer-dominant strategy past a given move at $\mc{h}\in\his{B}{}$ guarantees Fixer-domination of any other Fixer strategy prescribing the same move at \mc{h}.

\begin{lem}\label{nextRound}
If $\mc{h}\in\his{B}{k}$ and $\phi$ and $\phi'$ are Fixer strategies satisfying $\phi(\mc{h})=\phi'(\mc{h})$, with $\phi$ Fixer-dominant at any history $\mc{h^*}\in\his{B}{}$ properly extending \mc{h}, then $\phi$ Fixer-dominates $\phi'$ at \mc{h}.
\end{lem}
\begin{proof}
Let $\beta$ be an arbitrary Buster strategy; to show that $\phi$ Fixer-dominates $\phi'$ at \mc{h} we show there exists a Buster strategy $\beta'$ such that \mcp{h}{\beta,\phi} is Fixer-superior to \mcp{h}{\beta',\phi'}. Set $\beta'=\beta$ if either Buster wins \mc{h} in the $k$th round, in which case $\mcp{h}{\beta,\phi}=\mc{h}=\mcp{h}{\beta',\phi'}$, or Fixer wins \mcp{h}{\phi} in the $(k+1)$st round, in which case $\mcp{h}{\beta,\phi}=\mcp{h}{\phi}=\mcp{h}{\phi'}=\mcp{h}{\beta',\phi'}$. Since every complete history Fixer-dominates itself, we may therefore assume $\beta$ has Buster make a move in the $(k+1)$st round of \mcp{h}{\beta,\phi}, so we can define the history $\mc{h^*}=(\mcp{h}{\phi},\beta(\mcp{h}{\phi}))$ properly extending \mc{h}. Since $\phi$ is Fixer-dominant at \mc{h^*} by hypothesis, $\phi$ Fixer-dominates $\phi'$ at \mc{h^*}. Therefore there exists a Buster strategy $\beta'$ such that \mcp{h^*}{\beta,\phi} is Fixer-superior to \mcp{h^*}{\beta',\phi'}; without loss of generality we may stipulate that $\beta'(\mcp{h}{\phi})=\beta(\mcp{h}{\phi})$ because the value of $\beta'$ at \mcp{h}{\phi} does not affect \mcp{h^*}{\beta',\phi'}, as \mc{h^*} is \mcp{h}{\phi} with the extra Buster move $\beta(\mcp{h}{\phi})$. Since \mc{h^*} is simply \mc{h} appended with the common Fixer response $\phi(\mc{h})$ prescribed by both $\phi$ and $\phi'$ as well as the common subsequent Buster play $\beta(\mcp{h}{\phi})$ prescribed by both $\beta$ and $\beta'$, we have $\mcp{h}{\beta,\phi}=\mcp{h^*}{\beta,\phi}$ and $\mcp{h}{\beta',\phi'}=\mcp{h^*}{\beta',\phi'}$. Therefore \mcp{h}{\beta,\phi} is Fixer-superior to \mcp{h}{\beta',\phi'}, so $\phi$ Fixer-dominates $\phi'$ at \mc{h}.
\end{proof}

We now present a method for proving inductively that a Fixer strategy is Fixer-dominant.

\begin{prop}\label{fixDomStrat}
Suppose $\mc{h}\in\his{B}{}$ and $\phi$ and $\phi^{\#}$ are Fixer strategies that are Fixer-dominant at any history properly extending \mc{h}. If there exists a Fixer strategy $\phi^*$ satisfying $\phi^*(\mc{h})=\phi(\mc{h})$ and Fixer-dominating $\phi^{\#}$ at \mc{h}, then $\phi$ Fixer-dominates $\phi'$ at \mc{h} for any Fixer strategy $\phi'$ satisfying $\phi'(\mc{h})=\phi^{\#}(\mc{h})$.
\end{prop}
\begin{proof}
Let $\phi'$ be any Fixer strategy satisfying $\phi'(\mc{h})=\phi^{\#}(\mc{h})$. By hypothesis there exists a Fixer strategy $\phi^*$ satisfying $\phi^*(\mc{h})=\phi(\mc{h})$ and Fixer-dominating $\phi^{\#}$ at \mc{h}. Since $\phi$ and $\phi^*$ are Fixer strategies with $\phi$ Fixer-dominant at any history properly extending \mc{h}, by Lemma \ref{nextRound} $\phi$ Fixer-dominates $\phi^*$ at \mc{h}. Similarly, since $\phi^{\#}$ and $\phi'$ are Fixer strategies with $\phi^{\#}$ Fixer-dominant at any history properly extending \mc{h}, by Lemma \ref{nextRound} $\phi^{\#}$ Fixer-dominates $\phi'$ at \mc{h}. Thus $\phi$ Fixer-dominates $\phi^*$ at \mc{h}, $\phi^*$ Fixer-dominates $\phi^{\#}$ at \mc{h}, and $\phi^{\#}$ Fixer-dominates $\phi'$ at \mc{h}; by Proposition \ref{fixDomTransitivityCor} $\phi$ Fixer-dominates $\phi'$ at \mc{h}.
\end{proof}

Call a Buster strategy \emph{granular} if it restricts Buster to removing only singletons. Since Buster removing multiple edges at once is, in a sense, equivalent to removing them one-at-a-time while taking away Buster's ability to change course partway through based on Fixer's early responses, it makes sense that granular strategies are advantageous for Buster. Hence we can state below that in order to show a particular Fixer move \ess{F}{k} at some history $\mc{h}\in\his{B}{k}$ can lead to a Fixer strategy $\phi^*$ that satisfies $\phi^*(\mc{h})=\ess{F}{k}$ and Fixer-dominates some given Fixer strategy $\phi^{\#}$ at \mc{h}, it suffices to construct a Fixer strategy $\phi^1$ such that $\phi^1(\mc{h})=\ess{F}{k}$ and for any granular Buster strategy $\beta^1$ there exists a Buster strategy $\beta^{\#}$ for which \mcp{h}{\beta^1,\phi^1} is Fixer-superior to \mcp{h}{\beta^{\#},\phi^{\#}}; we leave the proof details to Appendix \ref{singStratAppendix}.

\begin{prop}\label{singStrat}
Suppose \ess{F}{k} is a Fixer move at some history $\mc{h}\in\his{B}{k}$ and $\phi^{\#}$ is a Fixer strategy such that for any granular Buster strategy $\beta^1$ there exist a Fixer strategy $\phi^1$ and Buster strategy $\beta^{\#}$ for which $\phi^1(\mc{h})=\ess{F}{k}$ and \mcp{h}{\beta^1,\phi^1} is Fixer-superior to \mcp{h}{\beta^{\#},\phi^{\#}}. Then there exists a Fixer strategy $\phi^*$ satisfying $\phi^*(\mc{h})=\ess{F}{k}$ and Fixer-dominating $\phi^{\#}$ at \mc{h}.
\end{prop}

Combining Propositions \ref{fixDomStrat} and \ref{singStrat}, to show that a Fixer strategy $\phi$ that is Fixer-dominant at any history properly extending \mc{h} Fixer-dominates some given Fixer strategy $\phi'$ at \mc{h}, it would suffice to produce Fixer strategies $\phi^1$ and $\phi^{\#}$ such that $\phi^1(\mc{h})=\phi(\mc{h})$, $\phi^{\#}(\mc{h})=\phi'(\mc{h})$, $\phi^{\#}$ is Fixer-dominant at any history properly extending \mc{h}, and for any granular Buster strategy $\beta^1$ there exists a Buster strategy $\beta^{\#}$ for which \mcp{h}{\beta^1,\phi^1} is Fixer-superior to \mcp{h}{\beta^{\#},\phi^{\#}}. We shall employ this strategy to prove Proposition \ref{proofStrategy} in the next section.

Having just described how a granular approach is advantageous for Buster, we provide an analogue for Fixer by noting that Fixer is best off repairing graphs using only reserve edges that are bridges in the reconnected graph. Indeed, let $F\subseteq\ess{R}{k+1}$ be a subset of reserve edges that Fixer does not activate in the $k$th round. If the game ends before Fixer's turn in the $(k+1)$st round then under the same Buster play (or lack thereof) in the $(k+1)$st round the game would have also ended at the same point had Fixer included $F$ in her repairs, with Fixer receiving a higher payout by not having used the edges in $F$ in addition to those actually used in the $k$th round. Alternatively, if Fixer does have a turn in the $(k+1)$st round then Fixer can simply include $F$ in her response in that round. Hence we can state below that in order to show a Fixer strategy $\phi$ is Fixer-dominant at some history \mc{h}, it suffices to show that $\phi$ Fixer-dominates at \mc{h} all Fixer strategies $\phi'$ for which every edge of $\phi'(\mc{h})$ is a bridge in the reconnected graph, leaving the proof details to Appendix \ref{subsetSupPropAppendix}.

\begin{prop}\label{subsetSupProp}
Suppose $\mc{h}\in\his{B}{k}$ and $\phi$ is a Fixer strategy such that for every Fixer strategy $\phi'$ for which every edge of $\phi'(\mc{h})$ is a bridge in \es{G}{\mcp{h}{\phi'}}{k+1}, $\phi$ Fixer-dominates $\phi'$ at \mc{h}. Then $\phi$ is Fixer-dominant at \mc{h}.
\end{prop}

\section{Proof of Main Theorem} 
\label{proof}

To prove our main theorem, that any greedy Fixer strategy $\phi$ is Fixer-dominant at any history $\mc{h}\in\his{B}{k}$, we perform induction on $|\es{G}{\mc{h}}{k}|+|\es{R}{\mc{h}}{k}|$ (i.e. the total number of edges in the graph and reserve set); this will allow us to assume for our inductive hypothesis that any greedy Fixer strategy is Fixer-dominant at any proper extension of \mc{h}. Let $V$ be the vertex set of \es{G}{\mc{h}}{1}, with $|V|=n$, and without loss of generality assume $k=1$. Note that $|\es{G}{\mc{h}}{1}|\geq n-1$ (remembering that we let $|\es{G}{\mc{h}}{1}|$ denote the number of edges of \es{G}{\mc{h}}{1}) since \es{G}{\mc{h}}{1} is connected, so for our base case we consider $|\es{G}{\mc{h}}{1}|+|\es{R}{\mc{h}}{1}|=n-1$. In this case, $|\es{G}{\mc{h}}{1}|=n-1$ and $\es{R}{\mc{h}}{1}=\emptyset$, so $\emptyset\neq\es{B}{\mc{h}}{1}\subseteq\es{G}{\mc{h}}{1}$ implies Buster wins \mcp{\mc{h}}{\phi} in the first round because $(\es{G}{\mc{h}}{1}-\es{B}{\mc{h}}{1})\cup\es{R}{\mc{h}}{1}$ is disconnected due to it having at most $n-2$ edges. Hence there are no proper extensions of \mc{h}, so playing any Fixer strategy at \mc{h} involves no additional moves and leads to the payout of $\pay{\mc{h}}=-1$, and therefore all Fixer strategies are Fixer-dominant at \mc{h}. Thus we may suppose $|\es{G}{\mc{h}}{1}|+|\es{R}{\mc{h}}{1}|\geq n$ and inductively assume the following for the rest of the section, noting that for any proper extension $\mc{h'}\in\his{B}{k}$ of \mc{h} we have $|\es{G}{\mc{h'}}{k}|+|\es{R}{\mc{h'}}{k}|=|\es{G}{\mc{h'}}{k-1}|+|\es{R}{\mc{h'}}{k-1}|-|\es{B}{\mc{h'}}{k-1}|<|\es{G}{\mc{h}}{1}|+|\es{R}{\mc{h}}{1}|$.

\begin{hyp}\label{ih}
For any history $\mc{h'}\in\his{B}{k}$ such that \es{G}{\mc{h'}}{1} is a connected graph on $V$ and $|\es{G}{\mc{h'}}{k}|+|\es{R}{\mc{h'}}{k}|<|\es{G}{\mc{h}}{1}|+|\es{R}{\mc{h}}{1}|$, any greedy Fixer strategy is Fixer-dominant at \mc{h'}.
\end{hyp}

For the rest of this section, fix any greedy Fixer strategy $\phi$, any arbitrary Fixer strategy $\phi'$ such that every edge of $\phi'(\mc{h})$ is a bridge in $(\es{G}{\mc{h}}{1}-\es{B}{\mc{h}}{1})\cup\phi'(\mc{h})$, and any arbitrary Buster strategy $\beta$. By Proposition \ref{subsetSupProp}, in order to show $\phi$ is Fixer-dominant at \mc{h}, it suffices to show that $\phi$ Fixer-dominates $\phi'$ at \mc{h}; since $\beta$ is arbitrary, it further suffices to construct another Buster strategy $\beta'$ such that \mcp{h}{\beta,\phi} is Fixer-superior to \mcp{h}{\beta',\phi'}.

Let $M$ be the multigraph whose vertices are the components of $\es{G}{\mc{h}}{1}-\es{B}{\mc{h}}{1}$ and whose edges are the edges of \es{R}{\mc{h}}{1} (identifying each endpoint of the edges in \es{R}{\mc{h}}{1} with the component of $\es{G}{\mc{h}}{1}-\es{B}{\mc{h}}{1}$ within which it lies), so $\phi(\mc{h})$ is a minimum spanning tree of $M$, and $\phi'(\mc{h})$ is a spanning tree of $M$. We complete the proof by showing that for each number of components of $\es{G}{\mc{h}}{1}-\es{B}{\mc{h}}{1}$, $\phi$ Fixer-dominates $\phi'$ at \mc{h}, handling separately the cases of one component, two components, and at least three components in Subsections \ref{c1Proof}, \ref{c2Proof}, and \ref{c3Proof}, respectively. 

In Subsection \ref{c1Proof}, we observe that $\es{G}{\mc{h}}{1}-\es{B}{\mc{h}}{1}$ having a single component implies $\phi(\mc{h})=\phi'(\mc{h})=\emptyset$, allowing us to apply Inductive Hypothesis \ref{ih} to Lemma \ref{nextRound} (if $\phi$ and $\phi'$ are Fixer strategies satisfying $\phi(\mc{h})=\phi'(\mc{h})$, with $\phi$ greedy and therefore Fixer-dominant at any history properly extending \mc{h}, then $\phi$ Fixer-dominates $\phi'$ at \mc{h}) to conclude $\phi$ Fixer-dominates $\phi'$ at \mc{h}. 

In Subsection \ref{c2Proof}, $\es{G}{\mc{h}}{1}-\es{B}{\mc{h}}{1}$ having exactly two components dictates an approach based on combining Inductive Hypothesis \ref{ih} with Propositions \ref{fixDomStrat} and \ref{singStrat} (if $\phi^g$ is a Fixer strategy such that $\phi^g(\mc{h})=\phi'(\mc{h})$ and $\phi^g$ is greedy and therefore Fixer-dominant at any history properly extending \mc{h}, but nonetheless for every granular Buster strategy $\beta^1$ there exist a Fixer strategy $\phi^1$ and Buster strategy $\beta^g$ for which $\phi^1(\mc{h})=\phi(\mc{h})$ and \mcp{h}{\beta^1,\phi^1} is Fixer-superior to \mcp{h}{\beta^g,\phi^g}, then $\phi(\mc{h})$ must be a "better" move than $\phi'(\mc{h})$ for Fixer at \mc{h}, so $\phi$ being greedy and therefore Fixer-dominant at any proper extension of \mc{h} implies $\phi$ Fixer-dominates $\phi'$ at \mc{h}). Thus in order to show $\phi$ Fixer-dominates $\phi'$ at \mc{h}, we build the histories \mcp{h}{\beta^1,\phi^1} and \mcp{h}{\beta^g,\phi^g} described above so that \mcp{h}{\beta^1,\phi^1} is Fixer-superior to \mcp{h}{\beta^g,\phi^g}. We construct \mcp{h}{\beta^1,\phi^1} and \mcp{h}{\beta^g,\phi^g} round-by-round simultaneously, showing after each step that the collection of \es{G}{\mcp{h}{\beta^1,\phi^1}}{k}, \es{R}{\mcp{h}{\beta^1,\phi^1}}{k}, \es{G}{\mcp{h}{\beta^g,\phi^g}}{j}, and \es{R}{\mcp{h}{\beta^g,\phi^g}}{j} (where $k-1\leq j\leq k$) either satisfies the conditions of one of three given disjoint scenarios, or has Buster or Fixer win both games with \mcp{h}{\beta^1,\phi^1} Fixer-superior to \mcp{h}{\beta^g,\phi^g}. For each of the three scenarios, the graphs \es{G}{\mcp{h}{\beta^1,\phi^1}}{k} and \es{G}{\mcp{h}{\beta^g,\phi^g}}{j} are required to be similar or identical, as are the sets \es{R}{\mcp{h}{\beta^1,\phi^1}}{k} and \es{R}{\mcp{h}{\beta^g,\phi^g}}{j}, allowing \mcp{h}{\beta^1,\phi^1} to be Fixer-superior to \mcp{h}{\beta^g,\phi^g} once the games end.

In Subsection \ref{c3Proof}, we show that if $\es{G}{\mc{h}}{1}-\es{B}{\mc{h}}{1}$ has at least three components then we can construct a Fixer strategy $\phi''$ with the following properties. For particular nonempty edge sets $F\subseteq\phi(\mc{h})$ and $F'\subseteq\phi'(\mc{h})$, $\phi''$ will be defined so that $\phi''(\mc{h})=F\cup F'$ and $\phi''$ is greedy at all proper extensions of \mc{h}. Furthermore, a history $\mc{\widetilde{h}}\in\his{B}{1}$ will be constructed from \mc{h} by deleting a particular nonempty edge subset $B\subset\es{B}{\mc{h}}{1}$ from both \es{G}{\mc{h}}{1} and \es{B}{\mc{h}}{1}, shifting $F$ from \es{R}{\mc{h}}{1} to \es{G}{\mc{h}}{1}, and lowering Buster's limit by $|B|$; similarly, a history $\mc{\widetilde{h'}}\in\his{B}{1}$ will be constructed from \mc{h} by deleting a particular nonempty edge subset $B'\subset\es{B}{\mc{h}}{1}$ from both \es{G}{\mc{h}}{1} and \es{B}{\mc{h}}{1}, shifting $F'$ from \es{R}{\mc{h}}{1} to \es{G}{\mc{h}}{1}, and lowering Buster's limit by $|B'|$. By choosing $B$, $F$, $B'$, and $F'$ appropriately, we can define auxiliary Fixer strategies $\widetilde{\phi}$, $\widetilde{\phi'}$, and $\widetilde{\phi''}$ constructed, respectively, from the original Fixer strategies $\phi$, $\phi'$, and $\phi''$ by first setting $\widetilde{\phi}(\mc{\widetilde{h}})=\phi(\mc{h})-F$, $\widetilde{\phi''}(\mc{\widetilde{h}})=\phi''(\mc{h})-F$, $\widetilde{\phi''}(\mc{\widetilde{h'}})=\phi''(\mc{h})-F'$, and $\widetilde{\phi'}(\mc{\widetilde{h'}})=\phi'(\mc{h})-F'$ (so that after the first round the result of playing the original strategies at \mc{h} is the same as playing their corresponding auxiliary strategies at \mc{\widetilde{h}} and \mc{\widetilde{h'}}, i.e. $\es{G}{\mcp{h}{\phi}}{2}=\es{G}{\mcp{\widetilde{h}}{\widetilde{\phi}}}{2}$, $\es{R}{\mcp{h}{\phi}}{2}=\es{R}{\mcp{\widetilde{h}}{\widetilde{\phi}}}{2}$, $\es{G}{\mcp{h}{\phi''}}{2}=\es{G}{\mcp{\widetilde{h}}{\widetilde{\phi''}}}{2}=\es{G}{\mcp{\widetilde{h'}}{\widetilde{\phi''}}}{2}$, $\es{R}{\mcp{h}{\phi''}}{2}=\es{R}{\mcp{\widetilde{h}}{\widetilde{\phi''}}}{2}=\es{R}{\mcp{\widetilde{h'}}{\widetilde{\phi''}}}{2}$, $\es{G}{\mcp{h}{\phi'}}{2}=\es{G}{\mcp{\widetilde{h'}}{\widetilde{\phi'}}}{2}$, and $\es{R}{\mcp{h}{\phi'}}{2}=\es{R}{\mcp{\widetilde{h'}}{\widetilde{\phi'}}}{2}$) and then directly translating the Fixer moves that $\phi$, $\phi'$, and $\phi''$ prescribe at any given proper extension \mc{h^*} of \mc{h} to the proper extensions of \mc{\widetilde{h}} and \mc{\widetilde{h'}} that differ only from \mc{h^*} in the first round. Hence $\widetilde{\phi}$ is greedy at \mc{\widetilde{h}} and all its extensions while $\widetilde{\phi''}$ is greedy at \mc{\widetilde{h'}} and all its extensions, so noting $|\es{G}{\mc{\widetilde{h}}}{1}|+|\es{R}{\mc{\widetilde{h}}}{1}|=|\es{G}{\mc{h}}{1}|+|\es{R}{\mc{h}}{1}|-|B|$ and $|\es{G}{\mc{\widetilde{h'}}}{1}|+|\es{R}{\mc{\widetilde{h'}}}{1}|=|\es{G}{\mc{h}}{1}|+|\es{R}{\mc{h}}{1}|-|B'|$ and applying Inductive Hypothesis \ref{ih} we observe that $\widetilde{\phi}$ is Fixer-dominant at \mc{\widetilde{h}} while $\widetilde{\phi''}$ is Fixer-dominant at \mc{\widetilde{h'}}, in which case $\widetilde{\phi}$ Fixer-dominates $\widetilde{\phi''}$ at \mc{\widetilde{h}} while $\widetilde{\phi''}$ Fixer-dominates $\widetilde{\phi'}$ at \mc{\widetilde{h'}}, and thus by the definitions of $\widetilde{\phi}$, $\widetilde{\phi'}$, and $\widetilde{\phi''}$ relative to $\phi$, $\phi'$, and $\phi''$ we see that $\phi$ Fixer-dominates $\phi''$ at \mc{h} while $\phi''$ Fixer-dominates $\phi'$ at \mc{h}, allowing us to conclude by Proposition \ref{fixDomTransitivityCor} that $\phi$ Fixer-dominates $\phi'$ at \mc{h}.

\subsection{The case that $\es{G}{\mc{h}}{1}-\es{B}{\mc{h}}{1}$ has one component}
\label{c1Proof}

If $\es{G}{\mc{h}}{1}-\es{B}{\mc{h}}{1}$ is connected, then the only spanning tree of our multigraph $M$ is edgeless, so $\phi(\mc{h})=\phi'(\mc{h})=\emptyset$. Since $\es{B}{\mc{h}}{1}\neq\emptyset$, by Inductive Hypothesis \ref{ih} $\phi$ is Fixer-dominant at any proper history of \mc{h}, so by Lemma \ref{nextRound} $\phi$ Fixer-dominates $\phi'$ at \mc{h}. Hence $\phi$ is Fixer-dominant at \mc{h}.

\subsection{The case that $\es{G}{\mc{h}}{1}-\es{B}{\mc{h}}{1}$ has two components}
\label{c2Proof}

If $\es{G}{\mc{h}}{1}-\es{B}{\mc{h}}{1}$ has two components, then the spanning trees of our multigraph $M$ are the individual edges in \es{R}{\mc{h}}{1} joining the two components of $\es{G}{\mc{h}}{1}-\es{B}{\mc{h}}{1}$. Hence for two such edges $s$ and $s'$, where no such edge is cheaper than $s$, we have $\phi(\mc{h})=\{s\}$ and $\phi'(\mc{h})=\{s'\}$; if $s=s'$ then by Lemma \ref{nextRound} $\phi$ Fixer-dominates $\phi'$ at \mc{h} because by Inductive Hypothesis \ref{ih} $\phi$ is Fixer-dominant at any proper history of \mc{h}; thus we may assume $s\neq s'$. We establish Proposition \ref{proofStrategy}, which provides a strategy for proving that $\phi$ is Fixer-dominant at \mc{h}, then show how to execute that strategy in Proposition \ref{proofExecution}, allowing us to conclude in Corollary \ref{conclusionCorollary2} that $\phi$ Fixer-dominates $\phi'$ at \mc{h}. Details of the proof of Proposition \ref{proofExecution} are extracted to Appendix \ref{proofExecutionAppendix}.

Recall that a Buster strategy is granular if it restricts Buster to removing only singletons. Call a Fixer strategy $\phi^g$ \emph{almost greedy} if $\phi^g(\mc{h})=\phi'(\mc{h})=\{s'\}$ and $\phi^g$ is greedy for all other moves, with $s\in\phi^g(\mc{h'})$ for any history $\mc{h'}\in\his{B}{k}$ for which $k>1$ and $s$ is part of some greedy Fixer response to Buster's move in the $k$th round of \mc{h'}.

\begin{prop}\label{proofStrategy}
Suppose there exists an almost greedy Fixer strategy $\phi^g$ such that for every granular Buster strategy $\beta^1$ there exist a Fixer strategy $\phi^1$ and Buster strategy $\beta^g$ for which $\phi^1(\mc{h})=\phi(\mc{h})=\{s\}$ and \mcp{h}{\beta^1,\phi^1} is Fixer-superior to \mcp{h}{\beta^g,\phi^g}. Then $\phi$ Fixer-dominates $\phi'$ at \mc{h}.
\end{prop}
\begin{proof}
By Proposition \ref{singStrat} there exists a Fixer strategy $\phi^*$ satisfying $\phi^*(\mc{h})=\phi^1(\mc{h})=\{s\}$ and Fixer-dominating $\phi^g$ at \mc{h}. By Inductive Hypothesis \ref{ih}, since $\phi$ and $\phi^g$ are greedy at any history properly extending \mc{h}, we have $\phi$ and $\phi^g$ Fixer-dominant at any history properly extending \mc{h}. Since $\phi^*$ is a Fixer strategy satisfying $\phi^*(\mc{h})=\phi(\mc{h})$ and Fixer-dominating $\phi^g$ at \mc{h}, and $\phi'(\mc{h})=\phi^g(\mc{h})$, by Proposition \ref{fixDomStrat} $\phi$ Fixer-dominates $\phi'$ at \mc{h}.
\end{proof}

By Proposition \ref{proofStrategy}, in order to show $\phi$ Fixer-dominates $\phi'$ at \mc{h}, we are to start with some almost greedy Fixer strategy $\phi^g$ and arbitrary granular Buster strategy $\beta^1$, then construct a Fixer strategy $\phi^1$ and Buster strategy $\beta^g$ for which $\phi^1(\mc{h})=\phi(\mc{h})=\{s\}$ and \mcp{h}{\beta^1,\phi^1} is Fixer-superior to \mcp{h}{\beta^g,\phi^g}. We shall define $\phi^1$ and $\beta^g$ by playing \mcp{h}{\beta^1,\phi^1} and \mcp{h}{\beta^g,\phi^g} round-by-round simultaneously, with $\beta^g$ being granular. For positive integers $k$ and $j$, let \pair{\mcp{h}{\beta^1,\phi^1}}{k}{\mcp{h}{\beta^g,\phi^g}}{j} denote the pair of histories consisting of the first $k-1$ rounds of \mcp{h}{\beta^1,\phi^1} (awaiting Buster's $k$th round move) and the first $j-1$ rounds of \mcp{h}{\beta^g,\phi^g} (awaiting Buster's $j$th round move). We show that after each step of our side-by-side play that \pair{\mcp{h}{\beta^1,\phi^1}}{k}{\mcp{h}{\beta^g,\phi^g}}{j} (where $j$ equals either $k$ or $k-1$) either satisfies the conditions of one of Scenarios \ref{sce1}, \ref{sce2}, and \ref{sce3}, each detailed below, or has Buster or Fixer win both games with \mcp{h}{\beta^1,\phi^1} Fixer-superior to \mcp{h}{\beta^g,\phi^g}. Roughly speaking, \pair{\mcp{h}{\beta^1,\phi^1}}{k}{\mcp{h}{\beta^g,\phi^g}}{k} belongs to Scenario \ref{sce1} if the only difference between \mcp{h}{\beta^1,\phi^1} and \mcp{h}{\beta^g,\phi^g} is that in the former, $s$ is still in the graph and $s'$ is still in the reserve set (i.e. $s\in\es{G}{\mcp{h}{\beta^1,\phi^1}}{k}$ and $s'\in\es{R}{\mcp{h}{\beta^1,\phi^1}}{k}$) while in the latter they are switched (i.e. $s'\in\es{G}{\mcp{h}{\beta^g,\phi^g}}{k}$ and $s\in\es{R}{\mcp{h}{\beta^g,\phi^g}}{k}$); \pair{\mcp{h}{\beta^1,\phi^1}}{k}{\mcp{h}{\beta^g,\phi^g}}{k-1} belongs to Scenario \ref{sce2} if \mcp{h}{\beta^1,\phi^1} is a round ahead of \mcp{h}{\beta^g,\phi^g} and no longer contains $s$ in its graph, with the only difference between the two games (besides the number of rounds played) being \mcp{h}{\beta^g,\phi^g} having $s'$ in its reserve set (i.e. $s\notin\es{G}{\mcp{h}{\beta^1,\phi^1}}{k}=\es{G}{\mcp{h}{\beta^g,\phi^g}}{k-1}$ and $s\in\es{R}{\mcp{h}{\beta^g,\phi^g}}{k-1}$ with $\es{R}{\mcp{h}{\beta^1,\phi^1}}{k}=\es{R}{\mcp{h}{\beta^g,\phi^g}}{k-1}-\{s\}$); and \pair{\mcp{h}{\beta^1,\phi^1}}{k}{\mcp{h}{\beta^g,\phi^g}}{k} belongs to Scenario \ref{sce3} if the two games are identical (i.e. $\es{G}{\mcp{h}{\beta^1,\phi^1}}{k}=\es{G}{\mcp{h}{\beta^g,\phi^g}}{k}$ and $\es{R}{\mcp{h}{\beta^1,\phi^1}}{k}=\es{R}{\mcp{h}{\beta^g,\phi^g}}{k}$). Since each scenario keeps the graphs \es{G}{\mcp{h}{\beta^1,\phi^1}}{k} and \es{G}{\mcp{h}{\beta^g,\phi^g}}{j} as well as the reserve sets \es{R}{\mcp{h}{\beta^1,\phi^1}}{k} and \es{R}{\mcp{h}{\beta^g,\phi^g}}{j} similar or identical, we can show \mcp{h}{\beta^1,\phi^1} is Fixer-superior to \mcp{h}{\beta^g,\phi^g} once the games end.

\begin{sce}\label{sce1}
Scenario where Buster has not used $s$ in \mcp{h}{\beta^1,\phi^1} and Fixer has not used $s$ in \mcp{h}{\beta^g,\phi^g}: \pair{\mcp{h}{\beta^1,\phi^1}}{k}{\mcp{h}{\beta^g,\phi^g}}{k} belongs to this scenario if \mcp{h}{\beta^1,\phi^1} and \mcp{h}{\beta^g,\phi^g} each start the $k$th round with the following properties:
\begin{itemize}
\item $s\in\es{G}{\mcp{h}{\beta^1,\phi^1}}{k}$, $s'\in\es{G}{\mcp{h}{\beta^g,\phi^g}}{k}$, and $\es{G}{\mcp{h}{\beta^1,\phi^1}}{k}-\{s\}=\es{G}{\mcp{h}{\beta^g,\phi^g}}{k}-\{s'\}$ (i.e. \es{G}{\mcp{h}{\beta^g,\phi^g}}{k} is \es{G}{\mcp{h}{\beta^1,\phi^1}}{k} with $s$ replaced by $s'$)
\item $s'\in\es{R}{\mcp{h}{\beta^1,\phi^1}}{k}$, $s\in\es{R}{\mcp{h}{\beta^g,\phi^g}}{k}$, and $\es{R}{\mcp{h}{\beta^1,\phi^1}}{k}-\{s'\}=\es{R}{\mcp{h}{\beta^g,\phi^g}}{k}-\{s\}$ (i.e. \es{R}{\mcp{h}{\beta^g,\phi^g}}{k} is \es{R}{\mcp{h}{\beta^1,\phi^1}}{k} with $s'$ replaced by $s$)
\item $s$ and $s'$ are bridges in \es{G}{\mcp{h}{\beta^1,\phi^1}}{k} and \es{G}{\mcp{h}{\beta^g,\phi^g}}{k}, respectively, between the same subgraphs $X_k$ and $Y_k$, but perhaps in different spots (i.e. removing both edges from their graphs leaves the same two-component graphs); see Figure \ref{figsGT1GTgXkYk}
\item for every $r\in\es{R}{\mcp{h}{\beta^1,\phi^1}}{k}\cup\es{R}{\mcp{h}{\beta^g,\phi^g}}{k}$ such that $r$ joins $X_k$ to $Y_k$, $w(r)\geq w(s)$ (i.e. in either game, no reserve edge going between subgraphs $X_k$ and $Y_k$ can be cheaper than $s$)
\end{itemize}
\end{sce}

\begin{figure}[htb]
\centering
\subcaptionbox{\es{G}{\mcp{h}{\beta^1,\phi^1}}{k} \label{figGT1XkYk}}[9cm]
{
\begin{tikzpicture}
\Vertex[x=0,y=0,size=1.5,label=$X_k$,fontsize=\large]{xk}
\Vertex[x=2,y=0,size=1.5,label=$Y_k$,fontsize=\large]{yk}
\Edge[label=$s$,bend=45,position={below=.5mm},fontsize=\large](xk)(yk)
\end{tikzpicture}
}
\subcaptionbox{\es{G}{\mcp{h}{\beta^g,\phi^g}}{k} \label{figGTgXkYk}}[9cm]
{
\begin{tikzpicture}
\Vertex[x=0,y=0,size=1.5,label=$X_k$,fontsize=\large]{xk}
\Vertex[x=2,y=0,size=1.5,label=$Y_k$,fontsize=\large]{yk}
\Edge[label=$s'$,bend=-45,position={above=.5mm},fontsize=\large](xk)(yk)
\end{tikzpicture}
}
\caption{Graphs \es{G}{\mcp{h}{\beta^1,\phi^1}}{k} and \es{G}{\mcp{h}{\beta^g,\phi^g}}{k} from Scenario \ref{sce1}.}\label{figsGT1GTgXkYk}
\end{figure}

\begin{sce}\label{sce2}
Scenario where Buster has used $s$ in \mcp{h}{\beta^1,\phi^1} but Fixer has not used $s$ in \mcp{h}{\beta^g,\phi^g}, which is a round behind \mcp{h}{\beta^1,\phi^1}: \pair{\mcp{h}{\beta^1,\phi^1}}{k}{\mcp{h}{\beta^g,\phi^g}}{k-1} belongs to this scenario if \mcp{h}{\beta^1,\phi^1} starts the $k$th round and \mcp{h}{\beta^g,\phi^g} starts the $(k-1)$st round with the following properties:
\begin{itemize}
\item $s'\in\es{G}{\mcp{h}{\beta^1,\phi^1}}{k}=\es{G}{\mcp{h}{\beta^g,\phi^g}}{k-1}$ (i.e. the graphs are identical and contain $s'$)
\item $\es{R}{\mcp{h}{\beta^1,\phi^1}}{k}=\es{R}{\mcp{h}{\beta^g,\phi^g}}{k-1}-\{ s\}$ and $s\in\es{R}{\mcp{h}{\beta^g,\phi^g}}{k-1}$ (i.e. the reserve sets only differ by $s$ being in \es{R}{\mcp{h}{\beta^g,\phi^g}}{k-1} but not \es{R}{\mcp{h}{\beta^1,\phi^1}}{k})
\item $s'$ and $s$ are bridges in \es{G}{\mcp{h}{\beta^1,\phi^1}}{k} and $(\es{G}{\mcp{h}{\beta^g,\phi^g}}{k-1}-\{s'\})\cup\{s\}$, respectively, between the same connected subgraphs $X_k$ and $Y_k$, but perhaps in different spots (i.e. removing both edges from their respective graphs leaves the same graphs, each with two components); see Figure \ref{scenario2Figs}
\item for every $r\in\es{R}{\mcp{h}{\beta^g,\phi^g}}{k-1}$ such that $r$ bridges $X_k$ and $Y_k$, $w(r)\geq w(s)$ (i.e. in either game, no reserve edge bridging subgraphs $X_k$ and $Y_k$ can be cheaper than $s$)
\end{itemize}
\end{sce}

\begin{figure}[htb]
\centering
\subcaptionbox{$\es{G}{\mcp{h}{\beta^1,\phi^1}}{k}=\es{G}{\mcp{h}{\beta^g,\phi^g}}{k-1}$ \label{scenario2Fig1}}[9cm]
{
\begin{tikzpicture}
\Vertex[x=0,y=0,size=1.5,label=$X_k$,fontsize=\large]{xk}
\Vertex[x=2,y=0,size=1.5,label=$Y_k$,fontsize=\large]{yk}
\Edge[label=$s'$,bend=-45,position={above=.5mm},fontsize=\large](xk)(yk)
\end{tikzpicture}
}
\subcaptionbox{$(\es{G}{\mcp{h}{\beta^g,\phi^g}}{k-1}-\{s'\})\cup\{s\}$ \label{scenario2Fig2}}[9cm]
{
\begin{tikzpicture}
\Vertex[x=0,y=0,size=1.5,label=$X_k$,fontsize=\large]{xk}
\Vertex[x=2,y=0,size=1.5,label=$Y_k$,fontsize=\large]{yk}
\Edge[label=$s$,bend=45,position={below=.5mm},fontsize=\large](xk)(yk)
\end{tikzpicture}
}
\caption{Graphs $\es{G}{\mcp{h}{\beta^1,\phi^1}}{k}=\es{G}{\mcp{h}{\beta^g,\phi^g}}{k-1}$ and $(\es{G}{\mcp{h}{\beta^g,\phi^g}}{k-1}-\{s'\})\cup\{s\}$ from Scenario \ref{sce2}.}\label{scenario2Figs}
\end{figure}

\begin{sce}\label{sce3}
Scenario where \mcp{h}{\beta^1,\phi^1} and \mcp{h}{\beta^g,\phi^g} are in the same state: \pair{\mcp{h}{\beta^1,\phi^1}}{k}{\mcp{h}{\beta^g,\phi^g}}{k} belongs to this scenario if \mcp{h}{\beta^1,\phi^1} and \mcp{h}{\beta^g,\phi^g} each start the $k$th round with the following properties:
\begin{itemize}
\item $\es{G}{\mcp{h}{\beta^1,\phi^1}}{k}=\es{G}{\mcp{h}{\beta^g,\phi^g}}{k}$
\item $\es{R}{\mcp{h}{\beta^1,\phi^1}}{k}=\es{R}{\mcp{h}{\beta^g,\phi^g}}{k}$
\end{itemize}
\end{sce}

\begin{prop}\label{proofExecution}
For every almost greedy Fixer strategy $\phi^g$ and granular Buster strategy $\beta^1$, there exist a Fixer strategy $\phi^1$ and granular Buster strategy $\beta^g$ for which $\phi^1(\mc{h})=\phi(\mc{h})=\{s\}$ and \mcp{h}{\beta^1,\phi^1} is Fixer-superior to \mcp{h}{\beta^g,\phi^g}.
\end{prop}
\begin{proof}
Fix any almost greedy Fixer strategy $\phi^g$ and granular Buster strategy $\beta^1$. We define the part of a Fixer strategy $\phi^1$ acting against $\beta^1$ by first setting $\phi^1(\mc{h})=\phi(\mc{h})=\{s\}$ and then giving its moves for Fixer against $\beta^1$, with these moves being derived from those prescribed by $\phi^g$ against a Buster strategy $\beta^g$, with $\beta^g$ derived from $\beta^1$; these derivations are described in the subsequent paragraphs. Note that \pair{\mcp{h}{\beta^1,\phi^1}}{2}{\mcp{h}{\beta^g,\phi^g}}{2} denotes the pair of histories $((\es{G}{\mcp{h}{\beta^1,\phi^1}}{1}=\es{G}{\mc{h}}{1},\es{R}{\mcp{h}{\beta^1,\phi^1}}{1}=\es{R}{\mc{h}}{1}),\es{B}{\mcp{h}{\beta^1,\phi^1}}{1}=\es{B}{\mc{h}}{1},\es{F}{\mcp{h}{\beta^1,\phi^1}}{1}=\phi(\mc{h})=\{s\},(\es{G}{\mcp{h}{\beta^1,\phi^1}}{2}=(\es{G}{\mc{h}}{1}-\es{B}{\mc{h}}{1})\cup\{s\},\es{R}{\mcp{h}{\beta^1,\phi^1}}{2}=\es{R}{\mc{h}}{1}-\{s\}))$ and $((\es{G}{\mcp{h}{\beta^g,\phi^g}}{1}=\es{G}{\mc{h}}{1},\es{R}{\mcp{h}{\beta^g,\phi^g}}{1}=\es{R}{\mc{h}}{1}),\es{B}{\mcp{h}{\beta^g,\phi^g}}{1}=\es{B}{\mc{h}}{1},\es{F}{\mcp{h}{\beta^g,\phi^g}}{1}=\phi'(\mc{h})=\{s'\},(\es{G}{\mcp{h}{\beta^g,\phi^g}}{2}=(\es{G}{\mc{h}}{1}-\es{B}{\mc{h}}{1})\cup\{s'\},\es{R}{\mcp{h}{\beta^g,\phi^g}}{2}=\es{R}{\mc{h}}{1}-\{s'\}))$. Starting with $k=j=2$, we define $\phi^1$ and $\beta^g$ by extending \pair{\mcp{h}{\beta^1,\phi^1}}{k}{\mcp{h}{\beta^g,\phi^g}}{j} by increasing $k$ and/or $j$ in each of a series of actions until Buster of Fixer wins the game.

In order to analyze \mcp{h}{\beta^1,\phi^1} and \mcp{h}{\beta^g,\phi^g}, we categorize \pair{\mcp{h}{\beta^1,\phi^1}}{k}{\mcp{h}{\beta^g,\phi^g}}{j} into Scenarios \ref{sce1}, \ref{sce2}, and \ref{sce3}, each fully defined above. To complete the proof, we note via Proposition \ref{scenario1prop0} that \pair{\mcp{h}{\beta^1,\phi^1}}{2}{\mcp{h}{\beta^g,\phi^g}}{2} belongs to Scenario \ref{sce1}, and then show via Propositions \ref{scenario1prop1} through \ref{scenario3prop2} that each extension of \pair{\mcp{h}{\beta^1,\phi^1}}{k}{\mcp{h}{\beta^g,\phi^g}}{j} either lands in Scenario \ref{sce1}, \ref{sce2}, or \ref{sce3}, or results in \mcp{h}{\beta^1,\phi^1} being Fixer-superior to \mcp{h}{\beta^g,\phi^g} after Fixer or Buster wins both games; see Figure \ref{strategyFigure} for a flowchart and Appendix \ref{proofExecutionAppendix} for statements and proofs of Propositions \ref{scenario1prop0} through \ref{scenario3prop2}. Since the finiteness of \es{G}{\mc{h}}{1}, \es{R}{\mc{h}}{1}, and the limit $\ell$ implies that every game ends with either Fixer or Buster winning, this shows that \mcp{h}{\beta^1,\phi^1} is Fixer-superior to \mcp{h}{\beta^g,\phi^g}.

First suppose \pair{\mcp{h}{\beta^1,\phi^1}}{k}{\mcp{h}{\beta^g,\phi^g}}{k} belongs to Scenario \ref{sce1} for some $k\geq 2$. If Fixer wins \mcp{h}{\beta^1,\phi^1} in the $k$th round, then Fixer also wins \mcp{h}{\beta^g,\phi^g} in the $k$th round while leaving \mcp{h}{\beta^1,\phi^1} Fixer-superior to \mcp{h}{\beta^g,\phi^g} (see \ref{111} of Proposition \ref{scenario1prop1} for details). Otherwise, Buster makes a move \es{B}{\mcp{h}{\beta^1,\phi^1}}{k} in the $k$th round of \mcp{h}{\beta^1,\phi^1} according to $\beta^1$. If Buster wins \mcp{h}{\beta^1,\phi^1} in the $k$th round, define $\beta^g$ so that $\es{B}{\mcp{h}{\beta^g,\phi^g}}{k}=\es{B}{\mcp{h}{\beta^1,\phi^1}}{k}$, and then Buster also wins \mcp{h}{\beta^g,\phi^g} while leaving \mcp{h}{\beta^1,\phi^1} Fixer-superior to \mcp{h}{\beta^g,\phi^g} (see \ref{112} of Proposition \ref{scenario1prop1} for details). Now suppose neither Fixer nor Buster wins \mcp{h}{\beta^1,\phi^1} in the $k$th round. If $\es{B}{\mcp{h}{\beta^1,\phi^1}}{k}=\{s\}$, let $\phi^1$ have Fixer play $\es{F}{\mcp{h}{\beta^1,\phi^1}}{k}=\{s'\}$ so that \pair{\mcp{h}{\beta^1,\phi^1}}{k+1}{\mcp{h}{\beta^g,\phi^g}}{k} belongs to Scenario \ref{sce2} (see Proposition \ref{scenario1prop2} for details). Otherwise, define $\beta^g$ so that $\es{B}{\mcp{h}{\beta^g,\phi^g}}{k}=\es{B}{\mcp{h}{\beta^1,\phi^1}}{k}$, and let Fixer respond with a move \es{F}{\mcp{h}{\beta^g,\phi^g}}{k} in the $k$th round of \mcp{h}{\beta^g,\phi^g} according to $\phi^g$. If $\es{F}{\mcp{h}{\beta^g,\phi^g}}{k}\neq\{s\}$, define $\phi^1$ so that $\es{F}{\mcp{h}{\beta^1,\phi^1}}{k}=\es{F}{\mcp{h}{\beta^g,\phi^g}}{k}$ to leave \pair{\mcp{h}{\beta^1,\phi^1}}{k+1}{\mcp{h}{\beta^g,\phi^g}}{k+1} belonging to Scenario \ref{sce1} (see Proposition \ref{scenario1prop3} for details, with \ref{131} covering $\es{G}{\mcp{h}{\beta^1,\phi^1}}{k}-\es{B}{\mcp{h}{\beta^1,\phi^1}}{k}$ being connected, and \ref{132} and \ref{133} covering $\es{G}{\mcp{h}{\beta^1,\phi^1}}{k}-\es{B}{\mcp{h}{\beta^1,\phi^1}}{k}$ being disconnected, each dealing with a separate possibility of the position of $s'$ in \es{G}{\mcp{h}{\beta^g,\phi^g}}{k}). If $\es{F}{\mcp{h}{\beta^g,\phi^g}}{k}=\{s\}$, let $\phi^1$ have Fixer play $\es{F}{\mcp{h}{\beta^1,\phi^1}}{k}=\{s'\}$ so that \pair{\mcp{h}{\beta^1,\phi^1}}{k+1}{\mcp{h}{\beta^g,\phi^g}}{k+1} belongs to Scenario \ref{sce3} (see \ref{134} of Proposition \ref{scenario1prop3} for details).

Next suppose \pair{\mcp{h}{\beta^1,\phi^1}}{k}{\mcp{h}{\beta^g,\phi^g}}{k-1} belongs to Scenario \ref{sce2} for some $k\geq 3$. If Fixer wins \mcp{h}{\beta^1,\phi^1} in the $k$th round, define $\beta^g$ to have Buster play $\es{B}{\mcp{h}{\beta^g,\phi^g}}{k-1}=\{s'\}$ in \mcp{h}{\beta^g,\phi^g}, against which $\phi^g$ will have Fixer play $\es{F}{\mcp{h}{\beta^g,\phi^g}}{k-1}=\{s\}$ as a greedy response, in which case Fixer wins \mcp{h}{\beta^g,\phi^g} in the $k$th round with \mcp{h}{\beta^1,\phi^1} Fixer-superior to \mcp{h}{\beta^g,\phi^g} (see \ref{211} of Proposition \ref{scenario2prop1} for details). Otherwise, Buster makes a move \es{B}{\mcp{h}{\beta^1,\phi^1}}{k} in the $k$th round of \mcp{h}{\beta^1,\phi^1} according to $\beta^1$. If Buster wins $\mcp{h}{\beta^1,\phi^1}$ in the $k$th round such that $(\es{G}{\mcp{h}{\beta^1,\phi^1}}{k}-\es{B}{\mcp{h}{\beta^1,\phi^1}}{k})\cup\es{R}{\mcp{h}{\beta^1,\phi^1}}{k}\cup\{s\}$ is disconnected, define $\beta^g$ to have Buster play $\es{B}{\mcp{h}{\beta^g,\phi^g}}{k-1}=\{s'\}$ in \mcp{h}{\beta^g,\phi^g}, against which $\phi^g$ will have Fixer respond greedily with $\es{F}{\mcp{h}{\beta^g,\phi^g}}{k-1}=\{s\}$. Further define $\beta^g$ to have Buster play $\es{B}{\mcp{h}{\beta^g,\phi^g}}{k}=\es{B}{\mcp{h}{\beta^1,\phi^1}}{k}$ in \mcp{h}{\beta^g,\phi^g} if $\es{B}{\mcp{h}{\beta^1,\phi^1}}{k}\neq\{s'\}$, or play $\es{B}{\mcp{h}{\beta^g,\phi^g}}{k}=\{s\}$ in \mcp{h}{\beta^g,\phi^g} if $\es{B}{\mcp{h}{\beta^1,\phi^1}}{k}=\{s'\}$, both resulting in Buster winning $\mcp{h}{\beta^g,\phi^g}$ in the $k$th round and \mcp{h}{\beta^1,\phi^1} Fixer-superior to \mcp{h}{\beta^g,\phi^g} (see \ref{212} of Proposition \ref{scenario2prop1} for details). If Buster wins $\mcp{h}{\beta^1,\phi^1}$ in the $k$th round such that $(\es{G}{\mcp{h}{\beta^1,\phi^1}}{k}-\es{B}{\mcp{h}{\beta^1,\phi^1}}{k})\cup\es{R}{\mcp{h}{\beta^1,\phi^1}}{k}\cup\{s\}$ is connected, then define $\beta^g$ to have Buster play $\es{B}{\mcp{h}{\beta^g,\phi^g}}{k-1}=\es{B}{\mcp{h}{\beta^1,\phi^1}}{k}$ in \mcp{h}{\beta^g,\phi^g}, against which $\phi^g$ will have Fixer play $\es{F}{\mcp{h}{\beta^g,\phi^g}}{k-1}=\{s\}$ as a greedy response, which allows us to define $\beta^g$ to have Buster play $\es{B}{\mcp{h}{\beta^g,\phi^g}}{k}=\{s\}$ in \mcp{h}{\beta^g,\phi^g}, resulting in Buster winning $\mcp{h}{\beta^g,\phi^g}$ in the $k$th round and \mcp{h}{\beta^1,\phi^1} Fixer-superior to \mcp{h}{\beta^g,\phi^g} (see \ref{213} of Proposition \ref{scenario2prop1} for details). Now suppose neither Fixer nor Buster wins $\mcp{h}{\beta^1,\phi^1}$ in the $k$th round, so $|\mcp{h}{\beta^1,\phi^1}|>k$. Define $\beta^g$ to have Buster copy his move from \mcp{h}{\beta^1,\phi^1} into \mcp{h}{\beta^g,\phi^g} with $\es{B}{\mcp{h}{\beta^g,\phi^g}}{k-1}=\es{B}{\mcp{h}{\beta^1,\phi^1}}{k}$, and let Fixer respond with a move \es{F}{\mcp{h}{\beta^g,\phi^g}}{k-1} in the $(k-1)$st round of \mcp{h}{\beta^g,\phi^g} according to $\phi^g$. If $\es{F}{\mcp{h}{\beta^g,\phi^g}}{k-1}\neq\{s\}$, define $\phi^1$ to have Fixer copy that move into \mcp{h}{\beta^1,\phi^1} with $\es{F}{\mcp{h}{\beta^1,\phi^1}}{k}=\es{F}{\mcp{h}{\beta^g,\phi^g}}{k-1}$ so that \pair{\mcp{h}{\beta^1,\phi^1}}{k+1}{\mcp{h}{\beta^g,\phi^g}}{k} belongs to Scenario \ref{sce2} (see \ref{221} of Proposition \ref{scenario2prop2} for details). If $\es{F}{\mcp{h}{\beta^g,\phi^g}}{k-1}=\{s\}$, define $\beta^g$ to have Buster play $\es{B}{\mcp{h}{\beta^g,\phi^g}}{k}=\{s\}$, against which $\phi^g$ will have Fixer respond by creating a connected graph \es{G}{\mcp{h}{\beta^g,\phi^g}}{k+1} with some greedy \es{F}{\mcp{h}{\beta^g,\phi^g}}{k} in \mcp{h}{\beta^g,\phi^g}; then we can define $\phi^1$ to have Fixer play $\es{F}{\mcp{h}{\beta^1,\phi^1}}{k}=\es{F}{\mcp{h}{\beta^g,\phi^g}}{k}$ in \mcp{h}{\beta^1,\phi^1}, so that \pair{\mcp{h}{\beta^1,\phi^1}}{k+1}{\mcp{h}{\beta^g,\phi^g}}{k+1} belongs to Scenario \ref{sce3} (see \ref{222} of Proposition \ref{scenario2prop2} for details).

Finally suppose \pair{\mcp{h}{\beta^1,\phi^1}}{k}{\mcp{h}{\beta^g,\phi^g}}{k} belongs to Scenario \ref{sce3} for some $k\geq 3$. If Fixer wins \mcp{h}{\beta^1,\phi^1} in the $k$th round, then Fixer also wins \mcp{h}{\beta^g,\phi^g} in the $k$th round leaving \mcp{h}{\beta^1,\phi^1} Fixer-superior to \mcp{h}{\beta^g,\phi^g} (see \ref{311} of Proposition \ref{scenario3prop1} for details). Otherwise Buster makes a move \es{B}{\mcp{h}{\beta^1,\phi^1}}{k} in the $k$th round of \mcp{h}{\beta^1,\phi^1} according to $\beta^1$. Define $\beta^g$ so that $\es{B}{\mcp{h}{\beta^g,\phi^g}}{k}=\es{B}{\mcp{h}{\beta^1,\phi^1}}{k}$. If Buster wins \mcp{h}{\beta^1,\phi^1} in the $k$th round, then Buster also wins \mcp{h}{\beta^g,\phi^g} in the $k$th round leaving \mcp{h}{\beta^1,\phi^1} Fixer-superior to \mcp{h}{\beta^g,\phi^g} (see \ref{312} of Proposition \ref{scenario3prop1} for details). If neither Fixer nor Buster wins \mcp{h}{\beta^1,\phi^1} in the $k$th round, then Fixer responds with a move \es{F}{\mcp{h}{\beta^g,\phi^g}}{k} in the $k$th round of \mcp{h}{\beta^g,\phi^g} according to $\phi^g$. Define $\phi^1$ so that $\es{F}{\mcp{h}{\beta^1,\phi^1}}{k}=\es{F}{\mcp{h}{\beta^g,\phi^g}}{k}$, and then \pair{\mcp{h}{\beta^1,\phi^1}}{k+1}{\mcp{h}{\beta^g,\phi^g}}{k+1} belongs to Scenario \ref{sce3} (see Proposition \ref{scenario3prop2} for details).
\end{proof}

\begin{figure}[htb]
\centering
\begin{tikzpicture}
\Vertex[x=-8,y=0,size=2,label=Round 1,fontsize=\large]{r1}
\Vertex[x=-4,y=0,size=2,label=Scenario \ref{sce1},fontsize=\large]{s1}
\Vertex[x=0,y=0,size=2,label=Scenario \ref{sce2},fontsize=\large]{s2}
\Vertex[x=4,y=0,size=2,label=Scenario \ref{sce3},fontsize=\large]{s3}
\Vertex[x=0,y=-3,shape=rectangle,style={minimum width=6.5cm,minimum height=.75cm},label={\mcp{h}{\beta^1,\phi^1} Fixer-superior to \mcp{h}{\beta^g,\phi^g}},fontsize=\large]{end}
\Edge[Direct,fontsize=\large](r1)(s1)
\Edge[Direct,loopposition=90,fontsize=\large](s1)(s1)
\Edge[Direct,loopposition=90,fontsize=\large](s2)(s2)
\Edge[Direct,loopposition=90,fontsize=\large](s3)(s3)
\Edge[Direct,position={below=.5mm},fontsize=\large](s1)(s2)
\Edge[Direct,position={below=.5mm},fontsize=\large](s2)(s3)
\Edge[Direct,bend=30,position={below=.5mm},fontsize=\large](s1)(s3)
\Edge[Direct,bend=-15,style={dashed}](s1)(end)
\Edge[Direct,bend=15,style={dotted}](s1)(end)
\Edge[Direct,bend=-15,style={dashed}](s2)(end)
\Edge[Direct,bend=15,style={dotted}](s2)(end)
\Edge[Direct,bend=-15,style={dashed}](s3)(end)
\Edge[Direct,bend=15,style={dotted}](s3)(end)
\end{tikzpicture}
\caption{Proof of Proposition \ref{proofExecution}. The solid edges represent Buster playing and Fixer responding without either winning in the given round, the dashed edges represent Fixer winning in the given round, and the dotted edges represent Buster winning in the given round. The edge directed from Round 1 to Scenario \ref{sce1} is detailed in Proposition \ref{scenario1prop0}, the edges directed away from Scenario \ref{sce1} are detailed in Propositions \ref{scenario1prop1}, \ref{scenario1prop2}, and \ref{scenario1prop3}, the edges directed away from Scenario \ref{sce2} are detailed in Propositions \ref{scenario2prop1} and \ref{scenario2prop2}, and the edges directed away from Scenario \ref{sce3} are detailed in Propositions \ref{scenario3prop1} and \ref{scenario3prop2}.}\label{strategyFigure}
\end{figure}

\begin{cor}\label{conclusionCorollary2}
$\phi$ Fixer-dominates $\phi'$ at \mc{h}.
\end{cor}
\begin{proof}
Let $\phi^g$ be an almost greedy Fixer strategy. By Proposition \ref{proofExecution} for any granular Buster strategy $\beta^1$ there exist a Fixer strategy $\phi^1$ and Buster strategy $\beta^g$ for which $\phi^1(\mc{h})=\phi(\mc{h})=\{s\}$ and \mcp{h}{\beta^1,\phi^1} is Fixer-superior to \mcp{h}{\beta^g,\phi^g}. Thus $\phi$ Fixer-dominates $\phi'$ at \mc{h} by Proposition \ref{proofStrategy}.
\end{proof}

\subsection{The case that $\es{G}{\mc{h}}{1}-\es{B}{\mc{h}}{1}$ has at least three components}
\label{c3Proof}

If $\phi^*$ and $\phi^{\#}$ are Fixer strategies and $F$ is a set satisfying $F\subseteq\phi^*(\mc{h})\cap\phi^{\#}(\mc{h})$, then say $\phi^*$ and $\phi^{\#}$ \emph{agree on $F$ at \mc{h}}. If $\phi^{\#}$ is a Fixer strategy and $B\subset\es{B}{\mc{h}}{1}$ and $F\subseteq\phi^{\#}(\mc{h})$ are sets of edges such that $(\es{G}{\mc{h}}{1}-B)\cup F$ is connected, define the \emph{$B,F$-translations of $\phi^{\#}$ and \mc{h}} to be the Fixer strategy $\widetilde{\phi^{\#}}$ and history $\mc{\widetilde{h}}\in\his{B}{1}$ with limit $\widetilde{\ell}=\ell-|B|$ (where $\ell$ is the limit of \mc{h}) for which $\es{G}{\mc{\widetilde{h}}}{1}=(\es{G}{\mc{h}}{1}-B)\cup F$ (which is connected),  $\es{R}{\mc{\widetilde{h}}}{1}=\es{R}{\mc{h}}{1}-F$ (which does not intersect \es{G}{\mc{\widetilde{h}}}{1} because $\es{G}{\mc{h}}{1}\cap\es{R}{\mc{h}}{1}=\emptyset$), $\es{B}{\mc{\widetilde{h}}}{1}=\es{B}{\mc{h}}{1}-B$ (which is nonempty), $\widetilde{\phi^{\#}}(\mc{\widetilde{h}})=\phi^{\#}(\mc{h})-F$ (so that $\es{G}{\mcp{\widetilde{h}}{\widetilde{\phi^{\#}}}}{2}=\es{G}{\mcp{h}{\phi^{\#}}}{2}$ and $\es{R}{\mcp{\widetilde{h}}{\widetilde{\phi^{\#}}}}{2}=\es{R}{\mcp{h}{\phi^{\#}}}{2}$), and $\widetilde{\phi^{\#}}(\mc{\widetilde{h'}})=\phi^{\#}(\mc{h'})$ for any histories \mc{h'} and \mc{\widetilde{h'}} in \his{B}{} such that the first round of \mc{h'} matches that of \mcp{h}{\phi^{\#}}, the first round of \mc{\widetilde{h'}} matches that of \mcp{\widetilde{h}}{\widetilde{\phi^{\#}}}, and subsequent rounds of \mc{h'} and \mc{\widetilde{h'}} are identical to each other. Observe that if $\phi^*$ and $\phi^{\#}$ agree on $F$ at \mc{h} and the $B,F$-translations $\widetilde{\phi^{\#}}$ and \mc{\widetilde{h}} of $\phi^{\#}$ and \mc{h} are such that $\widetilde{\phi^{\#}}$ is Fixer-dominant at \mc{\widetilde{h}}, then by dint of $\phi^{\#}$ having Fixer copy the first part $F$ of her first round move in \mc{h} prescribed by $\phi^*$, before finishing with a Fixer-dominant strategy copied from $\widetilde{\phi^{\#}}$, we should see $\phi^{\#}$ Fixer-dominate $\phi^*$ at \mc{h}; hence we can state Lemma \ref{equivalenceLem} below, leaving the proof details to Appendix \ref{c3propAppendix}.

\begin{lem}\label{equivalenceLem}
Suppose $\phi^*$ and $\phi^{\#}$ are Fixer strategies that agree on $F$ at \mc{h}, and the $B,F$-translations $\widetilde{\phi^{\#}}$ and \mc{\widetilde{h}} of $\phi^{\#}$ and \mc{h} are such that $\widetilde{\phi^{\#}}$ is Fixer-dominant at \mc{\widetilde{h}}. Then $\phi^{\#}$ Fixer-dominates $\phi^*$ at \mc{h}.
\end{lem}

We construct a Fixer strategy $\phi''$ and use Lemma \ref{equivalenceLem} to show $\phi$ is Fixer-superior to $\phi''$ at \mc{h} and $\phi''$ is Fixer-superior to $\phi'$ at \mc{h}, which by the transitivity of Fixer-dominance would imply $\phi$ is Fixer-superior to $\phi'$ at \mc{h}. To do this we produce edges $b\in\es{B}{\mc{h}}{1}$, $c\in\es{B}{\mc{h}}{1}$, and $e\in\phi(\mc{h})$, as well as edge subset $F\subset\phi'(\mc{h})$, such that $\phi''$ agrees with $\phi$ on $\{e\}$ at \mc{h} as well as with $\phi'$ on $F$ at \mc{h}, and the $\{b\},\{e\}$-translations $\widetilde{\phi}$ and \mc{\widetilde{h}} of $\phi$ and \mc{h} as well as the $\es{B}{\mc{h}}{1}-\{c\},F$-translations $\widetilde{\phi''}$ and \mc{\widetilde{h'}} of $\phi''$ and \mc{h} are such that $\widetilde{\phi}$ is Fixer-dominant at \mc{\widetilde{h}} and $\widetilde{\phi''}$ is Fixer-dominant at \mc{\widetilde{h'}}.

Let $e$ be the cheapest edge of $\phi(\mc{h})$, which must exist because $\es{G}{\mc{h}}{1}-\es{B}{\mc{h}}{1}$ has multiple components. Let $e'$ be an edge in $\phi'(\mc{h})$ in a path through \es{G}{\mcp{h}{\phi'}}{2} between the endpoints of $e$, which must exist since $e$ is a bridge in \es{G}{\mcp{h}{\phi}}{2} due to the greediness of $\phi$, and $\es{G}{\mc{h}}{1}-\es{B}{\mc{h}}{1}\subset\es{G}{\mcp{h}{\phi}}{2}-\{e\}$ (both of which must be disconnected) while $\es{G}{\mcp{h}{\phi'}}{2}=(\es{G}{\mc{h}}{1}-\es{B}{\mc{h}}{1})\cup\phi'(\mc{h})$ is connected; see Figure \ref{fig1intersectionSupProp1}. By the following proposition, we can set $F=\phi'(\mc{h})-\{ e'\}$ and $\phi''(\mc{h})=F\cup\{ e\}$ , so $\phi''$ agrees both with $\phi$ on $\{e\}$ at \mc{h} as well as with $\phi'$ on $F$ at \mc{h}. Define $\phi''$ to be greedy at all proper extensions of \mc{h}.

\begin{figure}[htb]
\centering
\subcaptionbox{The dashed lines are the edges of $\phi(\mc{h})$, while the dotted lines are the edges of $\phi'(\mc{h})$.\label{fig1intersectionSupProp1}}[9cm]
{
\begin{tikzpicture}
\Vertex[x=0,y=0,size=1]{a}
\Vertex[x=0,y=2,size=1]{b}
\Vertex[x=1.7,y=-1,size=1]{c}
\Vertex[x=-1.7,y=-1,size=1]{d}
\Edge[label=$e$,position={left},fontsize=\large,style={loosely dashed}](a)(b)
\Edge[style={loosely dashed}](b)(c)
\Edge[style={loosely dashed}](c)(d)
\Edge[style={densely dotted}](b)(d)
\Edge[label=$e'$,position={above},fontsize=\large,style={densely dotted}](d)(a)
\Edge[style={densely dotted}](a)(c)
\end{tikzpicture}
}
\subcaptionbox{The graph \es{G}{\mc{\widetilde{h'}}}{1}, where the edges of $\phi'(\mc{h})-\{ e'\}$ are dotted.\label{fig2intersectionSupProp1}}[9cm]
{
\begin{tikzpicture}
\Vertex[x=0,y=0,size=1]{a}
\Vertex[x=0,y=2,size=1]{b}
\Vertex[x=1.7,y=-1,size=1]{c}
\Vertex[x=-1.7,y=-1,size=1]{d}
\Edge[style={densely dotted}](b)(d)
\Edge[label=$c$,position={right},fontsize=\large](b)(c)
\Edge[style={densely dotted}](a)(c)
\end{tikzpicture}
}
\caption{Two graphs from the proof that $\phi$ is Fixer-dominant at \mc{h} when $\es{G}{h}{1}-\es{B}{h}{1}$ has at least three components. The blobs are the components of $\es{G}{\mc{h}}{1}-\es{B}{\mc{h}}{1}$.}\label{figsForIntersectionSupProp1}
\end{figure}

\begin{prop}\label{c3prop2}
$\phi''(\mc{h})=F\cup\{ e\}$ is a legal move for Fixer.
\end{prop}
\begin{proof}
$\phi''(\mc{h})=F\cup\{ e\}\subseteq\phi'(\mc{h})\cup\phi(\mc{h})\subseteq\es{R}{\mc{h}}{1}$, and \es{G}{\mcp{h}{\phi''}}{2} is connected because $\es{G}{\mcp{h}{\phi''}}{2}=(\es{G}{\mcp{h}{\phi'}}{2}-\{e'\})\cup\{e\}$, where \es{G}{\mcp{h}{\phi'}}{2} is connected and $e$ lies in a path through $(\es{G}{\mcp{h}{\phi'}}{2}-\{e'\})\cup\{e\}$ between the endpoints of $e'$ (such a path can be obtained from the path containing $e'$ through \es{G}{\mcp{h}{\phi'}}{2} between the endpoints of $e$ by replacing $e'$ with $e$).
\end{proof}

Since \es{G}{\mc{h}}{1} is connected but the endpoints of $e$ lie in separate components of $\es{G}{\mc{h}}{1}-\es{B}{\mc{h}}{1}$ (as $\phi(\mc{h})$ is greedy and contains $e$), there exists some edge $b\in\es{B}{\mc{h}}{1}$ on a path through \es{G}{\mc{h}}{1} between the endpoints of $e$. Note that $(\es{G}{\mc{h}}{1}-\{b\})\cup\{e\}$ is connected, since if $\es{G}{\mc{h}}{1}-\{b\}$ is disconnected, then $b$ lies on every path through \es{G}{\mc{h}}{1} between the endpoints of $e$, so $e$ must join the two components of $\es{G}{\mc{h}}{1}-\{b\}$. Furthermore, $\{b\}\subset\es{B}{\mc{h}}{1}$ because $b\in\es{B}{\mc{h}}{1}$ and $|\es{B}{\mc{h}}{1}|\geq 2$, as $\es{G}{\mc{h}}{1}-\es{B}{\mc{h}}{1}$ has at least three components, and $\{e\}\subseteq\phi(\mc{h})$ because $e\in\phi(\mc{h})$. Hence we can let $\widetilde{\phi}$ and \mc{\widetilde{h}} be the $\{b\},\{e\}$-translations of $\phi$ and \mc{h}, so $\mc{\widetilde{h}}\in\his{B}{1}$ has limit $\widetilde{\ell}=\ell-1$ and satisfies $\es{G}{\mc{\widetilde{h}}}{1}=(\es{G}{\mc{h}}{1}-\{b\})\cup\{e\}$, $\es{R}{\mc{\widetilde{h}}}{1}=\es{R}{\mc{h}}{1}-\{ e\}$, and $\es{B}{\mc{\widetilde{h}}}{1}=\es{B}{\mc{h}}{1}-\{b\}$, while $\widetilde{\phi}(\mc{\widetilde{h}})=\phi(\mc{h})-\{e\}$. By the following proposition, $\widetilde{\phi}(\mc{\widetilde{h}})=\phi(\mc{h})-\{e\}$ is a greedy play for Fixer to reconnect $\es{G}{\mc{\widetilde{h}}}{1}-\es{B}{\mc{\widetilde{h}}}{1}$.

\begin{prop}\label{c3prop4}
$\widetilde{\phi}(\mc{\widetilde{h}})=\phi(\mc{h})-\{e\}$ is a greedy play for Fixer.
\end{prop}
\begin{proof}
$\es{G}{\mc{\widetilde{h}}}{1}-\es{B}{\mc{\widetilde{h}}}{1} = ((\es{G}{\mc{h}}{1}-\{b\})\cup\{e\})-(\es{B}{\mc{h}}{1}-\{b\}) = (\es{G}{\mc{h}}{1}-\es{B}{\mc{h}}{1})\cup\{e\}$ and $\phi(\mc{h})=\widetilde{\phi}(\mc{\widetilde{h}})\cup\{e\}$ a greedy play for Fixer to reconnect $\es{G}{\mc{h}}{1}-\es{B}{\mc{h}}{1}$.
\end{proof}

\begin{prop}\label{intersectionSupProp2}
$\phi$ Fixer-dominates $\phi''$ at \mc{h}.
\end{prop}
\begin{proof}
By Proposition \ref{c3prop4} $\widetilde{\phi}(\mc{\widetilde{h}})=\phi(\mc{h})-\{e\}$ is a greedy Fixer move, and thus $\widetilde{\phi}$ prescribes greedy moves for Fixer at every extension of \mc{\widetilde{h}} because $\widetilde{\phi}$ and \mc{\widetilde{h}} are the $\{b\},\{e\}$-translations of $\phi$ and \mc{h}, with $\phi$ prescribing greedy moves for Fixer at every extension of \mc{h}. Hence $\widetilde{\phi}$ is Fixer-dominant at \mc{\widetilde{h}} by Inductive Hypothesis \ref{ih}, as $\es{G}{\mc{\widetilde{h}}}{1}\cup\es{R}{\mc{\widetilde{h}}}{1}=(\es{G}{\mc{h}}{1}\cup\es{R}{\mc{h}}{1})-\{b\}$ where $b\in\es{B}{\mc{h}}{1}\subseteq\es{G}{\mc{h}}{1}$. Since $\phi$ and $\phi''$ agree on $\{e\}$ at \mc{h}, and the $\{b\},\{e\}$-translations $\widetilde{\phi}$ and \mc{\widetilde{h}} of $\phi$ and \mc{h} are such that $\widetilde{\phi}$ is Fixer-dominant at \mc{\widetilde{h}}, by Lemma \ref{equivalenceLem} $\phi$ Fixer-dominates $\phi''$ at \mc{h}.
\end{proof}

We complete the proof by showing that for some edge $c\in\es{B}{\mc{h}}{1}$ the $\es{B}{\mc{h}}{1}-\{c\},F$-translations $\widetilde{\phi''}$ and \mc{\widetilde{h'}} of $\phi''$ and \mc{h} are such that $\widetilde{\phi''}$ is Fixer-dominant at \mc{\widetilde{h'}}. Let $c\in\es{B}{\mc{h}}{1}$ join the two components of $\es{G}{\mcp{h}{\phi'}}{2}-\{e'\}$; see Figure \ref{fig2intersectionSupProp1}. Such an edge $c$ exists, as $e'$ being a bridge in \es{G}{\mcp{h}{\phi'}}{2} (as are all edges of $\phi'(\mc{h})$ by hypothesis) implies $\es{G}{\mcp{h}{\phi'}}{2}-\{e'\}$ has two components, which must have been joined by some edge in \es{B}{\mc{h}}{1} for \es{G}{\mc{h}}{1} to have been connected. Note that $(\es{G}{\mc{h}}{1}-(\es{B}{\mc{h}}{1}-\{c\}))\cup F$ is connected, since $(\es{G}{\mc{h}}{1}-(\es{B}{\mc{h}}{1}-\{c\}))\cup F = (\es{G}{\mc{h}}{1}-\es{B}{\mc{h}}{1})\cup F\cup\{c\} = (\es{G}{\mcp{h}{\phi'}}{2}-\{e'\})\cup\{c\}$, with \es{G}{\mcp{h}{\phi'}}{2} connected and $e'$ and $c$ connecting the same two components of $\es{G}{\mcp{h}{\phi'}}{2}-\{e'\}$. Furthermore $\es{B}{\mc{h}}{1}-\{c\}\subset\es{B}{\mc{h}}{1}$ because $c\in\es{B}{\mc{h}}{1}$, and $F\subseteq\phi''(\mc{h})$ because $\phi''(\mc{h})=F\cup\{ e\}$. Hence there exist $\es{B}{\mc{h}}{1}-\{c\},F$-translations $\widetilde{\phi''}$ and \mc{\widetilde{h'}} of $\phi''$ and \mc{h}, with $\mc{\widetilde{h'}}\in\his{B}{1}$ having limit $\widetilde{\ell'}=\ell-|\es{B}{\mc{h}}{1}|+1$ and satisfying $\es{G}{\mc{\widetilde{h'}}}{1}=(\es{G}{\mc{h}}{1}-(\es{B}{\mc{h}}{1}-\{c\}))\cup F$, $\es{R}{\mc{\widetilde{h'}}}{1}=\es{R}{\mc{h}}{1}-F$, and $\es{B}{\mc{\widetilde{h'}}}{1}=\{c\}$, while $\widetilde{\phi''}(\mc{\widetilde{h'}})=\phi''(\mc{h})-F=\{e\}$. By Proposition \ref{c3prop9} $\widetilde{\phi''}(\mc{\widetilde{h'}})=\{e\}$ is a greedy play for Fixer to reconnect $\es{G}{\mc{\widetilde{h'}}}{1}-\es{B}{\mc{\widetilde{h'}}}{1}$, a fact we use to show $\phi''$ Fixer-dominates $\phi'$ at \mc{h} in Proposition \ref{intersectionSupProp1}, from which we conclude $\phi$ Fixer-dominates $\phi'$ at \mc{h} in Corollary \ref{conclusionCor}.

\begin{prop}\label{c3prop9}
$\widetilde{\phi''}(\mc{\widetilde{h'}})=\{e\}$ is a greedy move for Fixer.
\end{prop}
\begin{proof}
We show for every $r\in\es{R}{\mc{\widetilde{h'}}}{1}$ such that $(\es{G}{\mc{\widetilde{h'}}}{1}-\es{B}{\mc{\widetilde{h'}}}{1})\cup\{r\}$ is connected, we have $w(r)\geq w(e)$, in which case $\widetilde{\phi''}(\mc{\widetilde{h'}})=\{ e\}$ is greedy. First note $r\in\es{R}{\mc{\widetilde{h'}}}{1}\subset\es{R}{\mc{h}}{1}$, and $r$ joins two components of $\es{G}{\mc{h}}{1}-\es{B}{\mc{h}}{1}$ (indeed, $r$ joins the two components of $\es{G}{\mc{\widetilde{h'}}}{1}-\es{B}{\mc{\widetilde{h'}}}{1}=(\es{G}{\mc{h}}{1}-(\es{B}{\mc{h}}{1}-\{c\}))\cup F-\{c\}=(\es{G}{\mc{h}}{1}-\es{B}{\mc{h}}{1})\cup F=\es{G}{\mcp{h}{\phi'}}{2}-\{e'\}$, and $e'\in\phi'(\mc{h})$ so $e'$ is a bridge in \es{G}{\mcp{h}{\phi'}}{2}, so $r$ must join two components of $\es{G}{\mcp{h}{\phi'}}{2}-\phi'(\mc{h})=\es{G}{\mc{h}}{1}-\es{B}{\mc{h}}{1}$). Viewing $\phi(\mc{h})$ as a minimum spanning tree of the multigraph $M$ whose vertices are the components of $\es{G}{\mc{h}}{1}-\es{B}{\mc{h}}{1}$ and whose edges are the edges of \es{R}{\mc{h}}{1} (identifying each endpoint of the edges in \es{R}{\mc{h}}{1} with the component of $\es{G}{\mc{h}}{1}-\es{B}{\mc{h}}{1}$ within which it lies), the spanning tree of $M$ constructed from $\phi(\mc{h})$ by adding $r$ and removing an edge in $\phi(\mc{h})$ in the path through $\phi(\mc{h})$ between the endpoints of $r$ cannot weigh less than $\phi(\mc{h})$ (due to the minimality of $w(\phi(\mc{h}))$), so $r$ cannot weigh less than the cheapest edge $e$ of $\phi(\mc{h})$. Hence $\widetilde{\phi''}(\mc{\widetilde{h'}})=\{ e\}$ is greedy.
\end{proof}

\begin{prop}\label{intersectionSupProp1}
$\phi''$ Fixer-dominates $\phi'$ at \mc{h}.
\end{prop}
\begin{proof}
By Proposition \ref{c3prop9} $\widetilde{\phi''}(\mc{\widetilde{h'}})=\{e\}$ is a greedy Fixer move, and thus $\widetilde{\phi''}$ prescribes greedy moves for Fixer at every extension of \mc{\widetilde{h'}} because $\widetilde{\phi''}$ and \mc{\widetilde{h'}} are the $\es{B}{\mc{h}}{1}-\{c\},F$-translations of $\phi''$ and \mc{h}, with $\phi''$ prescribing greedy moves for Fixer at every proper extension of \mc{h}, so $\widetilde{\phi''}$ is Fixer-dominant at \mc{\widetilde{h'}} by Inductive Hypothesis \ref{ih} because
$$
|\es{G}{\mc{\widetilde{h'}}}{1}|+|\es{R}{\mc{\widetilde{h'}}}{1}|
=
|(\es{G}{\mc{h}}{1}-(\es{B}{\mc{h}}{1}-\{c\}))\cup F|+|\es{R}{\mc{h}}{1}-F|
=
|\es{G}{\mc{h}}{1}|-|\es{B}{\mc{h}}{1}-\{c\}|+|F|+|\es{R}{\mc{h}}{1}|-|F|
=
|\es{G}{\mc{h}}{1}|-|\es{B}{\mc{h}}{1}-\{c\}|+|\es{R}{\mc{h}}{1}|
<
|\es{G}{\mc{h}}{1}|+|\es{R}{\mc{h}}{1}|
$$
as $|\es{B}{\mc{h}}{1}|\geq 2$ since $\es{G}{\mc{h}}{1}-\es{B}{\mc{h}}{1}$ has at least three components. Since $\phi'$ and $\phi''$ agree on $F$ at \mc{h}, and the $\es{B}{h}{1}-\{c\},F$-translations $\widetilde{\phi''}$ and \mc{\widetilde{h'}} of $\phi''$ and \mc{h} are such that $\widetilde{\phi''}$ is Fixer-dominant at \mc{\widetilde{h'}}, by Lemma \ref{equivalenceLem} $\phi''$ Fixer-dominates $\phi'$ at \mc{h}.
\end{proof}

\begin{cor}\label{conclusionCor}
$\phi$ Fixer-dominates $\phi'$ at \mc{h}.
\end{cor}
\begin{proof}
By Proposition \ref{intersectionSupProp2} $\phi$ Fixer-dominates $\phi''$ at \mc{h}, and by Proposition \ref{intersectionSupProp1} $\phi''$ Fixer-dominates $\phi'$ at \mc{h}, so by Proposition \ref{fixDomTransitivityCor} $\phi$ Fixer-dominates $\phi'$ at \mc{h}.
\end{proof}

\vskip 20pt
\begin{ack}
The author thanks R.T. Solo for the problem statement, and the referees for their helpful comments.
\end{ack}

\pagebreak

\appendix

\section{Details of Proof of Proposition \ref{singStrat}}
\label{singStratAppendix}

We provide details of the proof of Proposition \ref{singStrat}, which states that if \ess{F}{p} is a Fixer move at some history $\mc{h}\in\his{B}{p}$ having limit $\ell$, and $\phi^{\#}$ is a Fixer strategy such that for any granular Buster strategy $\beta^1$ (i.e. $\beta^1$ only has Buster remove singletons) there exist a Fixer strategy $\phi^1$ and Buster strategy $\beta^{\#}$ for which $\phi^1(\mc{h})=\ess{F}{p}$ and \mcp{h}{\beta^1,\phi^1} is Fixer-superior to \mcp{h}{\beta^{\#},\phi^{\#}}, then there exists a Fixer strategy $\phi^*$ satisfying $\phi^*(\mc{h})=\ess{F}{p}$ and Fixer-dominating $\phi^{\#}$ at \mc{h}.

\begin{proof}
We let $\phi^*$ be a Fixer strategy satisfying $\phi^*(\mc{h})=\ess{F}{p}$ and whose subsequent moves we are to define against an arbitrary Buster strategy $\beta^*$. We shall derive a particular granular Buster strategy $\beta^1$ from $\beta^*$, so by hypothesis of this proposition there will exist a Fixer strategy $\phi^1$ and Buster strategy $\beta^{\#}$ for which $\phi^1(\mc{h})=\ess{F}{p}$ and \mcp{h}{\beta^1,\phi^1} is Fixer-superior to \mcp{h}{\beta^{\#},\phi^{\#}}. To complete the proof, we construct $\phi^*$ from $\phi^1$ in such a way that \mcp{h}{\beta^*,\phi^*} is Fixer-superior to \mcp{h}{\beta^1,\phi^1}, as this would imply \mcp{h}{\beta^*,\phi^*} is Fixer-superior to \mcp{h}{\beta^{\#},\phi^{\#}} by Proposition \ref{fixSupTransitivityProp}; by the arbitrarity of $\beta^*$, $\phi^*$ would Fixer-dominate $\phi^{\#}$ at \mc{h}.

Noting that $\es{G}{\mcp{h}{\phi^*}}{p+1}=\es{G}{\mcp{h}{\phi^1}}{p+1}$ and $\es{R}{{\mcp{h}{\phi^*}}}{p+1}=\es{R}{{\mcp{h}{\phi^1}}}{p+1}$ because $\mc{h}\in\his{B}{p}$ and $\phi^*(\mc{h})=\ess{F}{p}=\phi^1(\mc{h})$, let $\rho$ be a function satisfying $\rho(p)=p$ and assume inductively for some $k>p$ that at the ends of round $k-1$ of \mcp{h}{\beta^*,\phi^*} and round $\rho(k-1)$ of \mcp{h}{\beta^1,\phi^1} we have $\es{G}{\mcp{h}{\beta^*,\phi^*}}{k}=\es{G}{\mcp{h}{\beta^1,\phi^1}}{\rho(k-1)+1}$ and $\es{R}{\mcp{h}{\beta^*,\phi^*}}{k}=\es{R}{\mcp{h}{\beta^1,\phi^1}}{\rho(k-1)+1}$, whereby $\sum _{j=1}^{k-1}|\es{B}{\mcp{h}{\beta^*,\phi^*}}{j}| = |\es{G}{\mc{h}}{1}|+|\es{R}{\mc{h}}{1}|-|\es{G}{\mcp{h}{\beta^*,\phi^*}}{k}|-|\es{R}{\mcp{h}{\beta^*,\phi^*}}{k}| = |\es{G}{\mc{h}}{1}|+|\es{R}{\mc{h}}{1}|-|\es{G}{\mcp{h}{\beta^1,\phi^1}}{\rho(k-1)+1}|-|\es{R}{\mcp{h}{\beta^1,\phi^1}}{\rho(k-1)+1}| = \sum_{j=1}^{\rho(k-1)}|\es{B}{\mcp{h}{\beta^1,\phi^1}}{j}|$. To complete the proof, we show that the definitions of $\beta^1$, $\phi^*$, and $\rho$ can be extended so that if Fixer or Buster wins \mcp{h}{\beta^*,\phi^*} in the $k$th round, then the same player also wins \mcp{h}{\beta^1,\phi^1} in the $\rho(k)$th round with \mcp{h}{\beta^*,\phi^*} Fixer-superior to \mcp{h}{\beta^1,\phi^1}, and if neither Fixer nor Buster wins \mcp{h}{\beta^*,\phi^*} in the $k$th round, then at the ends of round $k$ of \mcp{h}{\beta^*,\phi^*} and round $\rho(k)$ of \mcp{h}{\beta^1,\phi^1} we have $\es{G}{\mcp{h}{\beta^*,\phi^*}}{k+1}=\es{G}{\mcp{h}{\beta^1,\phi^1}}{\rho(k)+1}$ and $\es{R}{\mcp{h}{\beta^*,\phi^*}}{k+1}=\es{R}{\mcp{h}{\beta^1,\phi^1}}{\rho(k)+1}$ (meaning that once Fixer or Buster finally wins \mcp{h}{\beta^*,\phi^*} in the $j$th round, we'll have $\es{G}{\mcp{h}{\beta^*,\phi^*}}{j}=\es{G}{\mcp{h}{\beta^1,\phi^1}}{\rho(j-1)+1}$ and $\es{R}{\mcp{h}{\beta^*,\phi^*}}{j}=\es{R}{\mcp{h}{\beta^1,\phi^1}}{\rho(j-1)+1}$, and the same analysis of Fixer or Buster winning will be applicable, implying \mcp{h}{\beta^*,\phi^*} is Fixer-superior to \mcp{h}{\beta^1,\phi^1}). 

If Fixer wins \mcp{h}{\beta^*,\phi^*} in the $k$th round, then $\ell=\sum _{j=1}^{k-1}|\es{B}{\mcp{h}{\beta^*,\phi^*}}{j}| = \sum_{j=1}^{\rho(k-1)}|\es{B}{\mcp{h}{\beta^1,\phi^1}}{j}|$, so we can set $\rho(k)=\rho(k-1)+1$ so that Fixer wins \mcp{h}{\beta^1,\phi^1} in the $\rho(k)$th round, leaving $\es{R}{\mcp{h}{\beta^1,\phi^1}}{\rho(k)}=\es{R}{\mcp{h}{\beta^*,\phi^*}}{k}$ and thus $\pay{\mcp{h}{\beta^*,\phi^*}}=1+w(\es{R}{\mcp{h}{\beta^*,\phi^*}}{k})=1+w(\es{R}{\mcp{h}{\beta^1,\phi^1}}{\rho(k)})=\pay{\mcp{h}{\beta^1,\phi^1}}$, so \mcp{h}{\beta^*,\phi^*} is Fixer-superior to \mcp{h}{\beta^1,\phi^1}. 

If Buster wins \mcp{h}{\beta^*,\phi^*} in the $k$th round, then there exists a minimal subset $\{b_1,\ldots,b_m\}$ of \es{B}{\mcp{h}{\beta^*,\phi^*}}{k} such that $(\es{G}{\mcp{h}{\beta^*,\phi^*}}{k}-\{b_1,\ldots,b_m\})\cup\es{R}{\mcp{h}{\beta^*,\phi^*}}{k}$ is disconnected and $\ell \geq m+\sum _{j=1}^{k-1}|\es{B}{\mcp{h}{\beta^*,\phi^*}}{j}| = m+\sum_{j=1}^{\rho(k-1)}|\es{B}{\mcp{h}{\beta^1,\phi^1}}{j}|$. Set $\rho(k)=\rho(k-1)+m$ and define $\beta^1$ to satisfy $\es{B}{\mcp{h}{\beta^1,\phi^1}}{\rho(k-1)+i}=\{b_i\}$ for $1\leq i\leq m$, noting that $\beta^1$ can have Buster play $\{b_1\},\ldots,\{b_m\}$ in rounds $\rho(k-1)+1,\ldots,\rho(k-1)+m$ no matter what Fixer does because $\bigcup_{i=1}^{m}\{b_i\}\subseteq\es{B}{\mcp{h}{\beta^*,\phi^*}}{k}\subseteq\es{G}{\mcp{h}{\beta^*,\phi^*}}{k}=\es{G}{\mcp{h}{\beta^1,\phi^1}}{\rho(k-1)+1}$, with \mcp{h}{\beta^*,\phi^*} and \mcp{h}{\beta^1,\phi^1} both inheriting the limit $\ell$ from \mc{h}, so that Buster wins \mcp{h}{\beta^1,\phi^1} in the $\rho(k)$th round (as $(\es{G}{\mcp{h}{\beta^1,\phi^1}}{\rho(k)}-\es{B}{\mcp{h}{\beta^1,\phi^1}}{\rho(k)})\cup\es{R}{\mcp{h}{\beta^1,\phi^1}}{\rho(k)} = (\es{G}{\mcp{h}{\beta^1,\phi^1}}{\rho(k-1)+1}-\{b_1,\ldots,b_m\})\cup\es{R}{\mcp{h}{\beta^1,\phi^1}}{\rho(k-1)+1} = (\es{G}{\mcp{h}{\beta^*,\phi^*}}{k}-\{b_1,\ldots,b_m\})\cup\es{R}{\mcp{h}{\beta^*,\phi^*}}{k}$, which is disconnected) while leaving $\es{R}{\mcp{h}{\beta^1,\phi^1}}{\rho(k)}\subseteq\es{R}{\mcp{h}{\beta^*,\phi^*}}{k}$ and thus $\pay{\mcp{h}{\beta^*,\phi^*}}=-1-w(\es{R}{\mc{h}}{1}-\es{R}{\mcp{h}{\beta^*,\phi^*}}{k})\geq -1-w(\es{R}{\mc{h}}{1}-\es{R}{\mcp{h}{\beta^1,\phi^1}}{\rho(k)})=\pay{\mcp{h}{\beta^1,\phi^1}}$, so \mcp{h}{\beta^*,\phi^*} is Fixer-superior to \mcp{h}{\beta^1,\phi^1}. 

If neither Fixer nor Buster wins \mcp{h}{\beta^*,\phi^*} in the $k$th round, then $\beta^*$ has Buster play $\es{B}{\mcp{h}{\beta^*,\phi^*}}{k}=\{b_1,\ldots,b_m\}$ such that $(\es{G}{\mcp{h}{\beta^*,\phi^*}}{k}-\{b_1,\ldots,b_m\})\cup\es{R}{\mcp{h}{\beta^*,\phi^*}}{k}$ is connected and $\ell \geq m
+\sum _{j=1}^{k-1}|\es{B}{\mcp{h}{\beta^*,\phi^*}}{j}| = m+\sum_{j=1}^{\rho(k-1)}|\es{B}{\mcp{h}{\beta^1,\phi^1}}{j}|$. Set $\rho(k)=\rho(k-1)+m$ and define $\beta^1$ to satisfy $\es{B}{\mcp{h}{\beta^1,\phi^1}}{\rho(k-1)+i}=\{b_i\}$ for $1\leq i\leq m$ (which, similar to the previous case, are legal plays for Buster no matter what Fixer does) so that neither Fixer nor Buster wins \mcp{h}{\beta^1,\phi^1} in the $\rho(k)$th round or earlier (as $(\es{G}{\mcp{h}{\beta^1,\phi^1}}{\rho(k)}-\es{B}{\mcp{h}{\beta^1,\phi^1}}{\rho(k)})\cup\es{R}{\mcp{h}{\beta^1,\phi^1}}{\rho(k)} = (\es{G}{\mcp{h}{\beta^1,\phi^1}}{\rho(k-1)+1}-\{b_1,\ldots,b_m\})\cup\es{R}{\mcp{h}{\beta^1,\phi^1}}{\rho(k-1)+1} = (\es{G}{\mcp{h}{\beta^*,\phi^*}}{k}-\{b_1,\ldots,b_m\})\cup\es{R}{\mcp{h}{\beta^*,\phi^*}}{k}$, which is connected) . With $\phi^1$ having Fixer respond in rounds $\rho(k-1)+1$ through $\rho(k-1)+m$ in \mcp{h}{\beta^1,\phi^1} with $\es{F}{\mcp{h}{\beta^1,\phi^1}}{\rho(k-1)+1},\ldots,\es{F}{\mcp{h}{\beta^1,\phi^1}}{\rho(k-1)+m}$, define $\phi^*$ so that $\es{F}{\mcp{h}{\beta^*,\phi^*}}{k}=\bigcup_{i=1}^{m}\es{F}{\mcp{h}{\beta^1,\phi^1}}{\rho(k-1)+i}$; note that these are legal plays for Buster and Fixer that leave $\es{G}{\mcp{h}{\beta^*,\phi^*}}{k+1}=\es{G}{\mcp{h}{\beta^1,\phi^1}}{\rho(k)+1}$ and $\es{R}{\mcp{h}{\beta^*,\phi^*}}{k+1}=\es{R}{\mcp{h}{\beta^1,\phi^1}}{\rho(k)+1}$. 

Thus we have \mcp{h}{\beta^*,\phi^*} Fixer-superior to \mcp{h}{\beta^1,\phi^1}, which itself is Fixer-superior to \mcp{h}{\beta^{\#},\phi^{\#}}. By Proposition \ref{fixSupTransitivityProp}, \mcp{h}{\beta^*,\phi^*} is Fixer-superior to \mcp{h}{\beta^{\#},\phi^{\#}}; by the arbitrarity of $\beta^*$, $\phi^*$ Fixer-dominates $\phi^{\#}$ at \mc{h}.
\end{proof}

\pagebreak

\section{Details of Proof of Proposition \ref{subsetSupProp}}
\label{subsetSupPropAppendix}

We provide details of the proof of Proposition \ref{subsetSupProp}, which states that if $\mc{h}\in\his{B}{k}$ and $\phi$ is a Fixer strategy such that for every Fixer strategy $\phi'$ for which every edge of $\phi'(\mc{h})$ is a bridge in \es{G}{\mcp{h}{\phi'}}{k+1}, $\phi$ Fixer-dominates $\phi'$ at \mc{h}, then $\phi$ is Fixer-dominant at \mc{h}.

\begin{proof}
Let $\phi''$ be an arbitrary Fixer strategy. To complete the proof, we construct a Fixer strategy $\phi'$ such that $\phi'(\mc{h})\subseteq\phi''(\mc{h})$, every edge of $\phi'(\mc{h})$ is a bridge in \es{G}{\mcp{h}{\phi'}}{k+1}, and $\phi'$ Fixer-dominates $\phi''$ at \mc{h}. Indeed, by the hypothesis of this proposition $\phi$ would Fixer-dominate $\phi'$ at \mc{h}, and with $\phi'$ Fixer-dominating $\phi''$ at \mc{h} we would also have $\phi$ Fixer-dominating $\phi''$ at \mc{h} by Proposition \ref{fixDomTransitivityCor}; hence $\phi$ would be Fixer-dominant by definition since $\phi''$ was arbitrary. To construct $\phi'$, we define its moves for Fixer against an arbitrary Buster strategy $\beta'$, while simultaneously constructing a Buster strategy $\beta''$ such that \mcp{h}{\beta',\phi'} is Fixer-superior to \mcp{h}{\beta'',\phi''}, from which Fixer-domination of $\phi''$ by $\phi'$ at \mc{h} follows. Note that for $1\leq j\leq k$ we automatically have $\es{G}{\mc{h}}{j}=\es{G}{\mcp{h}{\beta',\phi'}}{j}=\es{G}{\mcp{h}{\beta'',\phi''}}{j}$, $\es{R}{\mc{h}}{j}=\es{R}{\mcp{h}{\beta',\phi'}}{j}=\es{R}{\mcp{h}{\beta'',\phi''}}{j}$, and $\es{B}{\mc{h}}{j}=\es{B}{\mcp{h}{\beta',\phi'}}{j}=\es{B}{\mcp{h}{\beta'',\phi''}}{j}$.

First suppose $(\es{G}{\mc{h}}{k}-\es{B}{\mc{h}}{k})\cup\es{R}{\mc{h}}{k}$ is disconnected, in which case both $(\es{G}{\mcp{h}{\beta',\phi'}}{k}-\es{B}{\mcp{h}{\beta',\phi'}}{k})\cup\es{R}{\mcp{h}{\beta',\phi'}}{k}$ and $(\es{G}{\mcp{h}{\beta'',\phi''}}{k}-\es{B}{\mcp{h}{\beta'',\phi''}}{k})\cup\es{R}{\mcp{h}{\beta'',\phi''}}{k}$ are disconnected, so Buster wins both \mcp{h}{\beta',\phi'} and \mcp{h}{\beta'',\phi''} in the $k$th round, so $\pay{\mcp{h}{\beta',\phi'}}=-1-w(\es{R}{\mcp{h}{\beta',\phi'}}{1}-\es{R}{\mcp{h}{\beta',\phi'}}{k})=-1-w(\es{R}{\mc{h}}{1}-\es{R}{\mc{h}}{k})=-1-w(\es{R}{\mcp{h}{\beta'',\phi''}}{1}-\es{R}{\mcp{h}{\beta'',\phi''}}{k})=\pay{\mcp{h}{\beta'',\phi''}}$ and thus \mcp{h}{\beta',\phi'} is Fixer-superior to \mcp{h}{\beta'',\phi''}.

Now suppose $(\es{G}{\mc{h}}{k}-\es{B}{\mc{h}}{k})\cup\es{R}{\mc{h}}{k}$ is connected, so Buster wins neither \mcp{h}{\beta',\phi'} nor \mcp{h}{\beta'',\phi''} in the $k$th round. Define $\phi'(\mc{h})$ to be a minimal subset of $\phi''(\mc{h})$ such that $\es{G}{\mcp{h}{\phi'}}{k+1}=(\es{G}{\mc{h}}{k}-\es{B}{\mc{h}}{k})\cup\phi'(\mc{h})$ is connected, so every edge of $\phi'(\mc{h})$ is a bridge in \es{G}{\mcp{h}{\phi'}}{k+1}. Note that \es{G}{\mcp{h}{\phi'}}{k+1} and \es{G}{\mcp{h}{\phi''}}{k+1} were constructed from $\es{G}{\mc{h}}{k}-\es{B}{\mc{h}}{k}$ by adding $\phi'(\mc{h})$ and $\phi''(\mc{h})$, respectively, and $\es{R}{\mcp{h}{\phi'}}{k+1}$ and $\es{R}{\mcp{h}{\phi''}}{k+1}$ were constructed from $\es{R}{\mc{h}}{k}$ by removing $\phi'(\mc{h})$ and $\phi''(\mc{h})$, respectively. Hence $\phi'(\mc{h})\subseteq \phi''(\mc{h})$ implies $\es{G}{\mcp{h}{\phi'}}{k+1}\subseteq\es{G}{\mcp{h}{\phi''}}{k+1}$, $\es{R}{\mcp{h}{\phi''}}{k+1}\subseteq \es{R}{\mcp{h}{\phi'}}{k+1}$, and $\es{G}{\mcp{h}{\phi'}}{k+1}\cup\es{R}{\mcp{h}{\phi'}}{k+1}=\es{G}{\mcp{h}{\phi''}}{k+1}\cup\es{R}{\mcp{h}{\phi''}}{k+1}$.

Fixer wins \mcp{h}{\beta',\phi'} in the $(k+1)$st round if and only if Fixer also wins \mcp{h}{\beta'',\phi''} in the $(k+1)$st round because $\sum _{j=1}^{k}|\es{B}{\mcp{h}{\beta',\phi'}}{j}| = \sum _{j=1}^{k}|\es{B}{\mc{h}}{j}| = \sum _{j=1}^{k}|\es{B}{\mcp{h}{\beta'',\phi''}}{j}|$ and thus Buster must have reached the limit in both games if he did so in one, in which case $\pay{\mcp{h}{\beta',\phi'}}=1+w(\es{R}{\mcp{h}{\beta',\phi'}}{k+1})=1+w(\es{R}{\mc{h}}{k})-w(\phi'(\mc{h}))\geq 1+w(\es{R}{\mc{h}}{k})-w(\phi''(\mc{h}))=1+w(\es{R}{\mcp{h}{\beta'',\phi''}}{k+1})=\pay{\mcp{h}{\beta'',\phi''}}$ and thus \mcp{h}{\beta',\phi'} is Fixer-superior to \mcp{h}{\beta'',\phi''}.

If Fixer does not win \mcp{h}{\beta',\phi'} in the $(k+1)$st round, then $\beta'$ has Buster remove some set \es{B}{\mcp{h}{\beta',\phi'}}{k+1} of edges from \es{G}{\mcp{h}{\beta',\phi'}}{k+1} in the $(k+1)$st round of \mcp{h}{\beta',\phi'}. Since $\es{G}{\mcp{h}{\phi'}}{k+1}\subseteq\es{G}{\mcp{h}{\phi''}}{k+1}$, we can define the Buster strategy $\beta''$ to satisfy $\es{B}{\mcp{h}{\beta'',\phi''}}{k+1}=\es{B}{\mcp{h}{\beta',\phi'}}{k+1}$. Note that this implies $(\es{G}{\mcp{h}{\beta',\phi'}}{k+1}-\es{B}{\mcp{h}{\beta',\phi'}}{k+1})\cup\es{R}{\mcp{h}{\beta',\phi'}}{k+1} = (\es{G}{\mcp{h}{\beta',\phi'}}{k+1}\cup\es{R}{\mcp{h}{\beta',\phi'}}{k+1})-\es{B}{\mcp{h}{\beta',\phi'}}{k+1} = (\es{G}{\mcp{h}{\beta'',\phi''}}{k+1}\cup\es{R}{\mcp{h}{\beta'',\phi''}}{k+1})-\es{B}{\mcp{h}{\beta'',\phi''}}{k+1} = (\es{G}{\mcp{h}{\beta'',\phi''}}{k+1}-\es{B}{\mcp{h}{\beta'',\phi''}}{k+1})\cup\es{R}{\mcp{h}{\beta'',\phi''}}{k+1}$.

Buster wins \mcp{h}{\beta',\phi'} in the $(k+1)$st round if and only if Buster also wins \mcp{h}{\beta'',\phi''} in the $(k+1)$st round because $(\es{G}{\mcp{h}{\beta',\phi'}}{k+1}-\es{B}{\mcp{h}{\beta',\phi'}}{k+1})\cup\es{R}{\mcp{h}{\beta',\phi'}}{k+1}$ and $(\es{G}{\mcp{h}{\beta'',\phi''}}{k+1}-\es{B}{\mcp{h}{\beta'',\phi''}}{k+1})\cup\es{R}{\mcp{h}{\beta'',\phi''}}{k+1}$ are the same graph and thus must both be disconnected if one is, in which case $\pay{\mcp{h}{\beta',\phi'}}=-1-w(\es{R}{\mcp{h}{\beta',\phi'}}{1}-\es{R}{\mcp{h}{\beta',\phi'}}{k+1})=-1-w(\es{R}{\mc{h}}{1}-\es{R}{\mc{h}}{k})-w(\phi'(\mc{h}))\geq-1-w(\es{R}{\mc{h}}{1}-\es{R}{\mc{h}}{k})-w(\phi''(\mc{h}))=-1-w(\es{R}{\mcp{h}{\beta'',\phi''}}{1}-\es{R}{\mcp{h}{\beta'',\phi''}}{k+1})=\pay{\mcp{h}{\beta'',\phi''}}$ and thus \mcp{h}{\beta',\phi'} is Fixer-superior to \mcp{h}{\beta'',\phi''}.

Thus we may suppose neither Fixer nor Buster wins either \mcp{h}{\beta',\phi'} or \mcp{h}{\beta'',\phi''} before the $(k+2)$nd round. After letting $\phi''$ have Fixer repair $\es{G}{\mcp{h}{\beta'',\phi''}}{k+1}-\es{B}{\mcp{h}{\beta'',\phi''}}{k+1}$ with the set \es{F}{\mcp{h}{\beta'',\phi''}}{k+1} in the $(k+1)$st round of \mcp{h}{\beta'',\phi''} to create the connected graph \es{G}{\mcp{h}{\beta'',\phi''}}{k+2}, let $\phi'$ have Fixer repair $\es{G}{\mcp{h}{\beta',\phi'}}{k+1}-\es{B}{\mcp{h}{\beta',\phi'}}{k+1}$ with the set $\es{F}{\mcp{h}{\beta',\phi'}}{k+1}=\es{F}{\mcp{h}{\beta'',\phi''}}{k+1}\cup (\es{F}{\mcp{h}{\beta'',\phi''}}{k}-\es{F}{\mcp{h}{\beta',\phi'}}{k})\subseteq\es{R}{\mcp{h}{\beta'',\phi''}}{k+1}\cup (\es{F}{\mcp{h}{\beta'',\phi''}}{k}-\es{F}{\mcp{h}{\beta',\phi'}}{k})=\es{R}{\mcp{h}{\beta',\phi'}}{k+1}$ in the $(k+1)$st round of \mcp{h}{\beta',\phi'} to create the connected graph \es{G}{\mcp{h}{\beta',\phi'}}{k+2}. Note that in comparing \mcp{h}{\beta',\phi'} to \mcp{h}{\beta'',\phi''}, the same total set of edges were removed by Buster over rounds $k$ and $k+1$ (since $\es{B}{\mcp{h}{\beta',\phi'}}{k}=\es{B}{\mc{h}}{k}=\es{B}{\mcp{h}{\beta'',\phi''}}{k}$ and $\es{B}{\mcp{h}{\beta',\phi'}}{k+1}=\es{B}{\mcp{h}{\beta'',\phi''}}{k+1}$), and the same total set of reserve edges were activated by Fixer over rounds $k$ and $k+1$ (since $\es{F}{\mcp{h}{\beta',\phi'}}{k}=\phi'(\mc{h})\subseteq \phi''(\mc{h})=\es{F}{\mcp{h}{\beta'',\phi''}}{k}$ and $\es{F}{\mcp{h}{\beta',\phi'}}{k+1}=\es{F}{\mcp{h}{\beta'',\phi''}}{k+1}\cup (\es{F}{\mcp{h}{\beta'',\phi''}}{k}-\es{F}{\mcp{h}{\beta',\phi'}}{k})$), so \es{F}{\mcp{h}{\beta',\phi'}}{k+1} is a legal move by Fixer that leaves $\es{G}{\mcp{h}{\beta',\phi'}}{k+2}=\es{G}{\mcp{h}{\beta'',\phi''}}{k+2}$ and $\es{R}{\mcp{h}{\beta',\phi'}}{k+2}=\es{R}{\mcp{h}{\beta'',\phi''}}{k+2}$. Furthermore, $\phi'$ and $\beta''$ can continue to be defined for each round $j\geq k+2$ so that the $j$th rounds of \mcp{h}{\beta',\phi'} and \mcp{h}{\beta'',\phi''} are identical by letting $\beta'$ give a Buster play \es{B}{\mcp{h}{\beta',\phi'}}{j} in the $j$th round of \mcp{h}{\beta',\phi'} and letting $\beta''$ have Buster copy that move into the $j$th round of \mcp{h}{\beta'',\phi''} by playing $\es{B}{\mcp{h}{\beta'',\phi''}}{j}=\es{B}{\mcp{h}{\beta',\phi'}}{j}$ (unless Fixer wins \mcp{h}{\beta',\phi'} and thus also \mcp{h}{\beta'',\phi''} in the $j$th round), then letting $\phi''$ give a Fixer response \es{F}{\mcp{h}{\beta'',\phi''}}{j} in \mcp{h}{\beta'',\phi''} and letting $\phi'$ copy that Fixer response into \mcp{h}{\beta',\phi'} by playing $\es{F}{\mcp{h}{\beta',\phi'}}{j}=\es{F}{\mcp{h}{\beta'',\phi''}}{j}$ (unless Buster wins \mcp{h}{\beta'',\phi''} and thus also \mcp{h}{\beta',\phi'} in the $j$th round). Then the same player wins both \mcp{h}{\beta',\phi'} and \mcp{h}{\beta'',\phi''}, with $\pay{\mcp{h}{\beta',\phi'}}=1+w(\es{R}{\mcp{h}{\beta',\phi'}}{|\mcp{h}{\beta',\phi'}|})=1+w(\es{R}{\mcp{h}{\beta'',\phi''}}{|\mcp{h}{\beta'',\phi''}|})=\pay{\mcp{h}{\beta'',\phi''}}$ if Fixer wins both and $\pay{\mcp{h}{\beta',\phi'}}=-1-w(\es{R}{\mcp{h}{\beta',\phi'}}{1}-\es{R}{\mcp{h}{\beta',\phi'}}{|\mcp{h}{\beta',\phi'}|})=-1-w(\es{R}{\mcp{h}{\beta'',\phi''}}{1}-\es{R}{\mcp{h}{\beta'',\phi''}}{|\mcp{h}{\beta'',\phi''}|})=\pay{\mcp{h}{\beta'',\phi''}}$ if Buster wins both, so \mcp{h}{\beta',\phi'} is Fixer-superior to \mcp{h}{\beta'',\phi''}. Thus $\phi'$ Fixer-dominates $\phi''$ at \mc{h}, so $\phi$ is Fixer-dominant at \mc{h}.
\end{proof}

\pagebreak

\section{Details of Proof of Proposition \ref{proofExecution}}
\label{proofExecutionAppendix}

In this appendix we provide the details omitted from the proof of Proposition \ref{proofExecution} from Subsection \ref{c2Proof}, where we had fixed a granular Buster strategy $\beta^1$ and almost greedy Fixer strategy $\phi^g$, then used these to construct a granular Buster strategy $\beta^g$ and Fixer strategy $\phi^1$ such that $\phi^1(\mc{h})=\phi(\mc{h})=\{s\}$ and \mcp{h}{\beta^1,\phi^1} is Fixer-superior to \mcp{h}{\beta^g,\phi^g}. See that subsection for the definitions of granular and almost greedy, as well as for the conditions that must be satisfied at the start of the $k$th round of \mcp{h}{\beta^1,\phi^1} and $j$th round of \mcp{h}{\beta^g,\phi^g} in order for \pair{\mcp{h}{\beta^1,\phi^1}}{k}{\mcp{h}{\beta^g,\phi^g}}{j} to belong to Scenarios \ref{sce1}, \ref{sce2}, or \ref{sce3}, (recalling that \pair{\mcp{h}{\beta^1,\phi^1}}{k}{\mcp{h}{\beta^g,\phi^g}}{j} denotes the pair of histories consisting of the first $k-1$ rounds of \mcp{h}{\beta^1,\phi^1} and the first $j-1$ rounds of \mcp{h}{\beta^g,\phi^g}, awaiting Buster's $k$th round move in the former and Buster's $j$th round move in the latter). We note that the Buster strategy $\beta^g$ happened to be granular, making the following lemma useful for analysis.

\begin{lem}\label{gameLength}
Suppose Buster strategies $\beta^1$ and $\beta^g$ are both granular, and $\phi^1$ and $\phi^g$ are Fixer strategies. If Fixer wins \mcp{h}{\beta^1,\phi^1} in the $k$th round then Fixer can only win \mcp{h}{\beta^g,\phi^g} in exactly the $k$th round and must do so if \mcp{h}{\beta^g,\phi^g} reaches the $k$th round, and if \mcp{h}{\beta^1,\phi^1} reaches the $k$th round without Fixer winning then Fixer cannot win \mcp{h}{\beta^g,\phi^g} in the $k$th round or sooner.
\end{lem}
\begin{proof}
Let $\ell$ be the limit of the game. We have $|\es{B}{\mcp{h}{\beta^1,\phi^1}}{j}|=|\es{B}{\mcp{h}{\beta^g,\phi^g}}{j}|=|\es{B}{\mc{h}}{1}|$ if $j=1$ and $|\es{B}{\mcp{h}{\beta^1,\phi^1}}{j}|=|\es{B}{\mcp{h}{\beta^g,\phi^g}}{j}|=1$ if $2\leq j\leq\ell-|\es{B}{\mc{h}}{1}|+1$. If Fixer wins \mcp{h}{\beta^1,\phi^1} in the $k$th round, then $\ell=\sum _{j=1}^{k-1}|\es{B}{\mcp{h}{\beta^1,\phi^1}}{j}| = |\es{B}{\mc{h}}{1}|+k-2$, and Fixer wins \mcp{h}{\beta^g,\phi^g} if and only if $\ell = \sum_{j=1}^{|\mcp{h}{\beta^g,\phi^g}|-1}|\es{B}{\mcp{h}{\beta^g,\phi^g}}{j}| = |\es{B}{\mc{h}}{1}|+|\mcp{h}{\beta^g,\phi^g}|-2$, so Fixer can only win \mcp{h}{\beta^g,\phi^g} in exactly the $k$th round and must do so if \mcp{h}{\beta^g,\phi^g} reaches the $k$th round. If \mcp{h}{\beta^1,\phi^1} reaches the $k$th round without Fixer winning, then $\ell\geq\sum _{j=1}^{k}|\es{B}{\mcp{h}{\beta^1,\phi^1}}{j}| = |\es{B}{\mc{h}}{1}|+k-1$, and thus $\sum_{j=1}^{k}|\es{B}{\mcp{h}{\beta^g,\phi^g}}{j}| = |\es{B}{\mc{h}}{1}|+k-1\leq\ell$, so Fixer cannot win \mcp{h}{\beta^g,\phi^g} in the $k$th round or sooner.
\end{proof}

\begin{prop}\label{scenario1prop0}
\pair{\mcp{h}{\beta^1,\phi^1}}{2}{\mcp{h}{\beta^g,\phi^g}}{2} belongs to Scenario \ref{sce1}.
\end{prop}
\begin{proof}
Set $X_2$ and $Y_2$ as the two components of the graph $\es{G}{\mc{h}}{1}-\es{B}{\mc{h}}{1}$. Then
\begin{itemize}
\item $s\in\es{F}{\mcp{h}{\beta^1,\phi^1}}{1}\subseteq\es{G}{\mcp{h}{\beta^1,\phi^1}}{2}$, $s'\in\es{F}{\mcp{h}{\beta^g,\phi^g}}{1}\subseteq\es{G}{\mcp{h}{\beta^g,\phi^g}}{2}$, and $\es{G}{\mcp{h}{\beta^1,\phi^1}}{2}-\{s\}=\es{G}{\mcp{h}{\beta^1,\phi^1}}{1}-\es{B}{\mcp{h}{\beta^1,\phi^1}}{1}=\es{G}{\mcp{h}{\beta^g,\phi^g}}{1}-\es{B}{\mcp{h}{\beta^g,\phi^g}}{1}=\es{G}{\mcp{h}{\beta^g,\phi^g}}{2}-\{s'\}$
\item $s'\in\es{R}{\mcp{h}{\beta^1,\phi^1}}{1}-\{s\}=\es{R}{\mcp{h}{\beta^1,\phi^1}}{2}$, $s\in\es{R}{\mcp{h}{\beta^g,\phi^g}}{1}-\{s'\}=\es{R}{\mcp{h}{\beta^g,\phi^g}}{2}$, and $\es{R}{\mcp{h}{\beta^1,\phi^1}}{2}-\{s'\}=\es{R}{\mcp{h}{\beta^1,\phi^1}}{1}-\{s,s'\}=\es{R}{\mcp{h}{\beta^g,\phi^g}}{1}-\{s,s'\}=\es{R}{\mcp{h}{\beta^g,\phi^g}}{2}-\{s\}$
\item $s$ and $s'$ are bridges in \es{G}{\mcp{h}{\beta^1,\phi^1}}{2} and \es{G}{\mcp{h}{\beta^g,\phi^g}}{2}, respectively, between the same subgraphs $X_2$ and $Y_2$
\item $\es{F}{\mcp{h}{\beta^1,\phi^1}}{1}=\{s\}$ being greedy implies $w(r)\geq w(s)$ for every $r\in\es{R}{\mcp{h}{\beta^1,\phi^1}}{2}\cup\es{R}{\mcp{h}{\beta^g,\phi^g}}{2}$ such that $r$ joins $X_2$ to $Y_2$
\end{itemize}
and thus \pair{\mcp{h}{\beta^1,\phi^1}}{2}{\mcp{h}{\beta^g,\phi^g}}{2} belongs to Scenario \ref{sce1}.
\end{proof}

\begin{prop}\label{scenario1prop1}
Suppose \pair{\mcp{h}{\beta^1,\phi^1}}{k}{\mcp{h}{\beta^g,\phi^g}}{k} satisfies the conditions of Scenario \ref{sce1}. If Fixer or Buster wins \mcp{h}{\beta^1,\phi^1} in the $k$th round, then the same player wins \mcp{h}{\beta^g,\phi^g} in the $k$th round, with \mcp{h}{\beta^1,\phi^1} Fixer-superior to \mcp{h}{\beta^g,\phi^g}.
\end{prop}
\begin{proof}
We first consider Fixer winning \mcp{h}{\beta^1,\phi^1} in the $k$th round.

\begin{enumerate}[wide,labelindent=0pt,label=\textbf{Case \arabic*},ref=Case \arabic*]
\item\label{111}\hrulefill

If Fixer wins \mcp{h}{\beta^1,\phi^1} in the $k$th round, then Fixer also wins \mcp{h}{\beta^g,\phi^g} in the $k$th round by Lemma \ref{gameLength}.

\noindent\textbf{Result:} Fixer wins both \mcp{h}{\beta^1,\phi^1} and \mcp{h}{\beta^g,\phi^g} in the $k$th round, with \mcp{h}{\beta^1,\phi^1} Fixer-superior to \mcp{h}{\beta^g,\phi^g} because
$$
\pay{\mcp{h}{\beta^1,\phi^1}}
=
1+w(\es{R}{\mcp{h}{\beta^1,\phi^1}}{k})
=
1+w((\es{R}{\mcp{h}{\beta^g,\phi^g}}{k}-\{s\})\cup\{s'\})
=
1+w(\es{R}{\mcp{h}{\beta^g,\phi^g}}{k})-w(s)+w(s')
\geq
1+w(\es{R}{\mcp{h}{\beta^g,\phi^g}}{k})
=
\pay{\mcp{h}{\beta^g,\phi^g}}
$$

\noindent\hrulefill

We now consider Buster winning \mcp{h}{\beta^1,\phi^1} in the $k$th round.

\item\label{112}\hrulefill

\noindent$\boldsymbol{\es{B}{\mcp{h}{\beta^1,\phi^1}}{k}:}$ $\beta^1$ has Buster play some singleton set \es{B}{\mcp{h}{\beta^1,\phi^1}}{k} to win \mcp{h}{\beta^1,\phi^1} in the $k$th round by leaving $(\es{G}{\mcp{h}{\beta^1,\phi^1}}{k}-\es{B}{\mcp{h}{\beta^1,\phi^1}}{k})\cup\es{R}{\mcp{h}{\beta^1,\phi^1}}{k}$ disconnected. Note that $\es{B}{\mcp{h}{\beta^1,\phi^1}}{k}\neq\{s\}$, since otherwise we would have the contradiction of the connected graph \es{G}{\mcp{h}{\beta^g,\phi^g}}{k} being a spanning subgraph of the disconnected graph $(\es{G}{\mcp{h}{\beta^1,\phi^1}}{k}-\es{B}{\mcp{h}{\beta^1,\phi^1}}{k})\cup\es{R}{\mcp{h}{\beta^1,\phi^1}}{k}$, as $\es{G}{\mcp{h}{\beta^g,\phi^g}}{k}=(\es{G}{\mcp{h}{\beta^1,\phi^1}}{k}-\{s\})\cup\{s'\}\subseteq(\es{G}{\mcp{h}{\beta^1,\phi^1}}{k}-\es{B}{\mcp{h}{\beta^1,\phi^1}}{k})\cup\es{R}{\mcp{h}{\beta^1,\phi^1}}{k}$ if $\es{B}{\mcp{h}{\beta^1,\phi^1}}{k}=\{s\}$.

\noindent$\boldsymbol{\es{B}{\mcp{h}{\beta^g,\phi^g}}{k}:}$ $\beta^g$ has Buster play $\es{B}{\mcp{h}{\beta^g,\phi^g}}{k}=\es{B}{\mcp{h}{\beta^1,\phi^1}}{k}$ in \mcp{h}{\beta^g,\phi^g}, which is possible by Lemma \ref{gameLength} and because $\es{B}{\mcp{h}{\beta^1,\phi^1}}{k}\subseteq\es{G}{\mcp{h}{\beta^1,\phi^1}}{k}-\{s\}\subseteq\es{G}{\mcp{h}{\beta^g,\phi^g}}{k}$. This results in Buster winning \mcp{h}{\beta^g,\phi^g} in the $k$th round because
\begin{align*}
(\es{G}{\mcp{h}{\beta^g,\phi^g}}{k}-\es{B}{\mcp{h}{\beta^g,\phi^g}}{k})\cup\es{R}{\mcp{h}{\beta^g,\phi^g}}{k} &= (\es{G}{\mcp{h}{\beta^g,\phi^g}}{k}\cup\es{R}{\mcp{h}{\beta^g,\phi^g}}{k})-\es{B}{\mcp{h}{\beta^g,\phi^g}}{k} \\
&= ((\es{G}{\mcp{h}{\beta^1,\phi^1}}{k}-\{s\})\cup\{s'\}\cup(\es{R}{\mcp{h}{\beta^1,\phi^1}}{k}-\{s'\})\cup\{s\})-\es{B}{\mcp{h}{\beta^1,\phi^1}}{k} \\
&= (\es{G}{\mcp{h}{\beta^1,\phi^1}}{k}\cup\es{R}{\mcp{h}{\beta^1,\phi^1}}{k})-\es{B}{\mcp{h}{\beta^1,\phi^1}}{k} \\
&= (\es{G}{\mcp{h}{\beta^1,\phi^1}}{k}-\es{B}{\mcp{h}{\beta^1,\phi^1}}{k})\cup\es{R}{\mcp{h}{\beta^1,\phi^1}}{k}
\end{align*}
which is disconnected since Buster wins \mcp{h}{\beta^1,\phi^1} in the $k$th round.

\noindent\textbf{Result:} Buster wins both \mcp{h}{\beta^1,\phi^1} and \mcp{h}{\beta^g,\phi^g} in the $k$th round, with \mcp{h}{\beta^1,\phi^1} Fixer-superior to \mcp{h}{\beta^g,\phi^g} because
\begin{align*}
\pay{\mcp{h}{\beta^1,\phi^1}}
&=
-1-w(\es{R}{\mcp{h}{\beta^1,\phi^1}}{1}-\es{R}{\mcp{h}{\beta^1,\phi^1}}{k})
=
-1-w(\es{R}{\mcp{h}{\beta^g,\phi^g}}{1}-(\es{R}{\mcp{h}{\beta^g,\phi^g}}{k}-\{s\})\cup\{s'\}) \\
&=
-1-w(\es{R}{\mcp{h}{\beta^g,\phi^g}}{1}-\es{R}{\mcp{h}{\beta^g,\phi^g}}{k})-w(s)+w(s')
\geq
-1-w(\es{R}{\mcp{h}{\beta^g,\phi^g}}{1}-\es{R}{\mcp{h}{\beta^g,\phi^g}}{k})
=
\pay{\mcp{h}{\beta^g,\phi^g}}
\end{align*}

\noindent\hrulefill
\end{enumerate}

This completes the proof.
\end{proof}

\begin{prop}\label{scenario1prop2}
Suppose \pair{\mcp{h}{\beta^1,\phi^1}}{k}{\mcp{h}{\beta^g,\phi^g}}{k} satisfies the conditions of Scenario \ref{sce1}. If $\beta^1$ has Buster play $\es{B}{\mcp{h}{\beta^1,\phi^1}}{k}=\{s\}$ in \mcp{h}{\beta^1,\phi^1}, then $\phi^1$ has Fixer play $\es{F}{\mcp{h}{\beta^1,\phi^1}}{k}=\{s'\}$ so that \pair{\mcp{h}{\beta^1,\phi^1}}{k+1}{\mcp{h}{\beta^g,\phi^g}}{k} belongs to Scenario \ref{sce2}.
\end{prop}
\begin{proof}
Letting \mcp{h}{\beta^g,\phi^g} fall a round behind \mcp{h}{\beta^1,\phi^1} and setting $X_{k+1}=X_k$ and $Y_{k+1}=Y_k$, we have
\begin{itemize}
\item $s'\in\es{G}{\mcp{h}{\beta^1,\phi^1}}{k+1}=(\es{G}{\mcp{h}{\beta^1,\phi^1}}{k}-\{s\})\cup\{s'\}=\es{G}{\mcp{h}{\beta^g,\phi^g}}{k}$
\item $\es{R}{\mcp{h}{\beta^1,\phi^1}}{k+1}=\es{R}{\mcp{h}{\beta^1,\phi^1}}{k}-\{s'\}=\es{R}{\mcp{h}{\beta^g,\phi^g}}{k}-\{ s\}$ and $s\in\es{R}{\mcp{h}{\beta^g,\phi^g}}{k}$
\item $s'$ and $s$ are bridges in \es{G}{\mcp{h}{\beta^1,\phi^1}}{k+1} and $(\es{G}{\mcp{h}{\beta^g,\phi^g}}{k}-\{s'\})\cup\{s\}$, respectively, between $X_{k+1}$ and $Y_{k+1}$
\item for every $r\in\es{R}{\mcp{h}{\beta^g,\phi^g}}{k}$ joining $X_{k+1}$ to $Y_{k+1}$, $w(r)\geq w(s)$, since $r\in\es{R}{\mcp{h}{\beta^1,\phi^1}}{k}\cup\es{R}{\mcp{h}{\beta^g,\phi^g}}{k}$ and $r$ would have also joined $X_k$ to $Y_k$
\end{itemize}
and thus \pair{\mcp{h}{\beta^1,\phi^1}}{k+1}{\mcp{h}{\beta^g,\phi^g}}{k} belongs to Scenario \ref{sce2}.
\end{proof}

\begin{prop}\label{scenario1prop3}
Suppose \pair{\mcp{h}{\beta^1,\phi^1}}{k}{\mcp{h}{\beta^g,\phi^g}}{k} satisfies the conditions of Scenario \ref{sce1}. If $\beta^1$ has Buster play $\es{B}{\mcp{h}{\beta^1,\phi^1}}{k}=\{b\}\neq\{s\}$ in \mcp{h}{\beta^1,\phi^1} such that $|\mcp{h}{\beta^1,\phi^1}|>k$, then \pair{\mcp{h}{\beta^1,\phi^1}}{k+1}{\mcp{h}{\beta^g,\phi^g}}{k+1} belongs to either Scenario \ref{sce1} or \ref{sce3}.
\end{prop}
\begin{proof}
Note that $\beta^g$ can always have Buster play a singleton set \es{B}{\mcp{h}{\beta^g,\phi^g}}{k} in \mcp{h}{\beta^g,\phi^g} by Lemma \ref{gameLength}. We first consider the case $\es{G}{\mcp{h}{\beta^1,\phi^1}}{k}-\es{B}{\mcp{h}{\beta^1,\phi^1}}{k}$ is connected.

\begin{enumerate}[wide,labelindent=0pt,label=\textbf{Case \arabic*},ref=Case \arabic*]
\item\label{131}\hrulefill

\noindent$\boldsymbol{\es{B}{\mcp{h}{\beta^1,\phi^1}}{k}:}$ $\beta^1$ has Buster play $\es{B}{\mcp{h}{\beta^1,\phi^1}}{k}=\{b\}$ so that $\es{G}{\mcp{h}{\beta^1,\phi^1}}{k}-\{b\}$ is connected. Note that $b\notin\{s,s'\}$, since $s$ is a bridge in \es{G}{\mcp{h}{\beta^1,\phi^1}}{k} but $\es{G}{\mcp{h}{\beta^1,\phi^1}}{k}-\{b\}$ is connected, and $b\in\es{G}{\mcp{h}{\beta^1,\phi^1}}{k}$ but $s'\notin\es{G}{\mcp{h}{\beta^1,\phi^1}}{k}$.

\noindent$\boldsymbol{\es{B}{\mcp{h}{\beta^g,\phi^g}}{k}:}$ $\beta^g$ has Buster play $\es{B}{\mcp{h}{\beta^g,\phi^g}}{k}=\{b\}\subseteq\es{G}{\mcp{h}{\beta^1,\phi^1}}{k}-\{s\}\subset\es{G}{\mcp{h}{\beta^g,\phi^g}}{k}$.

\noindent$\boldsymbol{\es{F}{\mcp{h}{\beta^g,\phi^g}}{k}:}$ $\phi^g$ will have Fixer respond greedily with $\es{F}{\mcp{h}{\beta^g,\phi^g}}{k}=\emptyset$. Indeed, both $X_k-\{b\}$ and $Y_k-\{b\}$ are connected (since $X_k-\{b\}$ and $Y_k-\{b\}$ are vertex-disjoint subgraphs of the connected graph $\es{G}{\mcp{h}{\beta^1,\phi^1}}{k}-\{b\}$, whose only other edge is the bridge $s$ joining those two subgraphs), and thus $\es{G}{\mcp{h}{\beta^g,\phi^g}}{k}-\es{B}{\mcp{h}{\beta^g,\phi^g}}{k}$ is connected, as it consists solely of the connected subgraphs $X_k-\{b\}$ and $Y_k-\{b\}$ plus the bridge $s'$ between them. Therefore the only greedy response for Fixer in \mcp{h}{\beta^g,\phi^g} is $\es{F}{\mcp{h}{\beta^g,\phi^g}}{k}=\emptyset$.

\noindent$\boldsymbol{\es{F}{\mcp{h}{\beta^1,\phi^1}}{k}:}$ $\phi^1$ has Fixer play $\es{F}{\mcp{h}{\beta^1,\phi^1}}{k}=\emptyset$, which is possible because $\es{G}{\mcp{h}{\beta^1,\phi^1}}{k}-\{b\}$ is connected.

\noindent\textbf{Result:} \pair{\mcp{h}{\beta^1,\phi^1}}{k+1}{\mcp{h}{\beta^g,\phi^g}}{k+1} belongs to Scenario \ref{sce1} because setting $X_{k+1}=X_k-\{b\}$ and $Y_{k+1}=Y_k-\{b\}$ yields
\begin{itemize}
\item $s\in\es{G}{\mcp{h}{\beta^1,\phi^1}}{k}-\{b\}=\es{G}{\mcp{h}{\beta^1,\phi^1}}{k+1}$, $s'\in\es{G}{\mcp{h}{\beta^g,\phi^g}}{k}-\{b\}=\es{G}{\mcp{h}{\beta^g,\phi^g}}{k+1}$, and $\es{G}{\mcp{h}{\beta^1,\phi^1}}{k+1}-\{s\}=\es{G}{\mcp{h}{\beta^1,\phi^1}}{k}-\{b,s\}=\es{G}{\mcp{h}{\beta^g,\phi^g}}{k}-\{b,s'\}=\es{G}{\mcp{h}{\beta^g,\phi^g}}{k+1}-\{s'\}$
\item $s'\in\es{R}{\mcp{h}{\beta^1,\phi^1}}{k}=\es{R}{\mcp{h}{\beta^1,\phi^1}}{k+1}$, $s\in\es{R}{\mcp{h}{\beta^g,\phi^g}}{k}=\es{R}{\mcp{h}{\beta^g,\phi^g}}{k+1}$, and $\es{R}{\mcp{h}{\beta^1,\phi^1}}{k+1}-\{s'\}=\es{R}{\mcp{h}{\beta^1,\phi^1}}{k}-\{s'\}=\es{R}{\mcp{h}{\beta^g,\phi^g}}{k}-\{s\}=\es{R}{\mcp{h}{\beta^g,\phi^g}}{k+1}-\{s\}$
\item $s$ and $s'$ are bridges in \es{G}{\mcp{h}{\beta^1,\phi^1}}{k+1} and \es{G}{\mcp{h}{\beta^g,\phi^g}}{k+1}, respectively, between the connected subgraphs $X_{k+1}$ and $Y_{k+1}$, since they were bridges in \es{G}{\mcp{h}{\beta^1,\phi^1}}{k} and \es{G}{\mcp{h}{\beta^g,\phi^g}}{k} (of which \es{G}{\mcp{h}{\beta^1,\phi^1}}{k+1} and \es{G}{\mcp{h}{\beta^g,\phi^g}}{k+1} are respective subgraphs) between $X_k$ and $Y_k$ (which respectively have the same vertices as $X_{k+1}$ and $Y_{k+1}$)
\item for every $r\in\es{R}{\mcp{h}{\beta^1,\phi^1}}{k+1}\cup\es{R}{\mcp{h}{\beta^g,\phi^g}}{k+1}$ such that $r$ joins $X_{k+1}$ to $Y_{k+1}$, $w(r)\geq w(s)$ since $r\in\es{R}{\mcp{h}{\beta^1,\phi^1}}{k+1}\cup\es{R}{\mcp{h}{\beta^g,\phi^g}}{k+1}=\es{R}{\mcp{h}{\beta^1,\phi^1}}{k}\cup\es{R}{\mcp{h}{\beta^g,\phi^g}}{k}$ and $r$ would have also joined $X_k$ to $Y_k$ (because $X_k$ and $X_{k+1}$ share the same set of vertices, as do $Y_k$ and $Y_{k+1}$)
\end{itemize}

\noindent\hrulefill

For the remaining cases of this proof, by Propositions \ref{scenario1prop1} and \ref{scenario1prop2} as well as \ref{131} of this proof, we may suppose $|\mcp{h}{\beta^1,\phi^1}|>k$ and $\es{B}{\mcp{h}{\beta^1,\phi^1}}{k}=\{b\}\neq\{ s\}$ such that $\es{G}{\mcp{h}{\beta^1,\phi^1}}{k}-\{b\}$ is disconnected. Note that this means that in \mcp{h}{\beta^g,\phi^g} we can let $\beta^g$ have Buster play $\es{B}{\mcp{h}{\beta^g,\phi^g}}{k}=\{b\}$, to which $\phi^g$ will have Fixer respond greedily with $\es{F}{\mcp{h}{\beta^g,\phi^g}}{k}=\{f\}$ for some $f\in\es{R}{\mcp{h}{\beta^g,\phi^g}}{k}$ to create a connected graph $\es{G}{\mcp{h}{\beta^g,\phi^g}}{k+1}=(\es{G}{\mcp{h}{\beta^g,\phi^g}}{k}-\{b\})\cup\{f\}$. First, see that Buster can play $\es{B}{\mcp{h}{\beta^g,\phi^g}}{k}=\{b\}$ because $b\in\es{G}{\mcp{h}{\beta^1,\phi^1}}{k}-\{s\}\subset\es{G}{\mcp{h}{\beta^g,\phi^g}}{k}$. Next, see that $(\es{G}{\mcp{h}{\beta^1,\phi^1}}{k}-\es{B}{\mcp{h}{\beta^1,\phi^1}}{k})\cup\es{R}{\mcp{h}{\beta^1,\phi^1}}{k}=(\es{G}{\mcp{h}{\beta^1,\phi^1}}{k}\cup\es{R}{\mcp{h}{\beta^1,\phi^1}}{k})-\{b\}=(\es{G}{\mcp{h}{\beta^g,\phi^g}}{k}\cup\es{R}{\mcp{h}{\beta^g,\phi^g}}{k})-\{b\}=(\es{G}{\mcp{h}{\beta^g,\phi^g}}{k}-\es{B}{\mcp{h}{\beta^g,\phi^g}}{k})\cup\es{R}{\mcp{h}{\beta^g,\phi^g}}{k}$; since $|\mcp{h}{\beta^1,\phi^1}|>k$ implies $(\es{G}{\mcp{h}{\beta^1,\phi^1}}{k}-\es{B}{\mcp{h}{\beta^1,\phi^1}}{k})\cup\es{R}{\mcp{h}{\beta^1,\phi^1}}{k}$ is connected, $(\es{G}{\mcp{h}{\beta^g,\phi^g}}{k}-\es{B}{\mcp{h}{\beta^g,\phi^g}}{k})\cup\es{R}{\mcp{h}{\beta^g,\phi^g}}{k}$ is also connected, and thus $|\mcp{h}{\beta^g,\phi^g}|>k$. Hence there exists $f\in\es{R}{\mcp{h}{\beta^g,\phi^g}}{k}$ such that Fixer can greedily play $\es{F}{\mcp{h}{\beta^g,\phi^g}}{k}=\{f\}$ to create a connected graph $\es{G}{\mcp{h}{\beta^g,\phi^g}}{k+1}=(\es{G}{\mcp{h}{\beta^g,\phi^g}}{k}-\{b\})\cup\{f\}$.

Since $s$ bridges $X_k$ and $Y_k$ in \es{G}{\mcp{h}{\beta^1,\phi^1}}{k}, either $b\in X_k$ and $X_k-\{b\}$ is disconnected with two components, or $b\in Y_k$ and $Y_k-\{b\}$ is disconnected with two components. Without loss of generality, assume $b\in X_k$ and $X_k-\{b\}$ is disconnected with two components $X_k^1$ and $X_k^2$, with $s$ bridging $X_k^1$ and $Y_k$; see Figure \ref{figGT1Xk1Xk2Yk}. Figures \ref{figGTgXk1Yk} and \ref{figGTgXk2Yk} show the two possibilities for \es{G}{\mcp{h}{\beta^g,\phi^g}}{k}. We first consider the case that $s'$ bridges $X_k^1$ and $Y_k$ in \es{G}{\mcp{h}{\beta^g,\phi^g}}{k} (as in Figure \ref{figGTgXk1Yk}), and later the case that $s'$ bridges $X_k^2$ and $Y_k$ in \es{G}{\mcp{h}{\beta^g,\phi^g}}{k} (as in Figure \ref{figGTgXk2Yk}).

\begin{figure}[htb]
\centering
\subcaptionbox{\es{G}{\mcp{h}{\beta^1,\phi^1}}{k}\label{figGT1Xk1Xk2Yk}}[6cm]
{
\begin{tikzpicture}
\Vertex[x=0,y=1,size=1,label=$X^1_k$,fontsize=\large]{x1k}
\Vertex[x=0,y=-1,size=1,label=$X^2_k$,fontsize=\large]{x2k}
\Vertex[x=2,y=0,size=1.5,label=$Y_k$,fontsize=\large]{yk}
\Edge[label=$s$,bend=45,position={below=.5mm},fontsize=\large](x1k)(yk)
\Edge[label=$b$,position={left},fontsize=\large,style={dashed}](x1k)(x2k)
\end{tikzpicture}
}
\subcaptionbox{One \es{G}{\mcp{h}{\beta^g,\phi^g}}{k} possibility\label{figGTgXk1Yk}}[6cm]
{
\begin{tikzpicture}
\Vertex[x=0,y=1,size=1,label=$X^1_k$,fontsize=\large]{x1k}
\Vertex[x=0,y=-1,size=1,label=$X^2_k$,fontsize=\large]{x2k}
\Vertex[x=2,y=0,size=1.5,label=$Y_k$,fontsize=\large]{yk}
\Edge[label=$s'$,bend=-45,position={above=.5mm},fontsize=\large](x1k)(yk)
\Edge[label=$b$,position={left},fontsize=\large,style={dashed}](x1k)(x2k)
\end{tikzpicture}
}
\subcaptionbox{The other \es{G}{\mcp{h}{\beta^g,\phi^g}}{k} possibility\label{figGTgXk2Yk}}[6cm]
{
\begin{tikzpicture}
\Vertex[x=0,y=1,size=1,label=$X^1_k$,fontsize=\large]{x1k}
\Vertex[x=0,y=-1,size=1,label=$X^2_k$,fontsize=\large]{x2k}
\Vertex[x=2,y=0,size=1.5,label=$Y_k$,fontsize=\large]{yk}
\Edge[label=$s'$,bend=-45,position={above=.5mm},fontsize=\large](x2k)(yk)
\Edge[label=$b$,position={left},fontsize=\large,style={dashed}](x1k)(x2k)
\end{tikzpicture}
}
\caption{\es{G}{\mcp{h}{\beta^1,\phi^1}}{k} and the two possibilities for \es{G}{\mcp{h}{\beta^g,\phi^g}}{k}.}\label{figsThreePossibilities}
\end{figure}

\item\label{132}\hrulefill

\noindent$\boldsymbol{\es{B}{\mcp{h}{\beta^1,\phi^1}}{k}:}$ $\beta^1$ has Buster play $\es{B}{\mcp{h}{\beta^1,\phi^1}}{k}=\{b\}$, with $s'$ bridging $X_k^1$ and $Y_k$ in \es{G}{\mcp{h}{\beta^g,\phi^g}}{k} (as in Figure \ref{figGTgXk1Yk}).

\noindent$\boldsymbol{\es{B}{\mcp{h}{\beta^g,\phi^g}}{k}:}$ $\beta^g$ has Buster play $\es{B}{\mcp{h}{\beta^g,\phi^g}}{k}=\es{B}{\mcp{h}{\beta^1,\phi^1}}{k}=\{b\}$ (so $\es{G}{\mcp{h}{\beta^g,\phi^g}}{k}-\es{B}{\mcp{h}{\beta^g,\phi^g}}{k}$ is as in Figure \ref{figGTgXk1YkFixed}).

\noindent$\boldsymbol{\es{F}{\mcp{h}{\beta^g,\phi^g}}{k}:}$ $\phi^g$ has Fixer respond greedily with $\es{F}{\mcp{h}{\beta^g,\phi^g}}{k}=\{f\}$. Since we showed above that $\es{G}{\mcp{h}{\beta^g,\phi^g}}{k+1}=(\es{G}{\mcp{h}{\beta^g,\phi^g}}{k}-\{b\})\cup\{f\}$ is a connected graph, $f$ must have one endpoint in $X_k^2$ with the other in either $X_k^1$ or $Y_k$. 

\noindent$\boldsymbol{\es{F}{\mcp{h}{\beta^1,\phi^1}}{k}:}$ $\phi^1$ has Fixer play $\es{F}{\mcp{h}{\beta^1,\phi^1}}{k}=\es{F}{\mcp{h}{\beta^g,\phi^g}}{k}=\{f\}$. Note that $f\in\es{R}{\mcp{h}{\beta^1,\phi^1}}{k}$, since $f\in\es{R}{\mcp{h}{\beta^g,\phi^g}}{k}$ and $\es{R}{\mcp{h}{\beta^g,\phi^g}}{k}-\es{R}{\mcp{h}{\beta^1,\phi^1}}{k}=\{s\}\neq\{f\}$, as $s$ has one endpoint in $X^1_k$ and one in $Y_k$ (both vertex sets of which are disjoint from that of $X_k^2$). Furthermore, $\es{G}{\mcp{h}{\beta^1,\phi^1}}{k}-\{b\}$ consists of two components, one $X_k^2$ and the other $X_k^1$ and $Y_k$ being bridged by $s$ (see Figure \ref{figGT1Xk1Xk2Yk}), so Fixer can also play $\es{F}{\mcp{h}{\beta^1,\phi^1}}{k}=\{f\}$ in \mcp{h}{\beta^1,\phi^1} to create a connected graph $\es{G}{\mcp{h}{\beta^1,\phi^1}}{k+1}=(\es{G}{\mcp{h}{\beta^1,\phi^1}}{k}-\{b\})\cup\{f\}$.

\noindent\textbf{Result:} \pair{\mcp{h}{\beta^1,\phi^1}}{k+1}{\mcp{h}{\beta^g,\phi^g}}{k+1} belongs to Scenario \ref{sce1}.

If $f$ bridges $X_k^1$ and $X_k^2$, set $X_{k+1}=(X_k-\{b\})\cup\{f\}$ and $Y_{k+1}=Y_k$. Then $s$ and $s'$ are bridges in \es{G}{\mcp{h}{\beta^1,\phi^1}}{k+1} and \es{G}{\mcp{h}{\beta^g,\phi^g}}{k+1}, respectively, between $X_{k+1}$ and $Y_{k+1}$, as $s$ and $s'$ each bridged $X_k$ and $Y_k$, and the only new edge $f$ resides entirely inside $X_{k+1}$. For every $r\in\es{R}{\mcp{h}{\beta^1,\phi^1}}{k+1}\cup\es{R}{\mcp{h}{\beta^g,\phi^g}}{k+1}$ such that $r$ joins $X_{k+1}$ to $Y_{k+1}$, $w(r)\geq w(s)$ by hypothesis of this scenario, since $r\in\es{R}{\mcp{h}{\beta^1,\phi^1}}{k+1}\cup\es{R}{\mcp{h}{\beta^g,\phi^g}}{k+1}\subset\es{R}{\mcp{h}{\beta^1,\phi^1}}{k}\cup\es{R}{\mcp{h}{\beta^g,\phi^g}}{k}$ and $r$ would have also joined $X_k$ and $Y_k$ (because $X_k$ and $X_{k+1}$ share the same set of vertices, as do $Y_k$ and $Y_{k+1}$).

If $f$ bridges $X_k^2$ and $Y_k$, set $X_{k+1}=X_k^1$ and $Y_{k+1}=X_k^2\cup\{f\}\cup Y_k$. Then $s$ and $s'$ are bridges in \es{G}{\mcp{h}{\beta^1,\phi^1}}{k+1} and \es{G}{\mcp{h}{\beta^g,\phi^g}}{k+1}, respectively, between $X_{k+1}$ and $Y_{k+1}$, as $s$ and $s'$ each bridged $X_k^1$ and $Y_k$, and the only new edge $f$ resides entirely inside $Y_{k+1}$. For every $r\in\es{R}{\mcp{h}{\beta^1,\phi^1}}{k+1}\cup\es{R}{\mcp{h}{\beta^g,\phi^g}}{k+1}$ such that $r$ joins $X_{k+1}$ to $Y_{k+1}$, $w(r)\geq w(s)$ since $r\in\es{R}{\mcp{h}{\beta^1,\phi^1}}{k+1}\cup\es{R}{\mcp{h}{\beta^g,\phi^g}}{k+1}\subset\es{R}{\mcp{h}{\beta^1,\phi^1}}{k}\cup\es{R}{\mcp{h}{\beta^g,\phi^g}}{k}$ and either $r$ joined $X_k^1$ to $Y_k$ (so $w(r)\geq w(s)$ by hypothesis of this scenario), or $r$ joined $X_k^1$ to $X_k^2$ (so $w(r)\geq w(f)$ because $\es{F}{\mcp{h}{\beta^g,\phi^g}}{k}=\{f\}$ was greedy, and $w(f)\geq w(s)$ by hypothesis of this scenario, since $f\in\es{R}{\mcp{h}{\beta^1,\phi^1}}{k}\cup\es{R}{\mcp{h}{\beta^g,\phi^g}}{k}$ and $f$ joined $X_k$ to $Y_k$).

Thus \pair{\mcp{h}{\beta^1,\phi^1}}{k+1}{\mcp{h}{\beta^g,\phi^g}}{k+1} belongs to Scenario \ref{sce1} because
\begin{itemize}
\item $s\in\es{G}{\mcp{h}{\beta^1,\phi^1}}{k}-\{b\}\subset\es{G}{\mcp{h}{\beta^1,\phi^1}}{k+1}$, $s'\in\es{G}{\mcp{h}{\beta^g,\phi^g}}{k}-\{b\}\subset\es{G}{\mcp{h}{\beta^g,\phi^g}}{k+1}$, and $\es{G}{\mcp{h}{\beta^1,\phi^1}}{k+1}-\{s\} = (\es{G}{\mcp{h}{\beta^1,\phi^1}}{k}-\{b,s\})\cup\{f\} = (\es{G}{\mcp{h}{\beta^g,\phi^g}}{k}-\{b,s'\})\cup\{f\} = \es{G}{\mcp{h}{\beta^g,\phi^g}}{k+1}-\{s'\}$
\item $s'\in\es{R}{\mcp{h}{\beta^1,\phi^1}}{k}-\{f\}=\es{R}{\mcp{h}{\beta^1,\phi^1}}{k+1}$, $s\in\es{R}{\mcp{h}{\beta^g,\phi^g}}{k}-\{f\}=\es{R}{\mcp{h}{\beta^g,\phi^g}}{k+1}$, and $\es{R}{\mcp{h}{\beta^1,\phi^1}}{k+1}-\{s'\} = \es{R}{\mcp{h}{\beta^1,\phi^1}}{k}-\{f,s'\} = \es{R}{\mcp{h}{\beta^g,\phi^g}}{k}-\{f,s\} = \es{R}{\mcp{h}{\beta^g,\phi^g}}{k+1}-\{s\}$
\item $s$ and $s'$ are bridges in \es{G}{\mcp{h}{\beta^1,\phi^1}}{k+1} and \es{G}{\mcp{h}{\beta^g,\phi^g}}{k+1}, respectively, between $X_{k+1}$ and $Y_{k+1}$
\item for every $r\in\es{R}{\mcp{h}{\beta^1,\phi^1}}{k+1}\cup\es{R}{\mcp{h}{\beta^g,\phi^g}}{k+1}$ such that $r$ joins $X_{k+1}$ to $Y_{k+1}$, $w(r)\geq w(s)$
\end{itemize}

\noindent\hrulefill

Now suppose $s'$ bridges $X_k^2$ and $Y_k$ in \es{G}{\mcp{h}{\beta^g,\phi^g}}{k} (as in Figure \ref{figGTgXk2Yk}). Let $f$ be some cheapest edge in $\es{R}{\mcp{h}{\beta^g,\phi^g}}{k}$ such that $(\es{G}{\mcp{h}{\beta^g,\phi^g}}{k}-\{b\})\cup\{f\}$ is connected; then either $f$ has one endpoint in $X_k^1$ and the other in $X_k^2$ (as in Figure \ref{figGTgXk1Xk2Fixed}), or $f$ has one endpoint in $X_k^1$ and the other in $Y_k$ (as in Figure \ref{figGTgXk2YkXk1YkFixed}). We first consider the former case and then end the proof considering the latter.

\begin{figure}[htb]
\centering
\subcaptionbox{$s'$ bridges $X_k^1$ and $Y_k$\label{figGTgXk1YkFixed}}[6cm]
{
\begin{tikzpicture}
\Vertex[x=0,y=1,size=1,label=$X^1_k$,fontsize=\large]{x1k}
\Vertex[x=0,y=-1,size=1,label=$X^2_k$,fontsize=\large]{x2k}
\Vertex[x=2,y=0,size=1.5,label=$Y_k$,fontsize=\large]{yk}
\Edge[label=$s'$,bend=-45,position={above=.5mm},fontsize=\large](x1k)(yk)
\end{tikzpicture}
}
\subcaptionbox{$s'$ bridges $X_k^2$ and $Y_k$, and cheapest connecting edge $f\in\es{R}{\mcp{h}{\beta^g,\phi^g}}{k}$ joins $X_k^1$ to $X_k^2$\label{figGTgXk1Xk2Fixed}}[6cm]
{
\begin{tikzpicture}
\Vertex[x=0,y=1,size=1,label=$X^1_k$,fontsize=\large]{x1k}
\Vertex[x=0,y=-1,size=1,label=$X^2_k$,fontsize=\large]{x2k}
\Vertex[x=2,y=0,size=1.5,label=$Y_k$,fontsize=\large]{yk}
\Edge[label=$s'$,bend=-45,position={above=.5mm},fontsize=\large](x2k)(yk)
\Edge[label=$f$,position={left},fontsize=\large,style={loosely dashed}](x1k)(x2k)
\end{tikzpicture}
}
\subcaptionbox{$s'$ bridges $X_k^2$ and $Y_k$, and cheapest connecting edge $f\in\es{R}{\mcp{h}{\beta^g,\phi^g}}{k}$ joins $X_k^1$ to $Y_k$\label{figGTgXk2YkXk1YkFixed}}[6cm]
{
\begin{tikzpicture}
\Vertex[x=0,y=1,size=1,label=$X^1_k$,fontsize=\large]{x1k}
\Vertex[x=0,y=-1,size=1,label=$X^2_k$,fontsize=\large]{x2k}
\Vertex[x=2,y=0,size=1.5,label=$Y_k$,fontsize=\large]{yk}
\Edge[label=$s'$,bend=-45,position={above=.5mm},fontsize=\large](x2k)(yk)
\Edge[label=$f$,position={above},fontsize=\large,style={loosely dashed}](x1k)(yk)
\end{tikzpicture}
}
\caption{The three remaining possibilities for $\es{G}{\mcp{h}{\beta^g,\phi^g}}{k}-\es{B}{\mcp{h}{\beta^g,\phi^g}}{k}$, where in the latter two $f$ is the cheapest edge in \es{R}{\mcp{h}{\beta^g,\phi^g}}{k} such that $(\es{G}{\mcp{h}{\beta^g,\phi^g}}{k}-\es{B}{\mcp{h}{\beta^g,\phi^g}}{k})\cup\{f\}$ is connected.}\label{figsThreePossibilitiesFixed}
\end{figure}

\item\label{133}\hrulefill

\noindent$\boldsymbol{\es{B}{\mcp{h}{\beta^1,\phi^1}}{k}:}$ $\beta^1$ has Buster play $\es{B}{\mcp{h}{\beta^1,\phi^1}}{k}=\{b\}$, with $f$ having one endpoint in $X_k^1$ and the other in $X_k^2$.

\noindent$\boldsymbol{\es{B}{\mcp{h}{\beta^g,\phi^g}}{k}:}$ $\beta^g$ has Buster play $\es{B}{\mcp{h}{\beta^g,\phi^g}}{k}=\es{B}{\mcp{h}{\beta^1,\phi^1}}{k}=\{b\}$.

\noindent$\boldsymbol{\es{F}{\mcp{h}{\beta^g,\phi^g}}{k}:}$ $\phi^g$ has Fixer respond greedily with $\es{F}{\mcp{h}{\beta^g,\phi^g}}{k}=\{f\}$.

\noindent$\boldsymbol{\es{F}{\mcp{h}{\beta^1,\phi^1}}{k}:}$ $\phi^1$ has Fixer play $\es{F}{\mcp{h}{\beta^1,\phi^1}}{k}=\es{F}{\mcp{h}{\beta^g,\phi^g}}{k}=\{f\}$. Note that $f\in\es{R}{\mcp{h}{\beta^1,\phi^1}}{k}$, since $f\in\es{R}{\mcp{h}{\beta^g,\phi^g}}{k}$ and $\es{R}{\mcp{h}{\beta^g,\phi^g}}{k}-\es{R}{\mcp{h}{\beta^1,\phi^1}}{k}=\{s\}\neq\{f\}$, as $f$ has both endpoints in $X_k$ whereas $s$ has one in $Y_k$. Furthermore, $\es{G}{\mcp{h}{\beta^1,\phi^1}}{k}-\{b\}$ consists of two components, one $X_k^2$ and the other the featuring $X_k^1$ and $Y_k$ being bridged by $s$ (see Figure \ref{figGT1Xk1Xk2Yk}), so Fixer can also play $\es{F}{\mcp{h}{\beta^1,\phi^1}}{k}=\{f\}$ in \mcp{h}{\beta^1,\phi^1} to create a connected graph $\es{G}{\mcp{h}{\beta^1,\phi^1}}{k+1}=(\es{G}{\mcp{h}{\beta^1,\phi^1}}{k}-\{b\})\cup\{f\}$.

\noindent\textbf{Result:} \pair{\mcp{h}{\beta^1,\phi^1}}{k+1}{\mcp{h}{\beta^g,\phi^g}}{k+1} belongs to Scenario \ref{sce1} because setting $X_{k+1}=(X_k-\{b\})\cup\{f\}$ and $Y_{k+1}=Y_k$ yields
\begin{itemize}
\item $s\in\es{G}{\mcp{h}{\beta^1,\phi^1}}{k}-\{b\}\subset\es{G}{\mcp{h}{\beta^1,\phi^1}}{k+1}$, $s'\in\es{G}{\mcp{h}{\beta^g,\phi^g}}{k}-\{b\}\subset\es{G}{\mcp{h}{\beta^g,\phi^g}}{k+1}$, and $\es{G}{\mcp{h}{\beta^1,\phi^1}}{k+1}-\{s\} = (\es{G}{\mcp{h}{\beta^1,\phi^1}}{k}-\{b,s\})\cup\{f\} = (\es{G}{\mcp{h}{\beta^g,\phi^g}}{k}-\{b,s'\})\cup\{f\} = \es{G}{\mcp{h}{\beta^g,\phi^g}}{k+1}-\{s'\}$
\item $s'\in\es{R}{\mcp{h}{\beta^1,\phi^1}}{k}-\{f\}=\es{R}{\mcp{h}{\beta^1,\phi^1}}{k+1}$, $s\in\es{R}{\mcp{h}{\beta^g,\phi^g}}{k}-\{f\}=\es{R}{\mcp{h}{\beta^g,\phi^g}}{k+1}$, and $\es{R}{\mcp{h}{\beta^1,\phi^1}}{k+1}-\{s'\} = \es{R}{\mcp{h}{\beta^1,\phi^1}}{k}-\{f,s'\} = \es{R}{\mcp{h}{\beta^g,\phi^g}}{k}-\{f,s\} = \es{R}{\mcp{h}{\beta^g,\phi^g}}{k+1}-\{s\}$
\item $s$ and $s'$ are bridges in \es{G}{\mcp{h}{\beta^1,\phi^1}}{k+1} and \es{G}{\mcp{h}{\beta^g,\phi^g}}{k+1}, respectively, between $X_{k+1}$ and $Y_{k+1}$
\item for every $r\in\es{R}{\mcp{h}{\beta^1,\phi^1}}{k+1}\cup\es{R}{\mcp{h}{\beta^g,\phi^g}}{k+1}$ such that $r$ joins $X_{k+1}$ to $Y_{k+1}$, $w(r)\geq w(s)$ since $r\in\es{R}{\mcp{h}{\beta^1,\phi^1}}{k+1}\cup\es{R}{\mcp{h}{\beta^g,\phi^g}}{k+1}\subset\es{R}{\mcp{h}{\beta^1,\phi^1}}{k}\cup\es{R}{\mcp{h}{\beta^g,\phi^g}}{k}$ and $r$ would have also joined $X_k$ to $Y_k$ (because $X_k$ and $X_{k+1}$ share the same set of vertices, as do $Y_k$ and $Y_{k+1}$)
\end{itemize}

\noindent\hrulefill

Still supposing $s'$ bridges $X_k^2$ and $Y_k$ in \es{G}{\mcp{h}{\beta^g,\phi^g}}{k} (as in Figure \ref{figGTgXk2Yk}), with $f$ the cheapest edge in $\es{R}{\mcp{h}{\beta^g,\phi^g}}{k}$ such that $(\es{G}{\mcp{h}{\beta^g,\phi^g}}{k}-\{b\})\cup\{f\}$ is connected, we now consider the case that $f$ has one endpoint in $X_k^1$ and the other in $Y_k$ (as in Figure \ref{figGTgXk2YkXk1YkFixed}).

\item\label{134}\hrulefill

\noindent$\boldsymbol{\es{B}{\mcp{h}{\beta^1,\phi^1}}{k}:}$ $\beta^1$ has Buster play $\es{B}{\mcp{h}{\beta^1,\phi^1}}{k}=\{b\}$, with $f$ having one endpoint in $X_k^1$ and the other in $Y_k$.

\noindent$\boldsymbol{\es{B}{\mcp{h}{\beta^g,\phi^g}}{k}:}$ $\beta^g$ has Buster play $\es{B}{\mcp{h}{\beta^g,\phi^g}}{k}=\es{B}{\mcp{h}{\beta^1,\phi^1}}{k}=\{b\}$.

\noindent$\boldsymbol{\es{F}{\mcp{h}{\beta^g,\phi^g}}{k}:}$ $\phi^g$ has Fixer respond greedily with $\es{F}{\mcp{h}{\beta^g,\phi^g}}{k}=\{s\}$. Since the cheapest connecting reserve edge in \mcp{h}{\beta^g,\phi^g} is between $X_k^1$ and $Y_k$, $s$ is such an edge by hypothesis of this scenario, so Fixer will greedily play $\es{F}{\mcp{h}{\beta^g,\phi^g}}{k}=\{s\}$ in \mcp{h}{\beta^g,\phi^g} because $\phi^g$ is almost greedy.

\noindent$\boldsymbol{\es{F}{\mcp{h}{\beta^1,\phi^1}}{k}:}$ $\phi^1$ has Fixer play $\es{F}{\mcp{h}{\beta^1,\phi^1}}{k}=\{ s'\}$, which is possible because $\es{G}{\mcp{h}{\beta^1,\phi^1}}{k}-\{b\}$ consists of two components, one $X_k^2$ and the other $X_k^1$ and $Y_k$ being bridged by $s$ (see Figure \ref{figGT1Xk1Xk2Yk}), and $s'\in\es{R}{\mcp{h}{\beta^1,\phi^1}}{k}$ bridges $X_k^2$ and $Y_k$.

\noindent\textbf{Result:} \pair{\mcp{h}{\beta^1,\phi^1}}{k+1}{\mcp{h}{\beta^g,\phi^g}}{k+1} belongs to Scenario \ref{sce3} because
\begin{itemize}
\item $\es{G}{\mcp{h}{\beta^1,\phi^1}}{k+1}=(\es{G}{\mcp{h}{\beta^1,\phi^1}}{k}-\{b\})\cup\{s'\}=(\es{G}{\mcp{h}{\beta^g,\phi^g}}{k}-\{b\})\cup\{s\}=\es{G}{\mcp{h}{\beta^g,\phi^g}}{k+1}$
\item $\es{R}{\mcp{h}{\beta^1,\phi^1}}{k+1}=\es{R}{\mcp{h}{\beta^1,\phi^1}}{k}-\{s'\}=\es{R}{\mcp{h}{\beta^g,\phi^g}}{k}-\{s\}=\es{R}{\mcp{h}{\beta^g,\phi^g}}{k+1}$
\end{itemize}

\noindent\hrulefill
\end{enumerate}

This completes the proof.
\end{proof}

\begin{prop}\label{scenario2prop1}
Suppose \pair{\mcp{h}{\beta^1,\phi^1}}{k}{\mcp{h}{\beta^g,\phi^g}}{k-1} belongs to Scenario \ref{sce2}. If Fixer or Buster wins \mcp{h}{\beta^1,\phi^1} in the $k$th round, then the same player wins \mcp{h}{\beta^g,\phi^g} in the $k$th round, with \mcp{h}{\beta^1,\phi^1} Fixer-superior to \mcp{h}{\beta^g,\phi^g}.
\end{prop}
\begin{proof}
Note that in Scenario \ref{sce2}, by Lemma \ref{gameLength} Fixer does not win \mcp{h}{\beta^g,\phi^g} in the $(k-1)$st round, so $\beta^g$ can always have Buster play a singleton set \es{B}{\mcp{h}{\beta^g,\phi^g}}{k-1}.

Also note that if $\es{B}{\mcp{h}{\beta^g,\phi^g}}{k-1}=\{s'\}$, then as an almost greedy strategy $\phi^g$ will have Fixer respond in \mcp{h}{\beta^g,\phi^g} with $\es{F}{\mcp{h}{\beta^g,\phi^g}}{k-1}=\{s\}$ to leave $s'\in\es{G}{\mcp{h}{\beta^1,\phi^1}}{k}$, $s\in\es{G}{\mcp{h}{\beta^g,\phi^g}}{k}$, $\es{G}{\mcp{h}{\beta^g,\phi^g}}{k}-\{s\}=\es{G}{\mcp{h}{\beta^1,\phi^1}}{k}-\{s'\}$, and $\es{R}{\mcp{h}{\beta^g,\phi^g}}{k}=\es{R}{\mcp{h}{\beta^1,\phi^1}}{k}$. Indeed, since $\es{G}{\mcp{h}{\beta^g,\phi^g}}{k-1}-\{s'\}$ is the graph consisting of the components $X_k$ and $Y_k$, and $s\in\es{R}{\mcp{h}{\beta^g,\phi^g}}{k-1}$ is a bridge between $X_k$ and $Y_k$, $\es{F}{\mcp{h}{\beta^g,\phi^g}}{k-1}=\{s\}$ is a valid move by Fixer in \mcp{h}{\beta^g,\phi^g} (see Figure \ref{scenario2Fig2}, as $\es{G}{\mcp{h}{\beta^g,\phi^g}}{k}=(\es{G}{\mcp{h}{\beta^g,\phi^g}}{k-1}-\{s'\})\cup\{s\}$). Furthermore, $\es{F}{\mcp{h}{\beta^g,\phi^g}}{k-1}=\{s\}$ is a greedy move because for every $r\in\es{R}{\mcp{h}{\beta^g,\phi^g}}{k-1}$ such that $r$ bridges $X_k$ and $Y_k$, $w(r)\geq w(s)$; hence by our stipulation that $\phi^g$ have Fixer use $s$ whenever it is part of a greedy move, we have $\es{F}{\mcp{h}{\beta^g,\phi^g}}{k-1}=\{s\}$. Note that this leaves $s'\in\es{G}{\mcp{h}{\beta^1,\phi^1}}{k}$ (by hypothesis of this scenario), $s\in\es{F}{\mcp{h}{\beta^g,\phi^g}}{k-1}\subseteq\es{G}{\mcp{h}{\beta^g,\phi^g}}{k}$, $\es{G}{\mcp{h}{\beta^g,\phi^g}}{k}-\{s\}=\es{G}{\mcp{h}{\beta^g,\phi^g}}{k-1}-\{s'\}=\es{G}{\mcp{h}{\beta^1,\phi^1}}{k}-\{s'\}$, and $\es{R}{\mcp{h}{\beta^g,\phi^g}}{k}=\es{R}{\mcp{h}{\beta^g,\phi^g}}{k-1}-\{s\}=\es{R}{\mcp{h}{\beta^1,\phi^1}}{k}$.

\begin{enumerate}[wide,labelindent=0pt,label=\textbf{Case \arabic*},ref=Case \arabic*]
\item\label{211}\hrulefill

Fixer wins \mcp{h}{\beta^1,\phi^1} in the $k$th round.

\noindent$\boldsymbol{\es{B}{\mcp{h}{\beta^g,\phi^g}}{k-1}:}$ $\beta^g$ has Buster play $\es{B}{\mcp{h}{\beta^g,\phi^g}}{k-1}=\{s'\}\in\es{G}{\mcp{h}{\beta^g,\phi^g}}{k-1}$ in \mcp{h}{\beta^1,\phi^1}.

\noindent$\boldsymbol{\es{F}{\mcp{h}{\beta^g,\phi^g}}{k-1}:}$ $\phi^g$ will have Fixer will play $\es{F}{\mcp{h}{\beta^g,\phi^g}}{k-1}=\{s\}$ as a greedy response, as shown at the beginning of the proof, leaving $s'\in\es{G}{\mcp{h}{\beta^1,\phi^1}}{k}$, $s\in\es{G}{\mcp{h}{\beta^g,\phi^g}}{k}$, and $\es{G}{\mcp{h}{\beta^g,\phi^g}}{k}-\{s\}=\es{G}{\mcp{h}{\beta^1,\phi^1}}{k}-\{s'\}$ (i.e. the only difference between graphs is $s'$ in \es{G}{\mcp{h}{\beta^1,\phi^1}}{k} being replaced by $s$ in \es{G}{\mcp{h}{\beta^g,\phi^g}}{k}; see Figure \ref{scenario2Fig1} for \es{G}{\mcp{h}{\beta^1,\phi^1}}{k} and Figure \ref{scenario2Fig2} for $\es{G}{\mcp{h}{\beta^g,\phi^g}}{k}=(\es{G}{\mcp{h}{\beta^g,\phi^g}}{k-1}-\{s'\})\cup\{s\}$) as well as $\es{R}{\mcp{h}{\beta^g,\phi^g}}{k}=\es{R}{\mcp{h}{\beta^1,\phi^1}}{k}$. Fixer wins \mcp{h}{\beta^g,\phi^g} in the $k$th round by Lemma \ref{gameLength}.

\noindent\textbf{Result:} Fixer wins both \mcp{h}{\beta^1,\phi^1} and \mcp{h}{\beta^g,\phi^g} in the $k$th round, with \mcp{h}{\beta^1,\phi^1} Fixer-superior to \mcp{h}{\beta^g,\phi^g} because
$$
\pay{\mcp{h}{\beta^1,\phi^1}}
=
1+w(\es{R}{\mcp{h}{\beta^1,\phi^1}}{k})
=
1+w(\es{R}{\mcp{h}{\beta^g,\phi^g}}{k-1}-\{s\})
=
1+w(\es{R}{\mcp{h}{\beta^g,\phi^g}}{k})
=
\pay{\mcp{h}{\beta^g,\phi^g}}
$$

\noindent\hrulefill

Now we consider the case where $\beta^1$ has Buster win $\mcp{h}{\beta^1,\phi^1}$ in the $k$th round. We first consider the subcase of $(\es{G}{\mcp{h}{\beta^1,\phi^1}}{k}-\es{B}{\mcp{h}{\beta^1,\phi^1}}{k})\cup\es{R}{\mcp{h}{\beta^1,\phi^1}}{k}\cup\{s\}$ being disconnected, then finish the proof with the subcase that it is connected.

\item\label{212}\hrulefill

\noindent$\boldsymbol{\es{B}{\mcp{h}{\beta^1,\phi^1}}{k}:}$ $\beta^1$ has Buster win \mcp{h}{\beta^1,\phi^1} in the $k$th round by playing some singleton $\es{B}{\mcp{h}{\beta^1,\phi^1}}{k}$ such that $(\es{G}{\mcp{h}{\beta^1,\phi^1}}{k}-\es{B}{\mcp{h}{\beta^1,\phi^1}}{k})\cup\es{R}{\mcp{h}{\beta^1,\phi^1}}{k}\cup\{s\}$ is disconnected.

\noindent$\boldsymbol{\es{B}{\mcp{h}{\beta^g,\phi^g}}{k-1}:}$ $\beta^g$ has Buster play $\es{B}{\mcp{h}{\beta^g,\phi^g}}{k-1}=\{s'\}$ in \mcp{h}{\beta^g,\phi^g}.

\noindent$\boldsymbol{\es{F}{\mcp{h}{\beta^g,\phi^g}}{k-1}:}$ $\phi^g$ has Fixer respond greedily with $\es{F}{\mcp{h}{\beta^g,\phi^g}}{k-1}=\{s\}$, as shown at the beginning of the proof, leaving $s'\in\es{G}{\mcp{h}{\beta^1,\phi^1}}{k}$, $s\in\es{G}{\mcp{h}{\beta^g,\phi^g}}{k}$, and $\es{G}{\mcp{h}{\beta^g,\phi^g}}{k}-\{s\}=\es{G}{\mcp{h}{\beta^1,\phi^1}}{k}-\{s'\}$ (i.e. the only difference between graphs is $s'$ in \es{G}{\mcp{h}{\beta^1,\phi^1}}{k} being replaced by $s$ in \es{G}{\mcp{h}{\beta^g,\phi^g}}{k}; see Figure \ref{scenario2Fig1} for \es{G}{\mcp{h}{\beta^1,\phi^1}}{k} and Figure \ref{scenario2Fig2} for $\es{G}{\mcp{h}{\beta^g,\phi^g}}{k}=(\es{G}{\mcp{h}{\beta^g,\phi^g}}{k-1}-\{s'\})\cup\{s\}$) as well as $\es{R}{\mcp{h}{\beta^g,\phi^g}}{k}=\es{R}{\mcp{h}{\beta^1,\phi^1}}{k}$.

\noindent$\boldsymbol{\es{B}{\mcp{h}{\beta^g,\phi^g}}{k}:}$ $\beta^g$ has Buster play $\es{B}{\mcp{h}{\beta^g,\phi^g}}{k}=\es{B}{\mcp{h}{\beta^1,\phi^1}}{k}$ in \mcp{h}{\beta^g,\phi^g} if $\es{B}{\mcp{h}{\beta^1,\phi^1}}{k}\neq\{s'\}$, or play $\es{B}{\mcp{h}{\beta^g,\phi^g}}{k}=\{s\}$ in \mcp{h}{\beta^g,\phi^g} if $\es{B}{\mcp{h}{\beta^1,\phi^1}}{k}=\{s'\}$. We show in either case that Buster wins $\mcp{h}{\beta^g,\phi^g}$ in the $k$th round by showing that $(\es{G}{\mcp{h}{\beta^g,\phi^g}}{k}-\es{B}{\mcp{h}{\beta^g,\phi^g}}{k})\cup\es{R}{\mcp{h}{\beta^g,\phi^g}}{k}$ is a spanning subgraph of $(\es{G}{\mcp{h}{\beta^1,\phi^1}}{k}-\es{B}{\mcp{h}{\beta^1,\phi^1}}{k})\cup\es{R}{\mcp{h}{\beta^1,\phi^1}}{k}\cup\{s\}$ and thus disconnected as well. If $\es{B}{\mcp{h}{\beta^1,\phi^1}}{k}\neq\{s'\}$, then $\es{B}{\mcp{h}{\beta^g,\phi^g}}{k}=\es{B}{\mcp{h}{\beta^1,\phi^1}}{k}$ and 
$$
(\es{G}{\mcp{h}{\beta^g,\phi^g}}{k}-\es{B}{\mcp{h}{\beta^g,\phi^g}}{k})\cup\es{R}{\mcp{h}{\beta^g,\phi^g}}{k}=(((\es{G}{\mcp{h}{\beta^1,\phi^1}}{k}-\{s'\})\cup\{s\})-\es{B}{\mcp{h}{\beta^1,\phi^1}}{k})\cup\es{R}{\mcp{h}{\beta^1,\phi^1}}{k}
\subseteq(\es{G}{\mcp{h}{\beta^1,\phi^1}}{k}-\es{B}{\mcp{h}{\beta^1,\phi^1}}{k})\cup\es{R}{\mcp{h}{\beta^1,\phi^1}}{k}\cup\{s\}
$$
since $s\notin\es{B}{\mcp{h}{\beta^1,\phi^1}}{k}$, as $s\notin\es{G}{\mcp{h}{\beta^1,\phi^1}}{k}$ and $\es{B}{\mcp{h}{\beta^1,\phi^1}}{k}\subseteq\es{G}{\mcp{h}{\beta^1,\phi^1}}{k}$. If $\es{B}{\mcp{h}{\beta^1,\phi^1}}{k}=\{s'\}$, then $\es{B}{\mcp{h}{\beta^g,\phi^g}}{k}=\{s\}$ and
$$
(\es{G}{\mcp{h}{\beta^g,\phi^g}}{k}-\es{B}{\mcp{h}{\beta^g,\phi^g}}{k})\cup\es{R}{\mcp{h}{\beta^g,\phi^g}}{k} = (((\es{G}{\mcp{h}{\beta^1,\phi^1}}{k}-\{s'\})\cup\{s\})-\{s\})\cup\es{R}{\mcp{h}{\beta^1,\phi^1}}{k} 
= (\es{G}{\mcp{h}{\beta^1,\phi^1}}{k}-\es{B}{\mcp{h}{\beta^1,\phi^1}}{k})\cup\es{R}{\mcp{h}{\beta^1,\phi^1}}{k}
$$
since $s\notin\es{G}{\mcp{h}{\beta^1,\phi^1}}{k}$.

\noindent\textbf{Result:} Buster wins both $\mcp{h}{\beta^1,\phi^1}$ and $\mcp{h}{\beta^g,\phi^g}$ in the $k$th round, with \mcp{h}{\beta^1,\phi^1} Fixer-superior to \mcp{h}{\beta^g,\phi^g} because
$$
\pay{\mcp{h}{\beta^1,\phi^1}}
=
-1-w(\es{R}{\mcp{h}{\beta^1,\phi^1}}{1}-\es{R}{\mcp{h}{\beta^1,\phi^1}}{k})
=
-1-w(\es{R}{\mcp{h}{\beta^g,\phi^g}}{1}-(\es{R}{\mcp{h}{\beta^g,\phi^g}}{k-1}-\{s\}))
=
-1-w(\es{R}{\mcp{h}{\beta^g,\phi^g}}{1}-\es{R}{\mcp{h}{\beta^g,\phi^g}}{k})
=
\pay{\mcp{h}{\beta^g,\phi^g}}
$$

\noindent\hrulefill

Now we consider our final subcase of the proof, where $\beta^1$ has Buster win $\mcp{h}{\beta^1,\phi^1}$ in the $k$th round and $(\es{G}{\mcp{h}{\beta^1,\phi^1}}{k}-\es{B}{\mcp{h}{\beta^1,\phi^1}}{k})\cup\es{R}{\mcp{h}{\beta^1,\phi^1}}{k}\cup\{s\}$ is connected.

\item\label{213}\hrulefill

\noindent$\boldsymbol{\es{B}{\mcp{h}{\beta^1,\phi^1}}{k}:}$ $\beta^1$ has Buster win \mcp{h}{\beta^1,\phi^1} in the $k$th round by playing some singleton $\es{B}{\mcp{h}{\beta^1,\phi^1}}{k}$ such that $(\es{G}{\mcp{h}{\beta^1,\phi^1}}{k}-\es{B}{\mcp{h}{\beta^1,\phi^1}}{k})\cup\es{R}{\mcp{h}{\beta^1,\phi^1}}{k}\cup\{s\}$ is connected.

\noindent$\boldsymbol{\es{B}{\mcp{h}{\beta^g,\phi^g}}{k-1}:}$ $\beta^g$ has Buster play $\es{B}{\mcp{h}{\beta^g,\phi^g}}{k-1}=\es{B}{\mcp{h}{\beta^1,\phi^1}}{k}$ in \mcp{h}{\beta^g,\phi^g}, which is possible since $\es{B}{\mcp{h}{\beta^1,\phi^1}}{k}\subseteq\es{G}{\mcp{h}{\beta^1,\phi^1}}{k}=\es{G}{\mcp{h}{\beta^g,\phi^g}}{k-1}$. Note that Buster does not win $\mcp{h}{\beta^g,\phi^g}$ in the $(k-1)$st round since $(\es{G}{\mcp{h}{\beta^g,\phi^g}}{k-1}-\es{B}{\mcp{h}{\beta^g,\phi^g}}{k-1})\cup\es{R}{\mcp{h}{\beta^g,\phi^g}}{k-1}=(\es{G}{\mcp{h}{\beta^1,\phi^1}}{k}-\es{B}{\mcp{h}{\beta^1,\phi^1}}{k})\cup\es{R}{\mcp{h}{\beta^1,\phi^1}}{k}\cup\{s\}$, which is connected by hypothesis.

\noindent$\boldsymbol{\es{F}{\mcp{h}{\beta^g,\phi^g}}{k-1}:}$ $\phi^g$ has Fixer play $\es{F}{\mcp{h}{\beta^g,\phi^g}}{k-1}=\{s\}$, since if $s\notin\es{F}{\mcp{h}{\beta^g,\phi^g}}{k-1}$ then $(\es{G}{\mcp{h}{\beta^g,\phi^g}}{k-1}-\es{B}{\mcp{h}{\beta^g,\phi^g}}{k-1})\cup\es{F}{\mcp{h}{\beta^g,\phi^g}}{k-1}\subseteq(\es{G}{\mcp{h}{\beta^1,\phi^1}}{k}-\es{B}{\mcp{h}{\beta^1,\phi^1}}{k})\cup\es{R}{\mcp{h}{\beta^1,\phi^1}}{k}$, which is disconnected because Buster wins \mcp{h}{\beta^1,\phi^1} in the $k$th round.

\noindent$\boldsymbol{\es{B}{\mcp{h}{\beta^g,\phi^g}}{k}:}$ $\beta^g$ has Buster play $\es{B}{\mcp{h}{\beta^g,\phi^g}}{k}=\{s\}$ in \mcp{h}{\beta^g,\phi^g}, resulting in Buster winning $\mcp{h}{\beta^g,\phi^g}$ in the $k$th round, since
\begin{align*}
(\es{G}{\mcp{h}{\beta^g,\phi^g}}{k}-\es{B}{\mcp{h}{\beta^g,\phi^g}}{k})\cup\es{R}{\mcp{h}{\beta^g,\phi^g}}{k}&=(((\es{G}{\mcp{h}{\beta^g,\phi^g}}{k-1}-\es{B}{\mcp{h}{\beta^g,\phi^g}}{k-1})\cup\es{F}{\mcp{h}{\beta^g,\phi^g}}{k-1})-\{s\})\cup(\es{R}{\mcp{h}{\beta^g,\phi^g}}{k-1}-\es{F}{\mcp{h}{\beta^g,\phi^g}}{k-1}) \\
&=(\es{G}{\mcp{h}{\beta^1,\phi^1}}{k}-\es{B}{\mcp{h}{\beta^1,\phi^1}}{k})\cup(\es{R}{\mcp{h}{\beta^g,\phi^g}}{k-1}-\{s\}) \\
&=(\es{G}{\mcp{h}{\beta^1,\phi^1}}{k}-\es{B}{\mcp{h}{\beta^1,\phi^1}}{k})\cup\es{R}{\mcp{h}{\beta^1,\phi^1}}{k}
\end{align*}
which is disconnected because Buster wins $\mcp{h}{\beta^1,\phi^1}$ in the $k$th round.

\noindent\textbf{Result:} Buster wins both $\mcp{h}{\beta^1,\phi^1}$ and $\mcp{h}{\beta^g,\phi^g}$ in the $k$th round, with \mcp{h}{\beta^1,\phi^1} Fixer-superior to \mcp{h}{\beta^g,\phi^g} because
$$
\pay{\mcp{h}{\beta^1,\phi^1}}
=
-1-w(\es{R}{\mcp{h}{\beta^1,\phi^1}}{1}-\es{R}{\mcp{h}{\beta^1,\phi^1}}{k})
=
-1-w(\es{R}{\mcp{h}{\beta^g,\phi^g}}{1}-(\es{R}{\mcp{h}{\beta^g,\phi^g}}{k-1}-\{s\}))
=
-1-w(\es{R}{\mcp{h}{\beta^g,\phi^g}}{1}-\es{R}{\mcp{h}{\beta^g,\phi^g}}{k})
=
\pay{\mcp{h}{\beta^g,\phi^g}}
$$

\noindent\hrulefill
\end{enumerate}

This completes the proof.
\end{proof}

\begin{prop}\label{scenario2prop2}
Suppose \pair{\mcp{h}{\beta^1,\phi^1}}{k}{\mcp{h}{\beta^g,\phi^g}}{k-1} belongs to Scenario \ref{sce2}. If neither Fixer nor Buster wins \mcp{h}{\beta^1,\phi^1} in the $k$th round, then either \pair{\mcp{h}{\beta^1,\phi^1}}{k+1}{\mcp{h}{\beta^g,\phi^g}}{k} belongs to Scenario \ref{sce2} or \pair{\mcp{h}{\beta^1,\phi^1}}{k+1}{\mcp{h}{\beta^g,\phi^g}}{k+1} belongs to Scenario \ref{sce3}.
\end{prop}
\begin{proof}
By Lemma \ref{gameLength} Fixer does not win \mcp{h}{\beta^g,\phi^g} in either the $(k-1)$st or $k$th round, so after $\beta^1$ has Buster play some singleton set $\es{B}{\mcp{h}{\beta^g,\phi^g}}{k-1}$ in the $k$th round of \mcp{h}{\beta^1,\phi^1}, $\beta^g$ can always have Buster copy that move by playing $\es{B}{\mcp{h}{\beta^g,\phi^g}}{k-1}=\es{B}{\mcp{h}{\beta^1,\phi^1}}{k}\subseteq\es{G}{\mcp{h}{\beta^1,\phi^1}}{k}=\es{G}{\mcp{h}{\beta^g,\phi^g}}{k-1}$ in the $(k-1)$st round of \mcp{h}{\beta^g,\phi^g}. Furthermore, Buster would not win \mcp{h}{\beta^g,\phi^g} in the $(k-1)$st round, as Buster not winning \mcp{h}{\beta^1,\phi^1} in the $k$th round implies $(\es{G}{\mcp{h}{\beta^1,\phi^1}}{k}-\es{B}{\mcp{h}{\beta^1,\phi^1}}{k})\cup\es{R}{\mcp{h}{\beta^1,\phi^1}}{k}$ is connected, in which case so must be $(\es{G}{\mcp{h}{\beta^g,\phi^g}}{k-1}-\es{B}{\mcp{h}{\beta^g,\phi^g}}{k-1})\cup\es{R}{\mcp{h}{\beta^g,\phi^g}}{k-1}$ because $\es{G}{\mcp{h}{\beta^1,\phi^1}}{k}=\es{G}{\mcp{h}{\beta^g,\phi^g}}{k-1}$, $\es{B}{\mcp{h}{\beta^1,\phi^1}}{k}=\es{B}{\mcp{h}{\beta^g,\phi^g}}{k-1}$, and $\es{R}{\mcp{h}{\beta^1,\phi^1}}{k}\subseteq\es{R}{\mcp{h}{\beta^g,\phi^g}}{k-1}$. Hence $\beta^g$ can have Buster play a singleton set \es{B}{\mcp{h}{\beta^g,\phi^g}}{k} in the $k$th round of \mcp{h}{\beta^g,\phi^g}. We separate into two cases whether Fixer can respond greedily in \mcp{h}{\beta^g,\phi^g} with $\es{F}{\mcp{h}{\beta^g,\phi^g}}{k-1}=\{s\}$, starting with the case where Fixer cannot greedily play $\es{F}{\mcp{h}{\beta^g,\phi^g}}{k-1}=\{s\}$ and finishing with the case where Fixer can greedily play $\es{F}{\mcp{h}{\beta^g,\phi^g}}{k-1}=\{s\}$ (which Fixer will do in \mcp{h}{\beta^g,\phi^g} because $\phi^g$ is almost greedy).

\begin{enumerate}[wide,labelindent=0pt,label=\textbf{Case \arabic*},ref=Case \arabic*]
\item\label{221}\hrulefill

\noindent$\boldsymbol{\es{B}{\mcp{h}{\beta^1,\phi^1}}{k}:}$ $\beta^1$ has Buster play $\es{B}{\mcp{h}{\beta^1,\phi^1}}{k}=\{b\}$ for some edge $b\in\es{G}{\mcp{h}{\beta^1,\phi^1}}{k}$ so that Buster does not win $\mcp{h}{\beta^1,\phi^1}$ in the $k$th round.

\noindent$\boldsymbol{\es{B}{\mcp{h}{\beta^g,\phi^g}}{k-1}:}$ $\beta^g$ has Buster play $\es{B}{\mcp{h}{\beta^g,\phi^g}}{k-1}=\es{B}{\mcp{h}{\beta^1,\phi^1}}{k}=\{b\}$.

\noindent$\boldsymbol{\es{F}{\mcp{h}{\beta^g,\phi^g}}{k-1}:}$ $\phi^g$ has Fixer respond greedily in \mcp{h}{\beta^g,\phi^g} with a set $\es{F}{\mcp{h}{\beta^g,\phi^g}}{k-1}\neq\{s\}$, as no possible greedy response for Fixer in \mcp{h}{\beta^g,\phi^g} contains $s$ by assumption. We note that this implies $b\neq s'$, as if $b=s'$ then Fixer could have greedily played $\es{F}{\mcp{h}{\beta^g,\phi^g}}{k-1}=\{s\}$. Indeed, $\es{F}{\mcp{h}{\beta^g,\phi^g}}{k-1}=\{s\}$ would have been a valid Fixer response since $s$ is a bridge in $(\es{G}{\mcp{h}{\beta^g,\phi^g}}{k-1}-\{s'\})\cup\{s\}$ between connected subgraphs $X_k$ and $Y_k$ (see Figure \ref{scenario2Fig2}). Furthermore, it would have been a greedy response since $s'$ being a bridge between $X_k$ and $Y_k$ (see Figure \ref{scenario2Fig1}) implies that any greedy response would have to be a single edge in \es{R}{\mcp{h}{\beta^g,\phi^g}}{k-1} bridging $X_k$ and $Y_k$, and for every $r\in\es{R}{\mcp{h}{\beta^g,\phi^g}}{k-1}$ such that $r$ bridges $X_k$ and $Y_k$, $w(r)\geq w(s)$ by assumption of this scenario.

\noindent$\boldsymbol{\es{F}{\mcp{h}{\beta^1,\phi^1}}{k}:}$ $\phi^1$ has Fixer play $\es{F}{\mcp{h}{\beta^1,\phi^1}}{k}=\es{F}{\mcp{h}{\beta^g,\phi^g}}{k-1}$, which is possible because $\es{F}{\mcp{h}{\beta^g,\phi^g}}{k-1}\subseteq\es{R}{\mcp{h}{\beta^g,\phi^g}}{k-1}-\{s\}=\es{R}{\mcp{h}{\beta^1,\phi^1}}{k}$ and 
\begin{align*}
\es{G}{\mcp{h}{\beta^1,\phi^1}}{k+1} = (\es{G}{\mcp{h}{\beta^1,\phi^1}}{k}-\es{B}{\mcp{h}{\beta^1,\phi^1}}{k})\cup\es{F}{\mcp{h}{\beta^1,\phi^1}}{k} = (\es{G}{\mcp{h}{\beta^g,\phi^g}}{k-1}-\es{B}{\mcp{h}{\beta^g,\phi^g}}{k-1})\cup\es{F}{\mcp{h}{\beta^g,\phi^g}}{k-1} = \es{G}{\mcp{h}{\beta^g,\phi^g}}{k}
\end{align*}
which is connected.

\noindent\textbf{Result:} \pair{\mcp{h}{\beta^1,\phi^1}}{k+1}{\mcp{h}{\beta^g,\phi^g}}{k} belongs to Scenario \ref{sce2}. We first identify the two components $X_{k+1}$ and $Y_{k+1}$ of $\es{G}{\mcp{h}{\beta^1,\phi^1}}{k+1}-\{s'\}$ for which $s'$ and $s$ are bridges in \es{G}{\mcp{h}{\beta^1,\phi^1}}{k+1} and $(\es{G}{\mcp{h}{\beta^g,\phi^g}}{k}-\{s'\})\cup\{s\}$, respectively, between $X_{k+1}$ and $Y_{k+1}$, with every $r\in\es{R}{\mcp{h}{\beta^g,\phi^g}}{k}$ that bridges $X_{k+1}$ and $Y_{k+1}$ also satisfying $w(r)\geq w(s)$.

If $b$ is not a bridge in $X_k$ or $Y_k$, then both $X_k-\{b\}$ and $Y_k-\{b\}$ are connected, and furthermore $\es{G}{\mcp{h}{\beta^g,\phi^g}}{k-1}-\{b\}$ is connected because $s'$ is a bridge between $X_k-\{b\}$ and $Y_k-\{b\}$ while $b\neq s'$. Hence $\es{F}{\mcp{h}{\beta^1,\phi^1}}{k}=\es{F}{\mcp{h}{\beta^g,\phi^g}}{k-1}=\emptyset$ since that would be the only greedy move by Fixer in \mcp{h}{\beta^g,\phi^g}, so we can set $X_{k+1}=X_k-\{b\}$ and $Y_{k+1}=Y_k-\{b\}$, since $s$ and $s'$ were bridges in \es{G}{\mcp{h}{\beta^1,\phi^1}}{k} and $(\es{G}{\mcp{h}{\beta^g,\phi^g}}{k-1}-\{s'\})\cup\{s\}$, respectively, between $X_{k}$ and $Y_{k}$, and for every $r\in\es{R}{\mcp{h}{\beta^g,\phi^g}}{k}$ such that $r$ bridges $X_{k+1}$ and $Y_{k+1}$, $r$ bridged $X_k$ and $Y_k$ in which case $w(r)\geq w(s)$ by the assumptions of this scenario.

Thus without loss of generality we may assume $b$ is a bridge in $X_k$ between the two components $X^1_k$ and $X^2_k$ of $X_k-\{b\}$ (the case that $b$ is a bridge in $Y_k$ is identical), and $\es{F}{\mcp{h}{\beta^1,\phi^1}}{k}=\es{F}{\mcp{h}{\beta^g,\phi^g}}{k-1}=\{f\}$ for some edge $f\in\es{R}{\mcp{h}{\beta^1,\phi^1}}{k}$ either bridging $X^1_k$ and $X^2_k$, or bridging $Y_k$ and one of $X^1_k$ or $X^2_k$. If $f$ bridges $X^1_k$ and $X^2_k$ (see Figure \ref{scenario2Prop2Fig1}), then we may set $X_{k+1}=X^1_k\cup\{f\}\cup X^2_k$ and $Y_{k+1}=Y_k$, since $s$ and $s'$ were bridges in \es{G}{\mcp{h}{\beta^1,\phi^1}}{k} and $(\es{G}{\mcp{h}{\beta^g,\phi^g}}{k-1}-\{s'\})\cup\{s\}$, respectively, between $X_{k}$ and $Y_{k}$, and for every $r\in\es{R}{\mcp{h}{\beta^g,\phi^g}}{k}$ such that $r$ bridges $X_{k+1}$ and $Y_{k+1}$, $r$ bridged $X_k$ and $Y_k$ in which case $w(r)\geq w(s)$ by the assumptions of this scenario. If $f$ does not bridge $X^1_k$ and $X^2_k$, then without loss of generality we may assume $f$ bridges $X^1_k$ and $Y_k$ (see Figure \ref{scenario2Prop2Fig2}). Note that this implies $s$ and $s'$ are bridges between $X^2_{k}$ and $Y_{k}$, since $s$ bridging $X^1_k$ and $Y_k$ would contradict the assumption that no possible greedy response for Fixer in \mcp{h}{\beta^g,\phi^g} contains $s$, and $s'$ bridging $X^1_k$ and $Y_k$ would contradict the number of components of $\es{G}{\mcp{h}{\beta^g,\phi^g}}{k-1}-\es{B}{\mcp{h}{\beta^g,\phi^g}}{k-1}$ being reduced with the addition of $f$. Further note $w(f)\geq w(s)$ by the assumptions of this scenario, and $w(f)\leq w(r)$ for any $r\in\es{R}{\mcp{h}{\beta^g,\phi^g}}{k-1}$ bridging $X^1_k$ and $X^2_k$, since otherwise $\es{F}{\mcp{h}{\beta^g,\phi^g}}{k-1}=\{r\}$ would have been a cheaper valid response for Fixer in \mcp{h}{\beta^g,\phi^g}, contradicting $\es{F}{\mcp{h}{\beta^g,\phi^g}}{k-1}=\{f\}$ being greedy. Hence we can set $X_{k+1}=X^2_k$ and $Y_{k+1}=Y_k\cup\{f\}\cup X^1_k$, since $s$ and $s'$ were bridges in \es{G}{\mcp{h}{\beta^1,\phi^1}}{k} and $(\es{G}{\mcp{h}{\beta^g,\phi^g}}{k-1}-\{s'\})\cup\{s\}$, respectively, between $X^2_{k}$ and $Y_{k}$, and for every $r\in\es{R}{\mcp{h}{\beta^g,\phi^g}}{k}$ such that $r$ bridges $X_{k+1}$ and $Y_{k+1}$, either $r$ bridged $X_k$ and $Y_k$ in which case $w(r)\geq w(s)$ by the assumptions of this scenario, or $r$ bridged $X^1_k$ and $X^2_k$, in which case we've already shown $w(r)\geq w(f)\geq w(s)$.

Thus \pair{\mcp{h}{\beta^1,\phi^1}}{k+1}{\mcp{h}{\beta^g,\phi^g}}{k} belongs to Scenario \ref{sce2} because
\begin{itemize}
\item $s'\in\es{G}{\mcp{h}{\beta^1,\phi^1}}{k+1}=\es{G}{\mcp{h}{\beta^g,\phi^g}}{k}$
\item $\es{R}{\mcp{h}{\beta^1,\phi^1}}{k+1} = \es{R}{\mcp{h}{\beta^1,\phi^1}}{k}-\es{F}{\mcp{h}{\beta^1,\phi^1}}{k} = (\es{R}{\mcp{h}{\beta^g,\phi^g}}{k-1}-\{ s\})-\es{F}{\mcp{h}{\beta^g,\phi^g}}{k-1} = \es{R}{\mcp{h}{\beta^g,\phi^g}}{k}-\{ s\}$ and $s\in\es{R}{\mcp{h}{\beta^g,\phi^g}}{k-1}-\es{F}{\mcp{h}{\beta^g,\phi^g}}{k-1}=\es{R}{\mcp{h}{\beta^g,\phi^g}}{k}$
\item $s'$ and $s$ are bridges in \es{G}{\mcp{h}{\beta^1,\phi^1}}{k} and $(\es{G}{\mcp{h}{\beta^g,\phi^g}}{k-1}-\{s'\})\cup\{s\}$, respectively, between $X_{k+1}$ and $Y_{k+1}$
\item for every $r\in\es{R}{\mcp{h}{\beta^g,\phi^g}}{k}$ such that $r$ bridges $X_{k+1}$ and $Y_{k+1}$, $w(r)\geq w(s)$
\end{itemize}

\begin{figure}[htb]
\centering
\subcaptionbox{$f$ joins $X_k^1$ to $X_k^2$ ($s'$ bridges $Y_k$ to either $X_k^1$ or $X_k^2$, with the latter shown here)\label{scenario2Prop2Fig1}}[9cm]
{
\begin{tikzpicture}
\Vertex[x=0,y=1,size=1,label=$X^1_k$,fontsize=\large]{x1k}
\Vertex[x=0,y=-1,size=1,label=$X^2_k$,fontsize=\large]{x2k}
\Vertex[x=2,y=0,size=1.5,label=$Y_k$,fontsize=\large]{yk}
\Edge[label=$s'$,bend=-45,position={above=.5mm},fontsize=\large](x2k)(yk)
\Edge[label=$f$,position={left},fontsize=\large,style={loosely dashed}](x1k)(x2k)
\end{tikzpicture}
}
\subcaptionbox{$f$ joins $X_k^1$ to $Y_k$ (meaning $s'$ bridges $X_k^2$ and $Y_k$; the case where $f$ joins $X_k^2$ to $Y_k$ and $s'$ bridges $X_k^1$ and $Y_k$ is similar)\label{scenario2Prop2Fig2}}[9cm]
{
\begin{tikzpicture}
\Vertex[x=0,y=1,size=1,label=$X^1_k$,fontsize=\large]{x1k}
\Vertex[x=0,y=-1,size=1,label=$X^2_k$,fontsize=\large]{x2k}
\Vertex[x=2,y=0,size=1.5,label=$Y_k$,fontsize=\large]{yk}
\Edge[label=$s'$,bend=-45,position={above=.5mm},fontsize=\large](x2k)(yk)
\Edge[label=$f$,position={above},fontsize=\large,style={loosely dashed}](x1k)(yk)
\end{tikzpicture}
}
\caption{Possibilities for $\es{G}{\mcp{h}{\beta^1,\phi^1}}{k+1}=(\es{G}{\mcp{h}{\beta^1,\phi^1}}{k}-\es{B}{\mcp{h}{\beta^1,\phi^1}}{k})\cup\es{F}{\mcp{h}{\beta^1,\phi^1}}{k}=(\es{G}{\mcp{h}{\beta^1,\phi^1}}{k}-\{b\})\cup\{f\}$ in the proof of \ref{221} of Proposition \ref{scenario2prop2}, if $b$ is a bridge in $X_k$ between $X_k^1$ and $X_k^2$.}\label{scenario2Prop2Figs}
\end{figure}

\noindent\hrulefill

We finish with the case where Fixer can greedily play $\es{F}{\mcp{h}{\beta^g,\phi^g}}{k-1}=\{s\}$.

\item\label{222}\hrulefill

\noindent$\boldsymbol{\es{B}{\mcp{h}{\beta^1,\phi^1}}{k}:}$ $\beta^1$ has Buster play some singleton $\es{B}{\mcp{h}{\beta^1,\phi^1}}{k}$ so that Buster does not win $\mcp{h}{\beta^1,\phi^1}$ in the $k$th round.

\noindent$\boldsymbol{\es{B}{\mcp{h}{\beta^g,\phi^g}}{k-1}:}$ $\beta^g$ has Buster play $\es{B}{\mcp{h}{\beta^g,\phi^g}}{k-1}=\es{B}{\mcp{h}{\beta^1,\phi^1}}{k}$.

\noindent$\boldsymbol{\es{F}{\mcp{h}{\beta^g,\phi^g}}{k-1}:}$ $\phi^g$ has Fixer respond greedily in \mcp{h}{\beta^g,\phi^g} with $\es{F}{\mcp{h}{\beta^g,\phi^g}}{k-1}=\{s\}$.

\noindent$\boldsymbol{\es{B}{\mcp{h}{\beta^g,\phi^g}}{k}:}$ $\beta^g$ has Buster play $\es{B}{\mcp{h}{\beta^g,\phi^g}}{k}=\{s\}$, which is possible because $\{s\}=\es{F}{\mcp{h}{\beta^g,\phi^g}}{k-1}\subseteq\es{G}{\mcp{h}{\beta^g,\phi^g}}{k}$.

\noindent$\boldsymbol{\es{F}{\mcp{h}{\beta^g,\phi^g}}{k}:}$ $\phi^g$ will have Fixer respond by creating a connected graph \es{G}{\mcp{h}{\beta^g,\phi^g}}{k+1} with some greedy \es{F}{\mcp{h}{\beta^g,\phi^g}}{k} in \mcp{h}{\beta^g,\phi^g}. Fixer can respond with some $\es{F}{\mcp{h}{\beta^g,\phi^g}}{k}\subseteq\es{R}{\mcp{h}{\beta^g,\phi^g}}{k}$ such that $\es{G}{\mcp{h}{\beta^g,\phi^g}}{k+1}=(\es{G}{\mcp{h}{\beta^g,\phi^g}}{k}-\es{B}{\mcp{h}{\beta^g,\phi^g}}{k})\cup\es{F}{\mcp{h}{\beta^g,\phi^g}}{k}$ is connected (so Buster doesn't win \mcp{h}{\beta^g,\phi^g} in the $k$th round), because Buster's failure to win \mcp{h}{\beta^1,\phi^1} in the $k$th round implies $(\es{G}{\mcp{h}{\beta^1,\phi^1}}{k}-\es{B}{\mcp{h}{\beta^1,\phi^1}}{k})\cup\es{R}{\mcp{h}{\beta^1,\phi^1}}{k}$ is connected, and
\begin{align*}
(\es{G}{\mcp{h}{\beta^g,\phi^g}}{k}-\es{B}{\mcp{h}{\beta^g,\phi^g}}{k})\cup\es{R}{\mcp{h}{\beta^g,\phi^g}}{k}&=(\es{G}{\mcp{h}{\beta^g,\phi^g}}{k}\cup\es{R}{\mcp{h}{\beta^g,\phi^g}}{k})-\{s\}=((\es{G}{\mcp{h}{\beta^g,\phi^g}}{k-1}\cup\es{R}{\mcp{h}{\beta^g,\phi^g}}{k-1})-\es{B}{\mcp{h}{\beta^g,\phi^g}}{k-1})-\{s\} \\
&=(\es{G}{\mcp{h}{\beta^g,\phi^g}}{k-1}-\es{B}{\mcp{h}{\beta^g,\phi^g}}{k-1})\cup(\es{R}{\mcp{h}{\beta^g,\phi^g}}{k-1}-\{s\})=(\es{G}{\mcp{h}{\beta^1,\phi^1}}{k}-\es{B}{\mcp{h}{\beta^1,\phi^1}}{k})\cup\es{R}{\mcp{h}{\beta^1,\phi^1}}{k}
\end{align*}
so $(\es{G}{\mcp{h}{\beta^g,\phi^g}}{k}-\es{B}{\mcp{h}{\beta^g,\phi^g}}{k})\cup\es{R}{\mcp{h}{\beta^g,\phi^g}}{k}$ is connected as well.

\noindent$\boldsymbol{\es{F}{\mcp{h}{\beta^1,\phi^1}}{k}:}$ $\phi^1$ has Fixer play $\es{F}{\mcp{h}{\beta^1,\phi^1}}{k} = \es{F}{\mcp{h}{\beta^g,\phi^g}}{k} \subseteq\es{R}{\mcp{h}{\beta^g,\phi^g}}{k}=\es{R}{\mcp{h}{\beta^g,\phi^g}}{k-1}-\{s\} = \es{R}{\mcp{h}{\beta^1,\phi^1}}{k}$, which is possible because
\begin{align*}
\es{G}{\mcp{h}{\beta^1,\phi^1}}{k+1}&=(\es{G}{\mcp{h}{\beta^1,\phi^1}}{k}-\es{B}{\mcp{h}{\beta^1,\phi^1}}{k})\cup\es{F}{\mcp{h}{\beta^1,\phi^1}}{k}=(\es{G}{\mcp{h}{\beta^g,\phi^g}}{k-1}-\es{B}{\mcp{h}{\beta^g,\phi^g}}{k-1})\cup\es{F}{\mcp{h}{\beta^g,\phi^g}}{k} \\
&=((\es{G}{\mcp{h}{\beta^g,\phi^g}}{k-1}-\es{B}{\mcp{h}{\beta^g,\phi^g}}{k-1})\cup\es{F}{\mcp{h}{\beta^g,\phi^g}}{k-1})-\{s\})\cup\es{F}{\mcp{h}{\beta^g,\phi^g}}{k}=(\es{G}{\mcp{h}{\beta^g,\phi^g}}{k}-\es{B}{\mcp{h}{\beta^g,\phi^g}}{k})\cup\es{F}{\mcp{h}{\beta^g,\phi^g}}{k}=\es{G}{\mcp{h}{\beta^g,\phi^g}}{k+1}
\end{align*}
which we already showed was connected because there existed a valid Fixer move \es{F}{\mcp{h}{\beta^g,\phi^g}}{k} to prevent Buster from winning \mcp{h}{\beta^g,\phi^g} in the $k$th round.

\noindent\textbf{Result:} \pair{\mcp{h}{\beta^1,\phi^1}}{k+1}{\mcp{h}{\beta^g,\phi^g}}{k+1} belongs to Scenario \ref{sce3} because
\begin{itemize}
\item $\es{G}{\mcp{h}{\beta^1,\phi^1}}{k+1}=\es{G}{\mcp{h}{\beta^g,\phi^g}}{k+1}$
\item $
\es{R}{\mcp{h}{\beta^1,\phi^1}}{k+1}
=
\es{R}{\mcp{h}{\beta^1,\phi^1}}{k}-\es{F}{\mcp{h}{\beta^1,\phi^1}}{k}
=
(\es{R}{\mcp{h}{\beta^g,\phi^g}}{k-1}-\{s\})-\es{F}{\mcp{h}{\beta^g,\phi^g}}{k}
=
(\es{R}{\mcp{h}{\beta^g,\phi^g}}{k-1}-\es{F}{\mcp{h}{\beta^g,\phi^g}}{k-1})-\es{F}{\mcp{h}{\beta^g,\phi^g}}{k}
=
\es{R}{\mcp{h}{\beta^g,\phi^g}}{k+1}
$
\end{itemize}

\noindent\hrulefill
\end{enumerate}

This completes the proof.
\end{proof}

\begin{prop}\label{scenario3prop1}
Suppose \pair{\mcp{h}{\beta^1,\phi^1}}{k}{\mcp{h}{\beta^g,\phi^g}}{k} satisfies the conditions of Scenario \ref{sce3}. If Fixer or Buster wins \mcp{h}{\beta^1,\phi^1} in the $k$th round, then the same player wins \mcp{h}{\beta^g,\phi^g} in the $k$th round, with \mcp{h}{\beta^1,\phi^1} Fixer-superior to \mcp{h}{\beta^g,\phi^g}.
\end{prop}
\begin{proof}
We first consider Fixer winning \mcp{h}{\beta^1,\phi^1} in the $k$th round.

\begin{enumerate}[wide,labelindent=0pt,label=\textbf{Case \arabic*},ref=Case \arabic*]
\item\label{311}\hrulefill

If Fixer wins \mcp{h}{\beta^1,\phi^1} in the $k$th round, then Fixer also wins \mcp{h}{\beta^g,\phi^g} in the $k$th round by Lemma \ref{gameLength}.

\noindent\textbf{Result:} Fixer wins \mcp{h}{\beta^1,\phi^1} and \mcp{h}{\beta^g,\phi^g} in the $k$th round, with \mcp{h}{\beta^1,\phi^1} Fixer-superior to \mcp{h}{\beta^g,\phi^g} because
$$
\pay{\mcp{h}{\beta^1,\phi^1}}
=
1+w(\es{R}{\mcp{h}{\beta^1,\phi^1}}{k})
=
1+w(\es{R}{\mcp{h}{\beta^g,\phi^g}}{k})
=
\pay{\mcp{h}{\beta^g,\phi^g}}
$$

\noindent\hrulefill

We now consider Buster winning \mcp{h}{\beta^1,\phi^1} in the $k$th round.

\item\label{312}\hrulefill

\noindent$\boldsymbol{\es{B}{\mcp{h}{\beta^1,\phi^1}}{k}:}$ $\beta^1$ has Buster play \es{B}{\mcp{h}{\beta^1,\phi^1}}{k} to win $\mcp{h}{\beta^1,\phi^1}$ in the $k$th round.

\noindent$\boldsymbol{\es{B}{\mcp{h}{\beta^g,\phi^g}}{k}:}$ By Lemma \ref{gameLength}, Fixer does not win \mcp{h}{\beta^g,\phi^g} in the $k$th round, so $\beta^g$ can have Buster play $\es{B}{\mcp{h}{\beta^g,\phi^g}}{k}=\es{B}{\mcp{h}{\beta^1,\phi^1}}{k}\subseteq\es{G}{\mcp{h}{\beta^1,\phi^1}}{k}=\es{G}{\mcp{h}{\beta^g,\phi^g}}{k}$, which wins \mcp{h}{\beta^g,\phi^g} in the $k$th round because $(\es{G}{\mcp{h}{\beta^g,\phi^g}}{k}-\es{B}{\mcp{h}{\beta^g,\phi^g}}{k})\cup\es{R}{\mcp{h}{\beta^g,\phi^g}}{k}=(\es{G}{\mcp{h}{\beta^1,\phi^1}}{k}-\es{B}{\mcp{h}{\beta^1,\phi^1}}{k})\cup\es{R}{\mcp{h}{\beta^1,\phi^1}}{k}$, which is disconnected since Buster wins $\mcp{h}{\beta^1,\phi^1}$ in the $k$th round.

\noindent\textbf{Result:} Buster wins \mcp{h}{\beta^1,\phi^1} and \mcp{h}{\beta^g,\phi^g} in the $k$th round, with \mcp{h}{\beta^1,\phi^1} Fixer-superior to \mcp{h}{\beta^g,\phi^g} because
$$
\pay{\mcp{h}{\beta^1,\phi^1}}
=
-1-w(\es{R}{\mcp{h}{\beta^1,\phi^1}}{1}-\es{R}{\mcp{h}{\beta^1,\phi^1}}{k})
=
-1-w(\es{R}{\mcp{h}{\beta^g,\phi^g}}{1}-\es{R}{\mcp{h}{\beta^g,\phi^g}}{k})
=
\pay{\mcp{h}{\beta^g,\phi^g}}
$$

\noindent\hrulefill
\end{enumerate}

This completes the proof.
\end{proof}

\begin{prop}\label{scenario3prop2}
Suppose \pair{\mcp{h}{\beta^1,\phi^1}}{k}{\mcp{h}{\beta^g,\phi^g}}{k} satisfies the conditions of Scenario \ref{sce3}. If neither Fixer nor Buster wins \mcp{h}{\beta^1,\phi^1} in the $k$th round, then \pair{\mcp{h}{\beta^1,\phi^1}}{k+1}{\mcp{h}{\beta^g,\phi^g}}{k+1} belongs to Scenario \ref{sce3}.
\end{prop}
\begin{proof}
\noindent\hrulefill

\noindent$\boldsymbol{\es{B}{\mcp{h}{\beta^1,\phi^1}}{k}:}$ $\beta^1$ has Buster play some singleton set \es{B}{\mcp{h}{\beta^1,\phi^1}}{k} that does not win \mcp{h}{\beta^1,\phi^1} for Buster in the $k$th round.

\noindent$\boldsymbol{\es{B}{\mcp{h}{\beta^g,\phi^g}}{k}:}$ $\beta^g$ has Buster play $\es{B}{\mcp{h}{\beta^g,\phi^g}}{k}=\es{B}{\mcp{h}{\beta^1,\phi^1}}{k}$, as by Lemma \ref{gameLength} Fixer does not win \mcp{h}{\beta^g,\phi^g} in the $k$th round. Buster does not win \mcp{h}{\beta^g,\phi^g} in the $k$th round because $(\es{G}{\mcp{h}{\beta^g,\phi^g}}{k}-\es{B}{\mcp{h}{\beta^g,\phi^g}}{k})\cup\es{R}{\mcp{h}{\beta^g,\phi^g}}{k}=(\es{G}{\mcp{h}{\beta^1,\phi^1}}{k}-\es{B}{\mcp{h}{\beta^1,\phi^1}}{k})\cup\es{R}{\mcp{h}{\beta^1,\phi^1}}{k}$.

\noindent$\boldsymbol{\es{F}{\mcp{h}{\beta^g,\phi^g}}{k}:}$ $\phi^g$ has Fixer respond greedily with some \es{F}{\mcp{h}{\beta^g,\phi^g}}{k}.

\noindent$\boldsymbol{\es{F}{\mcp{h}{\beta^1,\phi^1}}{k}:}$ $\phi^1$ has Fixer play $\es{F}{\mcp{h}{\beta^1,\phi^1}}{k}=\es{F}{\mcp{h}{\beta^g,\phi^g}}{k}$, which is possible because $\es{R}{\mcp{h}{\beta^1,\phi^1}}{k}=\es{R}{\mcp{h}{\beta^g,\phi^g}}{k}$.

\noindent\textbf{Result:} \pair{\mcp{h}{\beta^1,\phi^1}}{k+1}{\mcp{h}{\beta^g,\phi^g}}{k+1} belongs to Scenario \ref{sce3} because
\begin{itemize}
\item $\es{G}{\mcp{h}{\beta^1,\phi^1}}{k+1}=(\es{G}{\mcp{h}{\beta^1,\phi^1}}{k}-\es{B}{\mcp{h}{\beta^1,\phi^1}}{k})\cup\es{F}{\mcp{h}{\beta^1,\phi^1}}{k}=(\es{G}{\mcp{h}{\beta^g,\phi^g}}{k}-\es{B}{\mcp{h}{\beta^g,\phi^g}}{k})\cup\es{F}{\mcp{h}{\beta^g,\phi^g}}{k}=\es{G}{\mcp{h}{\beta^g,\phi^g}}{k+1}$
\item $\es{R}{\mcp{h}{\beta^1,\phi^1}}{k+1}=\es{R}{\mcp{h}{\beta^1,\phi^1}}{k}-\es{F}{\mcp{h}{\beta^1,\phi^1}}{k}=\es{R}{\mcp{h}{\beta^g,\phi^g}}{k}-\es{F}{\mcp{h}{\beta^g,\phi^g}}{k}=\es{R}{\mcp{h}{\beta^g,\phi^g}}{k+1}$
\end{itemize}

\noindent\hrulefill

This completes the proof.
\end{proof}

\pagebreak

\section{Details of Proof of Lemma \ref{equivalenceLem}}
\label{c3propAppendix}

Recall the following definitions. If $\phi^*$ and $\phi^{\#}$ are Fixer strategies for which there exists a set $F$ satisfying $F\subseteq\phi^*(\mc{h})\cap\phi^{\#}(\mc{h})$, then $\phi^*$ and $\phi^{\#}$ agree on $F$ at $\mc{h}\in\his{B}{1}$. If $\phi^{\#}$ is a Fixer strategy and $B\subset\es{B}{\mc{h}}{1}$ and $F\subseteq\phi^{\#}(\mc{h})$ are sets of edges such that $(\es{G}{\mc{h}}{1}-B)\cup F$ is connected, then the $B,F$-translations of $\phi^{\#}$ and \mc{h} are the Fixer strategy $\widetilde{\phi^{\#}}$ and history $\mc{\widetilde{h}}\in\his{B}{1}$ with limit $\widetilde{\ell}=\ell-|B|$ (where $\ell$ is the limit of \mc{h}) for which $\es{G}{\mc{\widetilde{h}}}{1}=(\es{G}{\mc{h}}{1}-B)\cup F$, $\es{R}{\mc{\widetilde{h}}}{1}=\es{R}{\mc{h}}{1}-F$, $\es{B}{\mc{\widetilde{h}}}{1}=\es{B}{\mc{h}}{1}-B$, $\widetilde{\phi^{\#}}(\mc{\widetilde{h}})=\phi^{\#}(\mc{h})-F$ (so that $\es{G}{\mcp{\widetilde{h}}{\widetilde{\phi^{\#}}}}{2}=\es{G}{\mcp{h}{\phi^{\#}}}{2}$ and $\es{R}{\mcp{\widetilde{h}}{\widetilde{\phi^{\#}}}}{2}=\es{R}{\mcp{h}{\phi^{\#}}}{2}$), and $\widetilde{\phi^{\#}}(\mc{\widetilde{h'}})=\phi^{\#}(\mc{h'})$ for any histories \mc{h'} and \mc{\widetilde{h'}} in \his{B}{} such that the first round of \mc{h'} matches that of \mcp{h}{\phi^{\#}}, the first round of \mc{\widetilde{h'}} matches that of \mcp{\widetilde{h}}{\widetilde{\phi^{\#}}}, and subsequent rounds of \mc{h'} and \mc{\widetilde{h'}} are identical to each other. We are to prove that if $\phi^*$ and $\phi^{\#}$ are Fixer strategies that agree on $F$ at \mc{h}, and the $B,F$-translations $\widetilde{\phi^{\#}}$ and \mc{\widetilde{h}} of $\phi^{\#}$ and \mc{h} are such that $\widetilde{\phi^{\#}}$ is Fixer-dominant at \mc{\widetilde{h}}, then $\phi^{\#}$ Fixer-dominates $\phi^*$ at \mc{h}.

\begin{proof}
Let $\beta^{\#}$ be an arbitrary Buster strategy; to show $\phi^{\#}$ Fixer-dominates $\phi^*$ at \mc{h}, we need to show there exists a Buster strategy $\beta^*$ such that \mcp{h}{\beta^{\#},\phi^{\#}} is Fixer-superior to \mcp{h}{\beta^*,\phi^*}. Let the $B,F$-translations of $\phi^*$ and \mc{h} be the Fixer strategy $\widetilde{\phi^*}$ and history $\mc{\widetilde{h}}\in\his{B}{1}$, noting that \mc{\widetilde{h}} is the same history from the $B,F$-translations of $\phi^{\#}$ and \mc{h} and has limit $\widetilde{\ell}=\ell-|B|$ (where $\ell$ is the limit of \mc{h}) and satisfies $\es{G}{\mc{\widetilde{h}}}{1}=(\es{G}{\mc{h}}{1}-B)\cup F$, $\es{R}{\mc{\widetilde{h}}}{1}=\es{R}{\mc{h}}{1}-F$, and $\es{B}{\mc{\widetilde{h}}}{1}=\es{B}{\mc{h}}{1}-B$.

Let $\widetilde{\beta^{\#}}$ be a Buster strategy defined to match $\beta^{\#}$ past the first round (i.e. $\beta^{\#}(\mc{h^{\#}})=\widetilde{\beta^{\#}}(\mc{\widetilde{h^{\#}}})$ for any two histories $\mc{h^{\#}}\in\his{F}{}$ and $\mc{\widetilde{h^{\#}}}\in\his{F}{}$ for which \mc{h^{\#}} extends \mcp{h}{\phi^{\#}}, \mc{\widetilde{h^{\#}}} extends \mcp{\widetilde{h}}{\widetilde{\phi^{\#}}}, and \mc{h^{\#}} and \mc{\widetilde{h^{\#}}} only differ in the first round but are identical elsewhere); hence \mcp{h}{\beta^{\#},\phi^{\#}} is identical to \mcp{\widetilde{h}}{\widetilde{\beta^{\#}},\widetilde{\phi^{\#}}} past the first round due to $\widetilde{\phi^{\#}}$ and \mc{\widetilde{h}} being the $B,F$-translations of $\phi^{\#}$ and \mc{h}. Let $k^{\#}=|\mcp{h}{\beta^{\#},\phi^{\#}}|=|\mcp{\widetilde{h}}{\widetilde{\beta^{\#}},\widetilde{\phi^{\#}}}|$ and note that if Fixer wins \mcp{h}{\beta^{\#},\phi^{\#}} then
$$
\pay{\mcp{h}{\beta^{\#},\phi^{\#}}}
=
1+w(\es{R}{\mcp{h}{\beta^{\#},\phi^{\#}}}{k^{\#}})
=
1+w(\es{R}{\mcp{\widetilde{h}}{\widetilde{\beta^{\#}},\widetilde{\phi^{\#}}}}{k^{\#}})
=
\pay{\mcp{\widetilde{h}}{\widetilde{\beta^{\#}},\widetilde{\phi^{\#}}}}
$$
while if Buster wins \mcp{h}{\beta^{\#},\phi^{\#}} then
\begin{align*}
\pay{\mcp{h}{\beta^{\#},\phi^{\#}}}
&=
-1-w(\es{R}{\mcp{h}{\beta^{\#},\phi^{\#}}}{1}-\es{R}{\mcp{h}{\beta^{\#},\phi^{\#}}}{k^{\#}})
=
-1-w((\es{R}{\mcp{\widetilde{h}}{\widetilde{\beta^{\#}},\widetilde{\phi^{\#}}}}{1}\cup F)-\es{R}{\mcp{\widetilde{h}}{\widetilde{\beta^{\#}},\widetilde{\phi^{\#}}}}{k^{\#}}) \\
&=
-1-w(\es{R}{\mcp{\widetilde{h}}{\widetilde{\beta^{\#}},\widetilde{\phi^{\#}}}}{1}-\es{R}{\mcp{\widetilde{h}}{\widetilde{\beta^{\#}},\widetilde{\phi^{\#}}}}{k^{\#}})-w(F)
=
\pay{\mcp{\widetilde{h}}{\widetilde{\beta^{\#}},\widetilde{\phi^{\#}}}}-w(F)
\end{align*}
and no matter what the same player wins both \mcp{h}{\beta^{\#},\phi^{\#}} and \mcp{\widetilde{h}}{\widetilde{\beta^{\#}},\widetilde{\phi^{\#}}}.

Since $\widetilde{\phi^{\#}}$ is Fixer-dominant at \mc{\widetilde{h}} and therefore Fixer-dominates $\widetilde{\phi^*}$ at \mc{\widetilde{h}}, there exists a Buster strategy $\widetilde{\beta^*}$ such that \mcp{\widetilde{h}}{\widetilde{\beta^{\#}},\widetilde{\phi^{\#}}} is Fixer-superior to \mcp{\widetilde{h}}{\widetilde{\beta^*},\widetilde{\phi^*}}; note that this means $\pay{\mcp{\widetilde{h}}{\widetilde{\beta^{\#}},\widetilde{\phi^{\#}}}} \geq \pay{\mcp{\widetilde{h}}{\widetilde{\beta^*},\widetilde{\phi^*}}}$ and thus if Buster wins \mcp{\widetilde{h}}{\widetilde{\beta^{\#}},\widetilde{\phi^{\#}}} then Buster also wins \mcp{\widetilde{h}}{\widetilde{\beta^*},\widetilde{\phi^*}}.

Let $\beta^*$ be a Buster strategy defined to match $\widetilde{\beta^*}$ past the first round, i.e. $\beta^{*}(\mc{h^{*}})=\widetilde{\beta^*}(\mc{\widetilde{h^*}})$ for any two histories $\mc{h^{*}}\in\his{F}{}$ and $\mc{\widetilde{h^*}}\in\his{F}{}$ for which \mc{h^{*}} extends \mcp{h}{\phi^{*}}, \mc{\widetilde{h^*}} extends \mcp{\widetilde{h}}{\widetilde{\phi^{*}}}, and \mc{h^{*}} and \mc{\widetilde{h^*}} only differ in the first round but are identical elsewhere; hence \mcp{h}{\beta^{*},\phi^{*}} is identical to \mcp{\widetilde{h}}{\widetilde{\beta^*},\widetilde{\phi^{*}}} past the first round due to $\widetilde{\phi^{*}}$ and \mc{\widetilde{h}} being the $B,F$-translations of $\phi^{*}$ and \mc{h}. Let $k^{*}=|\mcp{h}{\beta^{*},\phi^{*}}|=|\mcp{\widetilde{h}}{\widetilde{\beta^*},\widetilde{\phi^{*}}}|$ and note that if Fixer wins \mcp{\widetilde{h}}{\widetilde{\beta^*},\widetilde{\phi^*}} then
$$
\pay{\mcp{\widetilde{h}}{\widetilde{\beta^*},\widetilde{\phi^*}}}
=
1+w(\es{R}{\mcp{\widetilde{h}}{\widetilde{\beta^*},\widetilde{\phi^{*}}}}{k^{*}})
=
1+w(\es{R}{\mcp{h}{\beta^*,\phi^*}}{k^*})
=
\pay{\mcp{h}{\beta^*,\phi^*}}
$$
while if Buster wins \mcp{\widetilde{h}}{\widetilde{\beta^*},\widetilde{\phi^*}} then
\begin{align*}
\pay{\mcp{\widetilde{h}}{\widetilde{\beta^*},\widetilde{\phi^*}}}
&=
-1-w(\es{R}{\mcp{\widetilde{h}}{\widetilde{\beta^*},\widetilde{\phi^{*}}}}{1}-\es{R}{\mcp{\widetilde{h}}{\widetilde{\beta^*},\widetilde{\phi^{*}}}}{k^*})
=
-1-w((\es{R}{\mcp{h}{\beta^{*},\phi^{*}}}{1}-F)-\es{R}{\mcp{h}{\beta^{*},\phi^{\*}}}{k^*}) \\
&=
-1-w(\es{R}{\mcp{h}{\beta^{*},\phi^{*}}}{1}-\es{R}{\mcp{h}{\beta^{*},\phi^{\*}}}{k^*})+w(F)
=
\pay{\mcp{h}{\beta^*,\phi^*}}+w(F)
\end{align*}
and no matter what the same player wins both \mcp{\widetilde{h}}{\widetilde{\beta^*},\widetilde{\phi^*}} and \mcp{h}{\beta^*,\phi^*}.

Hence if Fixer wins \mcp{h}{\beta^{\#},\phi^{\#}} and \mcp{h}{\beta^*,\phi^*} then Fixer also wins \mcp{\widetilde{h}}{\widetilde{\beta^{\#}},\widetilde{\phi^{\#}}} and \mcp{\widetilde{h}}{\widetilde{\beta^*},\widetilde{\phi^*}} with $\pay{\mcp{h}{\beta^{\#},\phi^{\#}}}=\pay{\mcp{\widetilde{h}}{\widetilde{\beta^{\#}},\widetilde{\phi^{\#}}}}\geq\pay{\mcp{\widetilde{h}}{\widetilde{\beta^*},\widetilde{\phi^*}}}=\pay{\mcp{h}{\beta^*,\phi^*}}$, if Fixer wins \mcp{h}{\beta^{\#},\phi^{\#}} but Buster wins \mcp{h}{\beta^*,\phi^*} then $\pay{\mcp{h}{\beta^{\#},\phi^{\#}}}>0>\pay{\mcp{h}{\beta^*,\phi^*}}$, and if Buster wins \mcp{h}{\beta^{\#},\phi^{\#}} then Buster also wins \mcp{\widetilde{h}}{\widetilde{\beta^{\#}},\widetilde{\phi^{\#}}}, \mcp{\widetilde{h}}{\widetilde{\beta^*},\widetilde{\phi^*}}, and \mcp{h}{\beta^*,\phi^*} with $\pay{\mcp{h}{\beta^{\#},\phi^{\#}}}=\pay{\mcp{\widetilde{h}}{\widetilde{\beta^{\#}},\widetilde{\phi^{\#}}}}-w(F)\geq\pay{\mcp{\widetilde{h}}{\widetilde{\beta^*},\widetilde{\phi^*}}}-w(F)=\pay{\mcp{h}{\beta^*,\phi^*}}$. Thus \mcp{h}{\beta^{\#},\phi^{\#}} is Fixer-superior to \mcp{h}{\beta^*,\phi^*}, implying $\phi^{\#}$ Fixer-dominates $\phi^*$ at \mc{h}.
\end{proof}

\end{document}